\documentclass[10pt]{article}

\usepackage{xcolor}\definecolor{darkgreen}{rgb}{0,0.45,0}
\usepackage[pagebackref,colorlinks,citecolor=darkgreen,linkcolor=darkgreen]{hyperref}
\usepackage{amsmath,amssymb,url,cleveref,stmaryrd,mathrsfs,mathtools,braket,phonetic,paralist}
\usepackage[all]{xy}
\newcommand{\pullback}[1][dr]{\save*!/#1-1.7pc/#1:(-1,1)@^{|-}\restore}
\usepackage{mathpartir}
\usepackage{tikz}
\usepackage[numbers,sort&compress]{natbib}
\usepackage{appendix}  
\usepackage{fancyvrb}\VerbatimFootnotes
\let\setof\Set
\let\jdeq\equiv
\def\oo{\ensuremath{\infty}}
\def\io{\ensuremath{(\oo,1)}}
\def\ty{\;\mathsf{type}}
\def\elt{\;\mathsf{elt}}
\def\m#1{\llbracket#1\rrbracket}
\def\N{\mathbb{N}}
\def\Z{\mathbb{Z}}
\def\Q{\mathbb{Q}}
\def\R{\mathbb{R}}
\def\bR{\mathbf{R}}
\def\Np{\mathbb{N}_+}
\def\inl{\mathsf{inl}}
\def\inr{\mathsf{inr}}
\def\case{\mathsf{case}}
\def\tprod{\textstyle\prod}
\def\tsm{\textstyle\sum}
\def\refl{\mathsf{refl}}
\def\proj#1{\mathsf{pr}_{#1}}
\def\J{\mathsf{J}}
\def\types{\vdash}
\def\U{\mathscr{U}}
\def\true{\mathrm{T}}
\def\false{\mathrm{F}}
\def\C{\mathscr{C}}
\def\T{\mathcal{T}}
\def\bool{\mathbf{2}}
\def\equiv{\mathsf{Equiv}}
\def\coeq{\mathsf{coeq}}
\def\classof#1{\langle #1\rangle}
\def\ceq{\mathsf{ceq}}
\def\cind{\mathsf{cind}}
\def\trans#1#2{{#1}_*(#2)}
\def\ap#1{\mathsf{ap}_{#1}}
\newcommand{\blank}{\mathord{\hspace{1pt}\text{--}\hspace{1pt}}}
\def\succ{\mathsf{succ}}
\def\hocirc{S^1}
\def\topcirc{\mathbb{S}^1}
\def\hosph{S^2}
\def\topsph{\mathbb{S}^2}
\def\topdisc{\mathbb{D}^2}
\def\unit{\mathbf{1}}
\def\base{\mathsf{base}}
\def\lloop{\mathsf{loop}}
\def\diag{\triangle}
\def\Id{\mathsf{Id}}

\def\Ctx{\mathbf{Ctx}}
\def\fT{\mathfrak{T}}

\def\Lang{\mathfrak{Lang}}
\def\CtxT{\mathbf{Ctx}(\fT)}
\def\cCtx{\underline{\mathscr{C}\!\mathit{tx}}}
\def\cCtxT{\cCtx(\fT)}
\def\sCtx{\underline{\mathsf{Ctx}}}
\def\sCtxT{\sCtx(\fT)}

\crefformat{section}{\S#2#1#3}
\Crefformat{section}{Section~#2#1#3}
\crefrangeformat{section}{\S\S#3#1#4--#5#2#6}
\Crefrangeformat{section}{Sections~#3#1#4--#5#2#6}
\crefmultiformat{section}{\S\S#2#1#3}{ and~#2#1#3}{, #2#1#3}{ and~#2#1#3}
\Crefmultiformat{section}{Sections~#2#1#3}{ and~#2#1#3}{, #2#1#3}{ and~#2#1#3}
\crefrangemultiformat{section}{\S\S#3#1#4--#5#2#6}{ and~#3#1#4--#5#2#6}{, #3#1#4--#5#2#6}{ and~#3#1#4--#5#2#6}
\Crefrangemultiformat{section}{Sections~#3#1#4--#5#2#6}{ and~#3#1#4--#5#2#6}{, #3#1#4--#5#2#6}{ and~#3#1#4--#5#2#6}

\crefname{figure}{Figure}{Figures}
\crefname{table}{Table}{Tables}

\numberwithin{equation}{section}

\newcommand{\trunc}[2]{\mathopen{}\left\Vert #2\right\Vert_{#1}\mathclose{}}

\newcommand{\tproj}[3][]{\mathopen{}\left|#3\right|_{#2}^{#1}\mathclose{}}

\newcommand{\brck}[1]{\trunc{}{#1}}

\newcommand{\bproj}[1]{\tproj{}{#1}}

\newcommand{\shape}{\ensuremath{\mathord{\raisebox{0.5pt}{\text{\rm\esh}}}}}

\title{Homotopy type theory: the logic of space}
\author{Michael Shulman\thanks{This material is based on research sponsored by The United States Air Force Research Laboratory under agreement number FA9550-15-1-0053.  The U.S.~Government is authorized to reproduce and distribute reprints for Governmental purposes notwithstanding any copyright notation thereon.  The views and conclusions contained herein are those of the author and should not be interpreted as necessarily representing the official policies or endorsements, either expressed or implied, of the United States Air Force Research Laboratory, the U.S.~Government, or Carnegie Mellon University.}}

\begin{document}
\maketitle

\section{Introduction to synthetic spaces}
\label{sec:introduction}

There are so many different notions of ``space'' (e.g.\ topological spaces, manifolds, schemes, stacks, and so on, as discussed in various other chapters of this book) that one might despair of finding any common thread tying them together.
However, one property shared by many notions of space is that they can be ``background'' structure.
For instance, many kinds of algebraic objects, such as groups, rings, lattices, Boolean algebras, etc., often come with ``extra space structure'' that is respected by all their operations.
In the case of groups, we have topological groups, Lie groups, sheaves of groups, \oo-groups, and so on.

For each kind of ``spatial group'', much of the theory of ordinary groups generalizes directly, with the ``extra space structure'' simply ``coming along for the ride''.
Additionally, many naturally arising groups, such as the real and complex numbers, matrix groups, the $p$-adic numbers, profinite groups, loop spaces, and so on, come ``naturally'' with spatial structure, and usually it would be ridiculous to study them without taking that spatial structure into account.
On the other hand, ``ordinary'' groups are the special case of ``spatial groups'' whose spatial structure is trivial (e.g.\ discrete); but certain natural constructions on groups, such as the Pontryagin dual, profinite completion, or delooping, take us out of the discrete world.
Thus, the theory of ``groups with spatial structure'' subsumes, and in a sense ``completes'', the study of ordinary groups.
Similar statements can be made about many other kinds of algebraic structure.

With this in mind, the idea of \emph{synthetic spaces} can be summarized as follows:
if all objects in mathematics come naturally with spatial structure, then it is perverse to insist on defining them first in terms of bare sets, as is the official foundational position of most mathematicians, and only \emph{later} equipping them with spatial structure. 
Instead, we can replace set theory with a different formal system whose basic objects are \emph{spaces}. 
Since spaces admit most of the same constructions that sets do (such as products, disjoint unions, exponential objects, and so on), we can develop mathematics in such a system with very few changes to its outward appearance, but all the desired spatial structure will automatically be present and preserved.
(In fact, as we will see in \cref{sec:why-spaces}, this can even be regarded as an explanation of \emph{why} many objects in mathematics come naturally with spatial structure.)
Moreover, if our formal system is sufficiently general, then its objects will be interpretable as many \emph{different} kinds of space; thus the same theorems about ``groups'' will apply to topological groups, Lie groups, sheaves of groups, \oo-groups, and so on.

A formal system with these properties is \emph{Martin--L\"{o}f dependent type theory}~\cite{martinlof:itt-pred,martinlof:itt}.
Originally conceived as a constructive foundation for mathematics where everything has ``computational'' content, it turns out to also admit ``spatial'' interpretations. 
This connection between constructivity/computability and topology/continuity goes back at least to Brouwer, and was gradually developed by many people.\footnote{Escard\'o~\cite{escardo:syntop-datatypes} cites ``Kleene, Kreisel, Myhill/Shepherdson, Rice/Shapiro, Nerode, Scott, Ershov, Plotkin, Smyth, Abramsky, Vickers, Weihrauch and no doubt many others''.}
It was originally restricted to particular topologies on computational data types, but eventually broadened to the realization that types could be interpreted as almost any kind of space. 
Categorically speaking, each ``kind of space'' forms a \emph{topos} or something like it (a category that shares many properties of the category of sets), and the interpretation proceeds by way of constructing a ``free topos'' from the syntax of type theory (see \cref{sec:type-theory}).

There are other formal systems that can be interpreted in toposes, such as Intuitionistic Higher-Order Logic.
Dependent type theory has some minor advantages of convenience, but more importantly, it has recently been recognized~\cite{aw:htpy-idtype,klv:ssetmodel} to also admit interpretations in \emph{higher} toposes.
More concretely, this means we can also interpret its basic objects as \emph{homotopy spaces}, a.k.a.\ \oo-groupoids.
The resulting collection of new axioms and techniques is known as \emph{homotopy type theory}~\cite{hottbook} or \emph{univalent foundations}~\cite{vv:unimath}.
It includes \emph{synthetic homotopy theory} which studies homotopical objects ``directly'' without the need for topological spaces, simplicial sets, or any other combinatorial gadget.
Like any new perspective on a subject, synthetic homotopy theory suggests new ways to attack problems; it has already led to new proofs of known theorems.

Classically, \oo-groupoids arose to prominence gradually, as repositories for the homotopy-theoretic information contained in a topological space; see Porter's chapter for an extensive survey.
As we will see, however, the synthetic viewpoint emphasizes that this structure of a ``homotopy space'' is essentially \emph{orthogonal} to other kinds of space structure, so that an object can be both ``homotopical'' and (for example) ``topological'' or ``smooth'' in unrelated ways.
This sort of mixed structure is visible in many other chapters of the present volume, such as those about Lie groupoids (Pradines), toposes (Joyal), and stacks (Mestrano--Simpson).
It is also central to many applications, such as differential cohomology and gauge field theory (e.g.\ Schreiber's chapter).
Finally, it describes cleanly how topological and smooth spaces give rise to homotopy ones (see \cref{sec:cohes-homot-type}).

This chapter is intended as a brief introduction to the above ideas: type theory, synthetic spaces, and homotopy type theory.
Of course many details will be left out, but I hope to convey a flavor of the subject, and leave the reader with some idea of what it means to talk about the \emph{logic of space}.

It should be emphasized that homotopy type theory, in particular, is a very new subject.
Many of its basic definitions are still in flux, and some of its expected fundamental theorems have not yet been completely proven.
In general I will focus on what is \emph{expected} to be true, 
in order to emphasize the possibilities opened up by these ideas; but I will endeavor not to lie, and to include some remarks on the current state of the art as well.

I will begin in \cref{sec:type-theory} with an introduction to type theory.
Then in \cref{sec:synthetic-topology,sec:homotopy-type-theory} I will discuss its spatial and homotopical aspects respectively and some of their applications.
Finally, in \cref{sec:cohes-homot-type} I will briefly mention how these aspects are combined. 
For further reading, I recommend~\cite{awodey:tt-and-htpy,apw:vvu-hott,pw:hottvvuf,shulman:synhott} and~\cite{hottbook}.

I would like to thank Peter LeFanu Lumsdaine, Steve Awodey, and Urs Schreiber, as well as the editors of the volume, for careful reading and helpful comments.

\section{Type theory}
\label{sec:type-theory}

\subsection{On syntax}
\label{sec:syntax}

Mathematicians often have a lot of difficulty understanding type theory (and the author was no exception).
One reason is that the usual presentation of type theory is heavy on \emph{syntax}, which most mathematicians are not used to thinking about.
Thus, it is appropriate to begin with a few general remarks about syntax, what it is, and its role in type theory and in mathematics more generally.%
\footnote{The ``algebraic'' perspective I will present is only one of many valid ways to look at type theory.
  It has been developed by~\cite{ls:itt-free-topos,dybjer:internal-tt,streicher:semtt,vv:csys-jfrel} among others.}

In general, \emph{syntax} refers to a system of formal symbols of some sort, whereas \emph{semantics} means the interpretation of those symbols as ``things''.
In the language of category theory, we can generally think of syntax as describing a \emph{free} (or \emph{presented}) object, and semantics as the morphisms out of that object determined by its universal property.

For instance, in group theory we may write a sequence of equations such as
\begin{equation}
  (g h) h^{-1} = g (h h^{-1}) = g e = g.\label{eq:groupeqn}
\end{equation}
Where does this computation take place?
One obvious answer is ``in an arbitrary group''.
But another is ``in the free group $F\langle g,h\rangle$ generated by two symbols $g$ and $h$.''
Since the elements of $F\langle g,h\rangle$ are literally strings of symbols (``words'') produced by multiplication and inversion from $g$ and $h$, strings such as ``$(g h)h^{-1}$'' \emph{are themselves} elements of $F\langle g,h\rangle$, and~\eqref{eq:groupeqn} holds as an equality between these elements, i.e.\ a statement \emph{in syntax}.
Now if we have any other group $G$ and two elements of it, there is a unique group homomorphism from $F\langle g,h\rangle$ to $G$ sending the letters $g$ and $h$ to the chosen elements of $G$.
This is the \emph{semantics} of our syntax, and it carries the equation~\eqref{eq:groupeqn} in $F\langle g,h\rangle$ to the analogous equation in $G$.
Such reasoning can be applied to arguments involving hypotheses, such as ``if $g^2 = e$, then $g^4 = (g^2)^2 = e^2 = e$'', by considering (in this case) the group $F\langle g \mid g^2=e\rangle$ \emph{presented} by one generator $g$ and one equation $g^2=e$.
(A free group, of course, has a presentation with no equations.)

In other words, we can regard an argument such as~\eqref{eq:groupeqn} either as a ``semantic'' statement about ``all groups'' or as a ``syntactic'' statement about a particular free or presented group.
The former is a consequence of the latter, by the universal property of free groups and presentations.

This may seem like mere playing with words,\footnote{No pun intended.} and the reader may wonder how such a viewpoint could ever gain us anything.
The reason is that often, we can say more about a free object than is expressed tautologically by its universal property.%
\footnote{This is dual to the familiar fact that studying a ``classifying space'' can yield insights about the objects it classifies --- a classifying space being a representing object for a contravariant functor, while a free or presented object represents a covariant one.}
Usually, this takes the form of an explicit and tractable construction of an object that is then \emph{proven} to be free.
Of course, \emph{any} construction of a free object must be proven correct, but such a proof can range from tautological to highly nontrivial.
The less trivial it is, the more potential benefit there is from working syntactically with the free object.

For instance, a ``tautological'' way to define $F\langle g,h\rangle$ is by ``throwing in freely'' the group operations of multiplication and inversion, obtaining formal ``words'' such as $(gg^{-1})(h^{-1}(h g))$, and then quotienting by an equivalence relation generated by the axioms of a group.
The universal property of a free group is then essentially immediate.
But a more interesting and useful construction of $F\langle g,h\rangle$ consists of ``reduced words'' in $g$, $h$, and their formal inverses (finite sequences in which no cancellation is possible), such as $g h g^{-1} g^{-1} h h g h^{-1}$, with multiplication by concatenation and cancellation.
The proof that this yields a free group is not entirely trivial (indeed, even the definition of the group multiplication is not completely trivial); but once we know it, it can simplify our lives.

As a fairly banal example of such a simplification, recall that the \emph{conjugation} of $h$ by $g$ is defined by $h^g = g h g^{-1}$.
Here is a proof that conjugation by $g$ is a group homomorphism:
\begin{equation}
  h^g\,k^g = (g h g^{-1})(g k g^{-1}) = g h k g^{-1} = (hk)^g\label{eq:conj}
\end{equation}
As straightforward as it is, this is not, technically, a complete proof from the usual axioms of a group.
For that, we would have to choose parenthesizations and use the associativity and unit axioms explicitly:
\begin{multline}\label{eq:conj-ugly}
  h^g k^g
  = ((g h) g^{-1})((g k) g^{-1})
  = (g (h g^{-1}))((g k) g^{-1})
  = ((g (h g^{-1}))(g k)) g^{-1}\\
  = (g ((h g^{-1})(g k))) g^{-1}
  = (g (h (g^{-1}(g k)))) g^{-1}
  = (g (h ((g^{-1} g) k))) g^{-1}\\
  = (g (h (e k))) g^{-1}
  = (g (h k)) g^{-1}
  = (h k)^g
\end{multline}
Of course, this would be horrific, so no one ever does it.
If mathematicians think about this sort of question at all, they usually call~\eqref{eq:conj} an ``acceptable abuse of notation''.
But with the above explicit description of free groups, we can make formal sense of~\eqref{eq:conj} as a calculation in $F\langle g,h,k\rangle$, wherein $g h g^{-1}$ and $g k g^{-1}$ are specific elements whose product is $g h k g^{-1}$.
Then we can extend this conclusion to every other group by freeness.
Note that if we tried to do the same thing with the ``tautological'' presentation of a free group, we would be forced to write down~\eqref{eq:conj-ugly} instead, so no simplification would result.

In general, there are several ways that a presentation of a free object might make our lives easier.
One is if its elements are ``canonical forms'', as for free groups (e.g.\ $g h k g^{-1}$ is the canonical form of $((g h) g^{-1})((g k) g^{-1})$).
This eliminates (or simplifies) the quotient by an equivalence relation required for ``tautological'' constructions.
Often there is a ``reduction'' algorithm to compute canonical forms, making equality in the free object computationally decidable.

Another potential advantage is if we obtain a ``version'' of a free object that is \emph{actually} simpler.
For instance, it might be stricter than the one given by a tautological construction.
This is particularly common in category theory and higher category theory, where it can be called a \emph{coherence theorem}.

Finally, a particular construction of a free object might also be \emph{psychologically} easier to work with, or at least suggest a different viewpoint that may lead to new insights.
The best example of this is type theory itself: though it also offers the advantages of canonical forms and strictness (see \cref{sec:compute,sec:semantics}), arguably its most important benefit is a way of thinking.

\subsection{Universes of mathematics}
\label{sec:univ-math}

What, then, \emph{is} type theory?\footnote{\label{fn:typetheory}Unfortunately, the phrase ``type theory'' has many different meanings.
On the one hand, type theory is a \emph{discipline} lying at the boundary of mathematics and computer science.
This discipline studies deductive systems that are themselves also known as type theories.
But in the context of mathematical foundations, such as here, ``type theory'' generally refers to a particular subclass of these deductive systems, which are more precisely called \emph{dependent type theories} (because they admit ``dependent types''; see below.)}
Roughly speaking, it is a particularly convenient construction of free objects for the theory of \emph{all of mathematics}.
Just as a group presented by $g$, $h$ and $g h = h g$ admits a unique homomorphism to any other group equipped with two commuting elements, type theory with certain structures presents ``a universe of mathematics'' with a unique ``mathematics-homomorphism'' to any other such universe of mathematics.

This discussion of ``universes of mathematics'' may sound odd; surely there is only one universe of mathematics?
Well, yes, mathematics is a whole; but it has been known since the early 20th century that some formal systems, such as Zermelo--Fraenkel set theory, can encode almost all of mathematics. 
To first approximation, by a ``universe of mathematics'' I mean a model of a formal system in which mathematics can be encoded.
Note that G\"odel's incompleteness theorem ensures that any such system has many different models.%
\footnote{Specifically, it shows that any sufficiently powerful formal system contains statements that are neither provable nor disprovable.
The \emph{completeness} theorem then implies that there must be some models in which these statements are true and some in which they are false.}
Thus, there are many ``universes of mathematics'' in this sense.

Often the incompleteness theorem
is seen as a bug, but from our point of view
it is actually a feature!
We can make positive use of it by recognizing that certain mathematical structures, like notions of space, 
happen to form new universes of mathematics by themselves.
In other words, starting from one universe of mathematics,%
\footnote{We may regard the starting universe as the ``true'' one, but there is no formal justification for this.
  We will come back to this in \cref{sec:conclusion}.}
we can construct another universe whose objects are, from the point of view of our original universe, ``spaces'' of some sort.
Thus, when a mathematician living in this new world constructs a bare function $A\to B$ between sets, the mathematician in the old world sees that it is in fact a \emph{continuous} function between \emph{spaces}.

This is admittedly a bit vague, so
let me pass to a second approximation of what I mean by a ``universe of mathematics'': a \emph{category}, or \io-category, with certain structure.
Our starting universe is then the category of sets (or perhaps the \io-category of \oo-groupoids).
Thus, type theory gives a way to construct free or presented objects in some \emph{category of structured categories}.
Such a free object is sometimes called the ``syntactic category'' or ``classifying category'' of the type theory. 
In the words of Scott~\cite{scott:ttalt}:
\begin{quote}
  \dots a category represents the ``algebra of types'', just as abstract rings give us the algebra of polynomials, originally understood to concern only integers or rationals.
\end{quote}

Now the usual way of working ``inside'' a particular category is to write all arguments in diagrammatic language.
For instance, if $G$ is a group object in a category (such as a topological group in the category of topological spaces, or a Lie group in the category of smooth manifolds), then the argument analogous to~\eqref{eq:groupeqn} would be the commutativity of the following diagram:
\begin{equation}
  \vcenter{\xymatrix@C=3.8pc{
      G\times G\ar[r]^-{1\times \Delta}\ar[d]_{\mathrm{proj}_1} &
      G\times G\times G\ar[r]^-{1\times 1\times \mathrm{inv}} &
      G\times G\times G\ar[r]^-{\mathrm{mult}\times 1}\ar[d]_{1\times \mathrm{mult}} &
      G\times G\ar[d]^-{\mathrm{mult}}\\
      G\ar[rr]^-{1\times \mathrm{id}} \ar@(dr,dl)[rrr]_{1} &&
      G\times G\ar[r]^-{\mathrm{mult}} &
      G
      }}\label{eq:groupdiag}
\end{equation}
Categorically trained mathematicians become quite adept at translating calculations like~\eqref{eq:groupeqn} into diagrams like~\eqref{eq:groupdiag}.
However, objectively I think it is hard to deny the relative simplicity of~\eqref{eq:groupeqn} compared to~\eqref{eq:groupdiag}.
The benefits are magnified further when we include additional simplifications like those in~\eqref{eq:conj}.

Type theory allows us to use equations like~\eqref{eq:groupeqn} and~\eqref{eq:conj} to prove things about all group objects in \emph{all} categories.
Its syntax involves elements with operations and equations, so we can speak and think as if we were talking about ordinary sets.
But it is nevertheless a description of a free category of a certain sort,\footnote{See \cref{sec:rules,sec:semantics} for some caveats to this statement.} so that its theorems can be uniquely mapped into any other similar category.
Thus, type theory supplies a different perspective on categories that is often more familiar and easier to work with.

To be a little more precise, the benefit here comes from the interplay between two modes of interacting with type theory.
On one hand, we can define and study the formal system of type theory \emph{inside} mathematics.
This enables us to talk about its having multiple models, and hence functioning as a syntax for categories, as described above.
But on the other hand, because type theory is sufficiently powerful to encode all of mathematics, we are also free to regard it as the ``ambient foundation'' for any mathematical theory.
Most modern mathematicians implicitly assume set theory as a foundation, but for the most part type theory is just as good (and, as we will see in \cref{sec:homotopy-type-theory}, it makes ``new kinds of mathematics'' possible as well).
Of course, real-world mathematics is rarely ``fully encoded'' into \emph{any} foundational system, but experience shows that it is always possible in principle, and nowadays with computer proof assistants it is becoming more common and feasible to do explicitly.

The point, then, is that any theorem in ``ordinary'' mathematics can be encoded using the second ``foundational'' point of view, obtaining a derivation in the formal system of type theory; but then we can switch to the first ``semantic'' point of view and conclude that \emph{that theorem} is actually true (suitably interpreted) in all categories with appropriate structure.
In this way, \emph{any} mathematical theorem is actually much more general than it appears.\footnote{However, as we will see in \cref{sec:constructive-logic}, it requires some care on the side of ordinary mathematics --- specifically, avoiding certain restrictive logical axioms --- to maximize this resulting generality.}

\subsection{Types versus sets}
\label{sec:types-versus-sets}

With those lengthy preliminaries out of the way, let's move on to what type theory actually looks like.
If it is to describe the free ``universe of mathematics'', type theory should be a formal system into which mathematics can be encoded.
The currently accepted formal system for encoding mathematics is Zermelo--Fraenkel set theory (ZFC), and mathematicians have a great deal of practice representing structures as sets.
Thus it makes sense that the basic objects of type theory, called \emph{types}, are very set-like --- with one important difference.


In ZFC, an assertion of membership like ``$x\in A$'' is a statement about two previously given objects $x$ and $A$, which might be true or false, and can be hypothesized, proven, or disproven.
In other words, the universe of ZFC is a vast undifferentiated collection of things called ``sets'', with a relation called ``membership'' that can be applied to any two of them.
By contrast, in type theory, the type to which an element belongs is ``part of its nature'', rather than something we can ask about and prove or disprove; two distinct types can never\footnote{As with almost any general statement about type theory, there are exceptions to this, but for the most part it is true.} share any elements in common.
To emphasize this difference, we write $x:A$, rather than $x\in A$, to mean that $x$ is an element of the type $A$.

This perspective on sets is like that of categorical or ``structural'' set theory, such as Lawvere's ETCS~\cite{lawvere:etcs-long,leinster:etcs}, which axiomatizes the \emph{category} of sets and functions.
It contrasts with membership-based or ``material'' set theory such as ZFC, which axiomatizes the \emph{class} of sets and its membership relation.
The structural approach generalizes better when thinking of the basic objects as spaces rather than bare sets, since the spatial relationships between points are specified by an ambient space: it doesn't make sense to ask whether two points are ``nearby'' unless we have fixed some space in which they both reside.

In principle it may be possible to use a more ZFC-like formal system for at least some of the same purposes as type theory (see e.g.~\cite{abss:long-version}), but the connection to spaces would become rather more tenuous.
Moreover, the structural perspective matches the usage of ``sets'' in most of mathematics.
Outside the formal theory of ZFC, the primary place where one element can belong to more than one set, or where elements of distinct sets are compared, is when the given sets are \emph{subsets} of some ambient set.
This situation is encoded in type theory by a notion of ``subset of $A$'' that, like ``element of $A$'', is a basic notion not reducible to something like ``set that happens to be a subset of $A$''; see \cref{sec:prop-as-types}.

While we are talking about ZFC and set theory, it is worth mentioning another reason type theory is often difficult for mathematicians.
Any formal system for encoding mathematics, be it ZFC, ETCS, or type theory, must by its nature be careful about many things that mathematicians usually gloss over.
Ordinary mathematical notation and writing is, technically speaking, ambiguous and full of gaps, trusting the human reader to draw ``obvious'' conclusions.
But to give a mathematical \emph{theory} of mathematics (and in particular, to prove things like ``type theory presents a free structured category''), we have to remove all such ambiguity and fill in all the gaps.
This causes the syntactic formulas of the formal system to appear quite verbose, and often barely comprehensible to a mathematician accustomed to informal mathematical language.

The important points are that this is true for \emph{all} formal systems, and that it should not bother us when doing ordinary mathematics.
The process of ``encoding'' mathematics into a formal system such as ZFC, ETCS, or type theory looks somewhat different depending on which formal system is chosen, but it is generally well understood.
In particular, no matter what formal system we choose, there is no need for its verbosity to infect ordinary mathematics; we remain free to ``abuse notation'' in the usual way.

I stress this point because one sometimes encounters a false impression that type theory requires ``heavier syntax'' than set-based mathematics, or that it forbids ``abuse of notation''.
This is probably partly because type theory is often presented in a very formal and syntactic way --- perhaps because many type theorists are logicians or computer scientists --- whereas most mathematicians' exposure to set theory has been fairly informal and intuitive.
Moreover, the particular notations used in type theory are somewhat unfamiliar to mathematicians, and take some practice to learn to read correctly.
But the syntax of type theory is \emph{intrinsically} no heavier or unabusable than that of set theory.
(Promoting a style of informal mathematics that matches the formal system of type theory was one of the explicit goals of~\cite{hottbook}.)

\subsection{Judgments and the classifying category}
\label{sec:type-term-judgments}

Finally, we are ready to describe the syntax of type theory and how it generates a category (which we will call the \emph{classifying category}; it is also called the \emph{syntactic category} and the \emph{category of contexts}).
Like the elements of a free group, the syntactic objects of type theory are ``words'' built out of operations.
In a free group there is only one sort of word, since a group involves only one collection of ``things'' (its elements).
But since type theory presents a category with both objects and morphisms, 
it has at least two sorts of ``words''.
Type theorists call a ``sort of word'' a \emph{judgment form}, and a particular word a \emph{judgment}.

The first judgment form is a \emph{type judgment}; it is written ``$B\ty$'' and pronounced ``$B$ is a type''.
Here $B$ is a syntactic expression like $\N\times (\R+\Q)$, 
in which $\times$ and $+$ are operations on types, formally analogous to the multiplication of elements represented by concatenation of words in a free group; we will come back to them in \cref{sec:rules}.
The objects of the ``classifying category'' generated by a type theory are\footnote{Well, not exactly; see below.} the syntactic expressions $B$ for which the judgment $B\ty$ can be produced by the rules (i.e.\ operations) to be described in \cref{sec:rules}.
For clarity, we will write $\m B$ when $B$ is regarded as an object of this category, and say that $B$ \emph{presents} the object $\m B$.

The second judgment form is a \emph{term judgment}, written ``$b:B$''.
Here $B$ is a syntactic expression for a type (i.e.\ we must also have ``$B\ty$'').
For instance, we might have $(3\cdot 2+1, \mathsf{inr}(\frac{3}{4}-17)) : \N\times (\R+\Q)$.
Here again, $\cdot,+,-,\mathsf{inr}$ and so on denote operations that will be described in \cref{sec:rules}.
We pronounce $b:B$ as ``$b$ is an element of $B$'' or ``$b$ is a point of $B$'' or ``$b$ is a term of type $B$'', emphasizing respectively the set-like, space-like, or syntactic character of $B$.

More generally, a term judgment can include a \emph{context}, consisting of a list of variables, each with a specified type, that may occur in the term $b$.
For instance, we might also write $(3x+1, \mathsf{inr}(\frac{3}{4}-y)) : \N\times (\R+\Q)$, which only makes sense in the context of $x:\N$ and $y:\Q$.
The traditional notation in type theory is to write the context as a list of variables with their types, joined by commas, and separate it from the judgment with the symbol $\types$ (called a turnstile).
Thus, the above judgment would be written
\[ x:\N, y:\Q \types (3x+1, \mathsf{inr}(\textstyle\frac{3}{4}-y)) : \N\times (\R+\Q) \]
Here the $\types$ is the ``outer relation'' that binds most loosely; then the commas on the left-hand side bind next most loosely, separating the (variable : type) pairs.
Thus, for emphasis it could be bracketed as
\[ ((x:\N), (y:\Q)) \types ((3x+1, \mathsf{inr}(\textstyle\frac{3}{4}-y)) : \N\times (\R+\Q)). \]
Often the Greek letter $\Gamma$ denotes an arbitrary context, so that $\Gamma, x:A \types b:B$ (to be parsed as $(\Gamma, (x:A)) \types (b:B)$) means that in some arbitrary context together with a variable $x$ of type $A$, we have a term $b$ of type $B$.%
\footnote{Technically, $\Gamma$, $b$, $A$, and $B$ here are ``meta-variables'', not to be confused with the ``variables'' such as $x$ in a context.  We will come back to this in \cref{sec:rules}.}

Term judgments $\Gamma \types a:A$ present \emph{morphisms} in the classifying category.
In the simplest case, $\Gamma$ contains only one variable, such as $x:A \types b:B$, and this morphism $\m b$ is from $\m A$ to $\m B$.
For the general case $\Gamma \types a:A$, we have to modify our definition of the classifying category by taking its objects to be \emph{contexts} rather than types, with our previous $\m A$ corresponding to $\m{x:A}$; then 
$\Gamma \types a:A$ presents a morphism from $\m \Gamma$ to $\m{A}$.
For this reason, the classifying category is also known as the \emph{category of contexts}; we denote it by $\Ctx$.

We stipulate that $\m\Gamma$ is a product of the types in $\Gamma$, so $\m{x:A,y:B}\cong \m{A}\times \m B$ and so on.
(In particular, the empty context yields a terminal object $\m{~}$.)
Thus, for instance, $x:\R, y:\R \types x y:\R$ yields the multiplication map $\R\times \R \to \R$.
The universal property of products implies that for contexts $\Gamma$ and $\Delta$, a morphism in $\Ctx$ from $\m\Gamma$ to a general context $\m\Delta$ must consist of a tuple of term judgments $\Gamma \types b_i : B_i$ for all variables $y_i:B_i$ occurring in $\Delta$.
If we also have $\Delta \types c:C$, we get another term judgment denoted
\[\Gamma \types c[b_1/y_1,\dots,b_m/y_m]:C\]
by \emph{substituting} each $b_i$ for $y_i$ in $c$; this presents the composite $\m \Gamma \to \m \Delta \to \m C$.
For instance, we have a morphism from $\m{x:\R}$ to $\m{z:\R,w:\R}$ defined by the terms $x:\R \types (x-1):\R$ and $x:\R \types (x+1):\R$; substituting it into $z:\R, w:\R \types z w:\R$ gives $x:\R \types (x-1)(x+1) : \R$.
That is,
\[(z w)[(x-1)/z,(x+1)/w] = (x-1)(x+1).\]

So far, we have described \emph{simple} type theory.
Next we allow type judgments ``$B\ty$'' to have a context as well, making $B$ into a \emph{dependent type} or \emph{type family}.
Intuitively, a dependent type ``$\Gamma \types B\ty$'' presents an object of the \emph{slice category} $\Ctx/\m\Gamma$, i.e.\ an object with a morphism to $\m\Gamma$.
We think of the ``fiber'' over a point of $\m\Gamma$ as the instance of $B$ corresponding to that point.

For instance, in informal mathematics we might speak of ``an arbitrary finite cyclic group $C_n$'', for $n: \Np$.
In type theory this becomes $n:\Np \types C_n \ty$, corresponding categorically to $\coprod_{n\in \Np} C_n$ with its projection to $\Np$.
Topologically, this is a \emph{bundle} over $\Np$, with the $C_n$ as its fibers.
Working in a category in the style of~\eqref{eq:groupdiag} requires manually translating from ``arbitrary objects'' to bundles; this is one of the least convenient aspects of categorical set theories and of the traditional way of doing mathematics in a topos.
The ability to talk directly about families of types and have them interpreted automatically as bundles is one of the most significant advantages of type theory.

A crucial fact is that \emph{substitution into a dependent type} presents the \emph{pullback} functor between slice categories.
For instance, we have a judgment $\types 3:\Np$ with no variables, yielding a context morphism from the terminal object $\m{~}$ to $\m{\Np}$.
Substitution into $n:\Np \types C_n \ty$ yields the non-dependent type $\types C_3\ty$, which is the pullback of $\coprod_{n\in \Np} C_n$ along the inclusion $3 :\m{~}\to \m{\Np}$:
\[ \xymatrix@-.7pc{ C_3 \ar[r] \ar[d] \pullback & \textstyle\coprod_{n\in \Np} C_n \ar[d] \\ 1 \ar[r] & \Np } \]
As an even simpler example, if $\types B\ty$ is a non-dependent type, we can substitute it along the unique context morphism from any $\Gamma$ to the empty context, yielding a ``trivially-dependent type'' $\Gamma \types B\ty$.
This presents the pullback of $\m B$ to the slice over $\m\Gamma$, i.e.\ the projection $\m\Gamma\times\m B \to \m\Gamma$ (a ``trivial bundle''):
\[ \xymatrix@-.7pc{ \m\Gamma\times\m B \ar[r] \ar[d] \pullback & \m B \ar[d] \\ \m \Gamma \ar[r] & \m{~} \mathrlap{\;=1} } \]

With dependent types, we can allow the type $B$ in a term judgment $\Gamma \types b:B$ to also depend on $\Gamma$.
For instance, the generators of the cyclic groups form a term judgment $n:\Np \types g_n : C_n$.
Such a judgment $\Gamma \types b:B$ represents a \emph{section} of the projection represented by the dependent type $\Gamma \types B\ty$: we ``select one point in each fiber''.
This includes the non-dependent case because morphisms $\m\Gamma\to \m B$ are equivalent to sections of the projection $\m\Gamma\times\m B \to \m\Gamma$.

\label{sec:identity-types}
An example that will be central to the story of this entire chapter is the \emph{diagonal} map $\Delta_{\m A} : \m A \to \m A \times \m A$.
We can regard this as an object of the slice category $\Ctx/(\m A\times \m A)$, or equivalently $\Ctx/\m{x:A,y:A}$; it is then presented by a dependent type called the \emph{equality type} or \emph{identity type}, written
\[x:A,y:A\types (x=y)\ty \qquad\text{or}\qquad x:A,y:A\types \Id(x,y)\ty.\]
We will explain this type in more detail in \cref{sec:myst-ident-types}.
For the moment, we observe that it reduces equalities of terms to existence of terms.
For instance, given $\Gamma \types a:A$ and $\Gamma \types b:A$ representing morphisms $\m a,\m b : \m\Gamma\to \m A$, substituting them into the equality type we get a dependent type $\Gamma \types (a=b)\ty$ that presents the pullback of $\Delta_{\m A}$ along $(\m a, \m b) : \m\Gamma \to \m A \times \m A$, or equivalently the \emph{equalizer} of $\m a$ and $\m b$.
Thus, a judgment $\Gamma \types e : a=b$ says that this equalizer has a section, or equivalently that $\m a=\m b$.
So our type and term judgments also suffice to present \emph{equality} of morphisms.

To be precise, in the presence of dependent types we extend our previous definition of the classifying category $\Ctx$ as follows.
Firstly, we also allow the types in a context to depend on the variables occurring earlier in the same context.
For instance, we can form the context $(n:\Np, x:C_n)$, and then in this context write $n:\Np, x:C_n \types x^2 : C_n$ for the operation that squares an arbitrary element of an arbitrary cyclic group.
Categorically, if $\Gamma\types B\ty$ presents an object of the slice over $\m\Gamma$, i.e.\ a morphism with codomain $\m\Gamma$, then the extended context $\m{\Gamma,x:B}$ is the \emph{domain} of this morphism.
This reduces to our previous $\m{\Gamma,x:B} = \m\Gamma\times \m B$ if $B$ is non-dependent.

Secondly, we take the objects of $\Ctx$ to be contexts in this generalized sense, and a morphism from $\m{\Gamma}$ to $\m{\Delta}$ to consist of term judgments for all $1\le i\le m$:
\[\Gamma \types b_i : B_i[b_1/y_1,\dots,b_{i-1}/y_{i-1}]\]
where $\Delta = (y_1:B_1, y_2:B_2,\dots, y_m:B_m)$, with $y_j$ potentially occurring in $B_i$ for $j<i$.
That is, we first give $\Gamma\types b_1:B_1$, presenting a morphism
\[\m{b_1} : \m\Gamma\to\m{B_1}. \]
Then we substitute $b_1$ for $y_1$ in $B_2$, obtaining a type $\Gamma\types B_2[b_1/y_1]$ and a corresponding extended context that presents the pullback
\[ \xymatrix{ \m{\Gamma, y_2:B_2[b_1/y_1]} \ar[r] \ar[d] \pullback & \m{y_1:B_1, y_2:B_2} \ar[d]\\
  \m\Gamma \ar[r]_{\m{b_1}} & \m{B_1}. } \]
Next we give $\Gamma \types b_2 : B_2[b_1/y_1]$, which presents a section of this pullback, or equivalently a morphism $\Gamma \to \m{y_1:B_1, y_2:B_2}$ making this triangle commute:
\[ \xymatrix@-.7pc{
& \m{y_0:B_0, y_1:B_1} \ar[d] \\
\m{\Gamma} \ar[r]_-{\m{b_0}} \ar[ur]^-{\m{b_1}} 
& \m{y_0:B_0}. } \]
Continuing in this way, the sequence of terms $(b_0,b_1,\dots,b_m)$ that represent a morphism $\m \Gamma \to \m\Delta$ individually present a tower of sections
\[ \xymatrix@-.7pc{
& \m{y_0:B_0, y_1:B_1,\dots, y_m:B_m} \mathrlap{\;=\m\Delta} \ar[d] \\
& \vdots \ar[d] \\
\m{\Gamma} \ar[dr]_-{\m{b_0}} \ar[r]^-{\m{b_1}} \ar[uur]^-{\m{b_m}} & \m{y_0:B_0, y_1:B_1} \ar[d] \\
& \m{y_0:B_0} } \]
with $\m{b_m}$ being the overall morphism $\m\Gamma\to\m\Delta$.
For instance, the ``squaring'' injections $i_n:C_n \hookrightarrow C_{2n}$, represented by term judgments
\[n:\Np \types 2n:\Np
\qquad\text{and}\qquad
n:\Np, x:C_n \types i_n(x) : C_{2n},\]
assemble into a morphism $\m{n:\Np, x:C_n} \to \m{m:\Np, y:C_m}$.
Categorically, this is a morphism $\coprod_n C_n \to \coprod_m C_m$ that sends the $n^{\mathrm{th}}$ summand to the $2n^{\mathrm{th}}$ summand.

Finally, we quotient these morphisms by an equivalence relation arising from the identity type.
In the simplest case where each context has only one type, we identify the morphisms presented by $x:A \types b_1:B$ and $x:A \types b_2:B$ if there is a term $x:A \types p:b_1=b_2$.
The case of morphisms between arbitary contexts is a generalization of this. 
(We will reconsider this last step in \crefrange{sec:semantics}{sec:equality}.)

This completes our definition of the classifying category of a type theory.
We can now \emph{define} the projection morphism $\m{\Gamma, z:C } \to \m{\Gamma}$ associated to a dependent type $\Gamma\types C\ty$, exhibiting $\m{\Gamma,z:C}$ as an object of the slice category over $\m\Gamma$, as we intended.
According to the above description of morphisms, this projection morphism should consist of a term in context $\Gamma,z:C$ for each type in $\Gamma$; we take these to be just the variables in $\Gamma$, ignoring $z$.

For instance, the projection map $\m{x:A, y:B,z:C} \to \m{x:A,y:B}$ is determined by the terms $x:A, y:B, z:C \types x:A$ and $x:A, y:B,z:C \types y:B$.
Similarly, a section of this projection consists of terms
\begin{align*}
  x:A, y:B &\types a:A\\
  x:A, y:B &\types b:B[a/x]\\
  x:A, y:B &\types c:C[a/x,b/y]
\end{align*}
such that the composite $\m{x:A,y:B}\to \m{x:A, y:B,z:C} \to \m{x:A,y:B}$ is the identity, i.e.\ that $a$ and $b$ are the same as $x$ and $y$.
Thus, such a section is simply determined by a term $x:A,y:B\types c:C$, as we intended.

Of course, not \emph{every} object of the slice category $\Ctx/\m{\Gamma}$ is of this form, but every object of $\Ctx/\m{\Gamma}$ is \emph{isomorphic} to one of this form.
Consider the simplest case when $\Gamma$ is a single type $B$, and we have an object of $\Ctx/\m{B}$ whose domain is also a single type $A$, equipped with a term $x:A \types f(x):B$.
Let $\Psi$ denote the context $(y:B, x:A, p:f(x)=y)$;
then $\m\Psi$ is the pullback
\[ \xymatrix@C=3pc{ \m \Psi \ar[r] \ar[d] \pullback & \m B \ar[d]^{\diag_{\m B}} \\
  \m B\times \m A \ar[r]_{1\times \m f} & \m B \times \m B } \]
using the identity type $y_1:B,y_2:B \types (y_1=y_2)\ty$ mentioned above to present the diagonal $\diag_{\m B}$.
It is easy to see categorically that such a pullback is isomorphic to $\m A$.
Thus, every object of $\Ctx/\m B$ is at least isomorphic to a composite of \emph{two} projections from dependent types
\[ \m{y:B,x:A,p:f(x)=y} \to \m{y:B,x:A} \to \m{y:B}. \]
Using the $\Sigma$-type to be defined in \cref{sec:rules}, we can reduce this to one such projection:
\[ \m{y:B, z:\textstyle\sum_{x:A} (f(x)=y)} \to \m{y:B}. \]
A similar argument works with $B$ replaced by any context $\Gamma$.
Thus we can assume that any object of a slice category is determined by a dependent type.




\subsection{Rules and universal properties}
\label{sec:rules}

In the previous section we described the judgment forms of type theory and how they present the classifying category, claiming that each judgment 
is analogous to a word in a free group.
In this section we will describe \emph{what} the judgments are for each judgment form, or more precisely how we can \emph{generate} them.

The words in a free group are generated by successive application of ``operations''.
For the tautological description of a free group, these operations are just the operations of a group: multiplication, inversion, and the identity (a nullary operation).
When describing an arbitrary group we think of these operations as defined on a fixed underlying set, but when generating a free group we instead think of each of them as a ``way to produce new elements'', usually represented as syntactic strings of symbols.
That is, the elements of the free group are all the syntactic strings obtainable by successive application of the rules
\begin{compactenum}[(i)]
\item Given elements $X$ and $Y$, we have an element $(X Y)$,
\item Given an element $X$, we have an element $(X^{-1})$, and
\item We have an element $e$.
\end{compactenum}
Formally, this is an \emph{inductive definition}: the elements are the smallest set of syntactic strings closed under the rules.
Usually we think of applying these operations starting with a set of \emph{generators}, but an equivalent description that generalizes better is to include each generator as another nullary operation:
\begin{compactenum}[(i)]\setcounter{enumi}{3}
\item For any generator $g$, we have an element $g$.
\end{compactenum}
That is, generators are a special case of operations.
This allows us to describe the ``reduced words'' version of a free group by a similar set of inductive operations:
\begin{compactenum}[(i)]
\item We have an element $e$,
\item For any generator $g$ and any element $X$ not ending with $g^{-1}$, we have an element $X g$, and
\item For any generator $g$ and any element $X$ not ending with $g$, we have an element $X g^{-1}$.
\end{compactenum}

The judgments in type theory are likewise generated inductively by operations, which are usually called \emph{rules}.
Categorically, they build new objects and morphisms from old ones, generally according to them some universal property.
For example, there is a rule saying that any two types have a coproduct (disjoint union).
This rule applies in any context (categorically, all slice categories inherit coproducts from the base category);
type theorists write it as
\begin{equation}
  \inferrule{\Gamma \types A\ty \\ \Gamma\types B\ty}{\Gamma\types (A+B)\ty}\label{eq:coprod-form}
\end{equation}
As with judgments, this notation takes practice to read.
The horizontal bar separates the ``inputs'' (called \emph{premises}), on top, from the ``output'' (or \emph{conclusion}), on the bottom.
Each input or output is a judgment-in-context, and the inputs are separated by wide spaces or linebreaks.
If the operations for the tautological description of a free group were written analogously, they would be
\begin{mathpar}
  \inferrule{X\elt\\Y\elt}{(X Y) \elt} \and \inferrule{X\elt}{(X^{-1})\elt} \and \inferrule{\quad}{e\elt}
  \and \inferrule{g\text{ is a generator}}{g\elt}
\end{mathpar}
Here ``$X \elt$'' is the judgment that $X$ is an element of the free group, analogous to the judgments ``$A\ty$'' and ``$x:A\types b:B$'' that $A$ is an object and $b$ a morphism in a free category.
Note that the nullary identity operation has \emph{no} premises.
Similarly, the operations for the reduced-words description would be
\begin{mathpar}
  \inferrule{\quad}{e\elt}\and
  \inferrule{X\elt \\ g\text{ is a generator} \\ X\text{ doesn't end with } g^{-1}}{(Xg)\elt}\and
  \inferrule{X\elt \\ g\text{ is a generator} \\ X\text{ doesn't end with } g}{(Xg^{-1})\elt}
\end{mathpar}

The variables $X$ and $Y$ are analogous to $\Gamma$, $A$, $B$, and $b$ in type theory.
We call the latter \emph{meta-variables} to distinguish them from the variables $x:A$ occurring \emph{in} a context $\Gamma$, which have no analogue in group theory.

Returning to the coproduct type $A+B$, for it to be worthy of the name ``coproduct'', it needs to have certain structure.
There should be injections from $\m A$ and $\m B$ into $\m {A+B}$, which 
it may seem natural to write as
\begin{equation}
  \inferrule{ }{x:A\types \inl(x) :(A+B)}
  \qquad \text{and}\qquad
  \inferrule{ }{y:B\types \inr(y) :(A+B)}.\label{eq:inl-nocut}
\end{equation}
(We omit $\Gamma \types A\ty$ and $\Gamma\types B\ty$ from the premises, since these are implied\footnote{Depending on technical details far beyond our present scope, this implication might be a theorem about type theory, or it might be just an unproblematic abuse of notation.} by mention of $A+B$.)
However, one usually uses the following rules instead:
\begin{equation}
  \inferrule{\Gamma \types a:A}{\Gamma\types \inl(a) :(A+B)}
  \qquad \text{and}\qquad
  \inferrule{\Gamma \types b:B}{\Gamma\types \inr(b) :(A+B)}.\label{eq:inl-cut}
\end{equation}
Intuitively, this corresponds to describing the morphism $\m\inl:\m A \to \m{A+B}$ indirectly in terms of its image under the Yoneda embedding: for any morphism $\m a :\m\Gamma \to \m A$, we have an induced morphism $\m {\inl(a)} :\m\Gamma \to \m {A+B}$.

The reason for this is analogous to the difference between the tautological and reduced-words description of a free group.
In the tautological free group, we insert the generators as nullary operations, so that we need multiplication as a separate rule, which requires quotienting by an equivalence relation to enforce desired properties of multiplication.
But with reduced words, the generator rules $X \leadsto Xg$ and $X \leadsto Xg^{-1}$ incorporate ``just enough multiplication'' that we can \emph{define} multiplication of reduced words, and \emph{prove} that it is associative and so on, without needing to pass to a quotient.
Similarly, if we used rules like~\eqref{eq:inl-nocut}, we would need a separate ``substitution rule'' such as
\[ \inferrule{\Gamma \types a:A \\ \Gamma,x:A \types b:B}{\Gamma \types b[a/x]:B[a/x]} \]
to obtain composition in the classifying category, and to then quotient by an equivalence relation to make composition associative and unital.
But if we instead use rules like~\eqref{eq:inl-cut} that incorporate just the right amount of composition, one can \emph{define} substitution as on operation on judgments and \emph{prove} that it is associative and so on, thereby reducing or eliminating the need to pass to a quotient in constructing the classifying category.
Type theorists refer to this as the \emph{admissibility of substitution}; it is closely related to \emph{cut-elimination}.

The rules~\eqref{eq:inl-cut} are called the \emph{introduction} rules for the coproduct (they ``introduce'' elements of $A+B$), whereas~\eqref{eq:coprod-form} is called the \emph{formation} rule.
We also have an \emph{elimination} rule, which expresses the ``existence'' part of the categorical universal property of a coproduct:
\begin{equation}
  \inferrule{\Gamma,x:A \types c_A: C \\ \Gamma,y:B \types c_B: C \\ \Gamma \types s:A+B}
{\Gamma \types \case(C,c_A,c_B,s) : C}\label{eq:coprod-nondepelim}
\end{equation}
That is, given morphisms $\m A \to \m C$ and $\m B \to \m C$, we have a morphism $\m {A+B} \to \m C$.
The notation $\case(C,c_A,c_B,s)$ suggests that it is defined by inspecting the element $s$ of $A+B$ and dividing into cases: if it is of the form $\inl(x)$, then we use $c_A$, whereas if it is of the form $\inr(y)$, then we use $c_B$.%
\footnote{Technically, we should really write something like $\case(C,x.c_A,y.c_B,s)$, to indicate which variables $x$ and $y$ are being used in the terms $c_A$ and $c_B$.}

More generally, we allow $C$ to be a dependent type:
\begin{equation}
  \inferrule{\Gamma, z:A+B \types C \ty\\ \Gamma \types s:A+B \\\\ \Gamma,x:A \types c_A: C[\inl(x)/z] \\ \Gamma,y:B \types c_B: C[\inr(y)/z] }
{\Gamma \types \case(C,c_A,c_B,s) : C[s/z]}\label{eq:coprod-depelim}
\end{equation}
Categorically, this says that given a map $\m C\to \m{A+B}$ and sections of its pullbacks to $\m A$ and $\m B$, we can define a section over $\m{A+B}$ by the universal property of $\m{A+B}$.
This generalization of the existence part of the universal property is actually an equivalent way to include the \emph{uniqueness} part of it.
On one hand, categorically, uniqueness is what tells us that the induced map $\m{A+B}\to \m{C}$ is in fact a section.
On the other hand, assuming~\eqref{eq:coprod-depelim}, if $z:A+B \types c:C$ and $z:A+B \types d:C$ have equal composites with $\inl$ and $\inr$, then we can express this using the the ``equality type'' from \cref{sec:identity-types}:
\begin{mathpar}
  x:A \types e_A : c[\inl(x)/z] = d[\inl(x)/z] \and y:B \types e_B : c[\inr(y)/z] = d[\inr(y)/z]
\end{mathpar}
and then use~\eqref{eq:coprod-depelim} to construct $z:A+B \types e:(c=d)$.

Finally, the universal property also requires that $\m{\case(C,c_A,c_B)} \circ \m\inl$ equals $\m{c_A}$, and similarly for $\inr$.
These are called \emph{computation rules}.
We postpone discussing the exact meaning of ``equals'' here until \cref{sec:compute,sec:equality}; notationally we write it as $\jdeq$.
\begin{equation*}
  \inferrule{\Gamma, z:A+B \types C \ty\\ \Gamma,x:A \types c_A: C[\inl(x)/z] \\ \Gamma,y:B \types c_B: C[\inr(y)/z] }
  {\Gamma \types \case(C,c_A,c_B,\inl(a)) \jdeq c_A[a/x]}
\end{equation*}

In conclusion, we have four groups of rules relating to coproducts: formation (how to build types), introduction (how to build elements of those types), elimination (how to use elements of those types to build elements of other types), and computation (how to combine introduction and elimination).
Most rules of type theory come in packages like this, associated to one ``type constructor'' (here the coproduct) and expressing some universal property.
Given any class of \emph{structured categories} determined by universal properties, we can obtain a corresponding type theory by choosing all the corresponding packages of rules.
By and large, the rules for each type constructor are self-contained, allowing them to be ``mixed and matched''; thus
unlike ZFC, type theory is not a fixed system of axioms or rules, but a ``modular'' framework for such systems.

\begin{table}
  \centering
  \begin{tabular}{c|c|c}
    \textbf{Type constructor} & \textbf{Universal property} & \textbf{See}\\\hline
    coproduct types $A+B$ & binary coproducts & \cref{sec:rules}\\
    empty type $\emptyset$ & initial object\\
    product types $A\times B$ & binary products\\
    unit type $\unit$ & terminal object\\
    natural numbers $\N$ & natural numbers object & \cref{sec:hits} \\
    identity type $(x=y)$ & diagonal/equalizer & \cref{sec:myst-ident-types}\\
    function type $A\to B$ & exponential object (cartesian closure) & \cref{sec:rules}\\
    dependent sum $\sum_{x:A} B$ & left adjoint to pullback\\
    dependent product $\prod_{x:A} B$ & right adjoint to pullback (lcc)\\
    proposition type $\Omega$ & subobject classifier (elementary topos) & \cref{sec:universes}\\
    universe type $\U$ & object classifier ($\infty$-topos) & \cref{sec:universes}\\
    coequalizer type $\coeq(f,g)$ & coequalizer & \cref{sec:hits}\\
  \end{tabular}
  \caption{Type constructors and their semantics}
  \label{tab:type-constructors}
\end{table}

The most common type constructors and their corresponding universal properties are shown in \cref{tab:type-constructors}.
We will discuss some of these further in later sections, as indicated.
Here we give the rules explicitly only for \emph{function types}, which correspond to categorical exponentials; see \cref{fig:function-types}.
The exponential object from $A$ to $B$ is often denoted by $B^A$ or ${}^A B$, but in type theory we denote it by $A\to B$; this way the notation for its elements is $f:A\to B$, matching the usual notation for ``$f$ is a function from $A$ to $B$''.
Again we have a formation rule saying when $A\to B$ is a type, an introduction rule saying how to produce terms in $A\to B$, an elimination rule saying how to use such terms (by applying them to an argument), and two computation rules.
\begin{figure}
  \centering
  \begin{mathpar}
  \inferrule{\Gamma\types A\ty \\ \Gamma \types B\ty}{\Gamma\types (A\to B)\ty}\and
  \inferrule{\Gamma,x:A\types b:B}{\Gamma\types \lambda x.M:A\to B}\and
  \inferrule{\Gamma\types f:A\to B \\ \Gamma \types a:A}{\Gamma\types f(a) : B}\and
  \inferrule{\Gamma,x:A\types b:B \\ \Gamma \types a:A}{\Gamma\types (\lambda x.b)(a) \jdeq b[a/x]}\and
  \inferrule{\Gamma\types f:A\to B}{\Gamma\types f \jdeq (\lambda x.f(x))}
\end{mathpar}
\caption{The rules for function types}
\label{fig:function-types}
\end{figure}
Categorically, the elimination rule yields 
an ``evaluation'' morphism $\m{A\to B} \times \m A \to \m B$, while the introduction rule says that any map $\m \Gamma \times \m A \to \m B$ has a ``transpose'' $\m \Gamma \to \m{A\to B}$.
The computation rules say that these compose correctly.%
\footnote{We need two computation rules because we can't prove the uniqueness part of the universal property for functions the way we did for coproducts.
  See also
  \cref{sec:voev-univ-axiom}.}

An important generalization of this is the \emph{dependent function type}, where the codomain $B$ is allowed to depend on the domain $A$.
For instance, the family of generators of cyclic groups $n:\Np \types g_n : C_n$ yields a dependent function, assigning to each $n$ the generator of $C_n$:
\[ \lambda n.g_n \;:\; \tprod_{n:\Np} C_n \]
Categorically, 
given $\Gamma\types A\ty$ and $\Gamma,x:A\types B\ty$, the type $\Gamma\types \tprod_{n:A} B \ty$ is obtained from $\m{\Gamma,x:A,y:B} \to \m{\Gamma,x:A}$ by applying the right adjoint of pullback along $\m{\Gamma,x:A} \to \m{\Gamma}$.
Such a right adjoint exists exactly when the category is locally cartesian closed (LCC).

Of the type constructors in \cref{tab:type-constructors}, LCC categories also have products $A\times B$, a terminal object $\unit$, diagonals represented by the identity type $(x=y)$,%
\footnote{When we discuss the identity type in \cref{sec:myst-ident-types}, we will see that there are multiple choices that can be made for its rules.
For present purposes, the rules include ``UIP'' (see \cref{sec:myst-ident-types}).}
and left adjoints to pullback represented by dependent sums $\sum_{x:A} B$.
The latter generalizes $A\times B$; its elements are pairs $(x,y)$ where the type of $y$ can depend on $x$.
This collection of type constructors, corresponding to locally cartesian closed categories, is one of the most ``standard'' type theories, sometimes called \textbf{extensional Martin-L\"of Type Theory (EMLTT)} without universes (though that phrase also sometimes includes coproducts and the empty type).



In general, choosing a particular collection of type constructors specifies the rules for a particular type theory $\fT$, and thereby the collection of derivable judgments.
From this we construct a classifying category $\CtxT$ as in \cref{sec:type-term-judgments}, which
one can prove to be \emph{initial} among categories with the corresponding structure.
For instance, the classifying category of extensional Martin-L\"of Type Theory, as above, is the initial locally cartesian closed category.
It follows that the types and terms we construct in type theory have unique interpretations as objects and morphisms in \emph{any} category with appropriate structure, by applying the unique structure-preserving functor out of the initial object.

As in the case of groups, we often want to generalize this by including ``generators'' in addition to operations, allowing us to reason about \emph{arbitrary} objects and morphisms that may not necessarily be constructible ``from nothing'' using the categorical structure present.
And as we did for groups, we can do this by adding standalone rules to our type theory.
For instance, if we add the rules
\begin{equation}\label{eq:XYfrules}
  \inferrule{ }{\Gamma\types \mathsf{X}\ty} \qquad
  \inferrule{ }{\Gamma\types \mathsf{Y}\ty} \qquad
  \inferrule{ }{\Gamma\types \mathsf{f}:\mathsf{X}\to \mathsf{Y}}  
\end{equation}
to EMLTT, we obtain a type theory whose classifying category is the free locally cartesian closed category generated by two objects $\mathsf X$ and $\mathsf{Y}$ and a morphism $\mathsf{f}:\mathsf{X}\to \mathsf{Y}$.
(Here $\mathsf{X}$, $\mathsf{Y}$, and $\mathsf{f}$ are ``constants'', distinct from both variables and meta-variables.)
Thus, given any other locally cartesian closed category $\C$ in which we have chosen two objects $A,B$ and a morphism $g:A\to B$, there is a unique map from this classifying category sending $\mathsf{X}$, $\mathsf{Y}$, and $\mathsf{f}$ to $A$, $B$, and $g$ respectively.
Thus, anything constructable in type theory with the additional rules~\eqref{eq:XYfrules} can be interpreted in $\C$ \emph{and yield a result relative to $A$, $B$, and $g$}.

Rules like~\eqref{eq:XYfrules} appear to correspond only to \emph{free} groups, whereas we generally also consider \emph{presented} groups, with relations in addition to generators. 
However, as we saw in \cref{sec:identity-types}, equalities can be represented by elements of identity types; thus here the ``presented'' case is subsumed by the ``free'' case.


Constant rules like~\eqref{eq:XYfrules} allow us to reason about small collections of data in arbitrary structured categories.
In addition, given a category $\C$, there is a way to reason about ``all of $\C$ and nothing else'', by adding a constant for \emph{every} object, morphism, and equality in $\C$.
This yields a type theory $\Lang(\C)$ called the \textbf{internal language} of $\C$.
It is closely related to the ``Mitchell--Benabou language'' and ``Kripke--Joyal semantics'' of a topos~\cite{mm:shv-gl,ptj:elephant,bell:topos-lst,goldblatt:topoi}, which are type theories of a sort, but unlike ours do not include general dependent types.

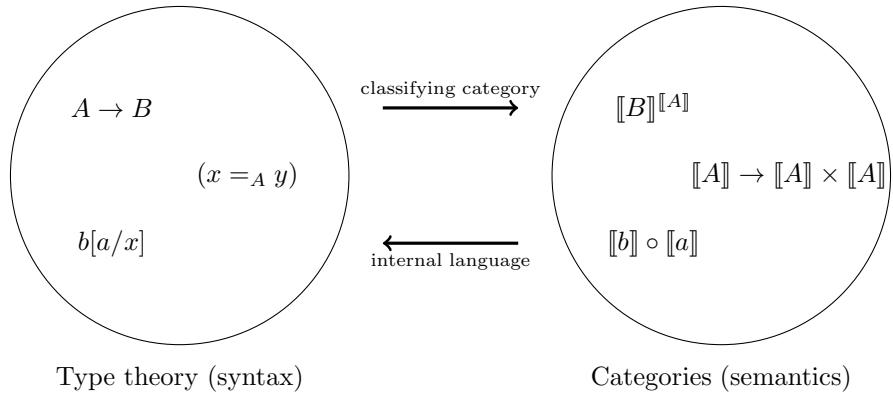
\begin{figure}
  \centering
  \begin{tikzpicture}[scale=.9]
    \node at (-4,0) {Type theory (syntax)};
    \node at (4,0) {Categories (semantics)};
    \draw (-4,3) circle [radius=2.5]; 
    \draw (4,3) circle [radius=2.5];
    \draw[->,very thick] (-1,4) -- node[auto] {\scriptsize classifying category} (1,4);
    \draw[->,very thick] (1,2) -- node[auto] {\scriptsize internal language} (-1,2);
    \node at (-5,4) {$A\to B$};
    \node at (3,4) {$\m{B}^{\m{A}}$};
    \node at (-3,3) {$(x=_Ay)$};
    \node at (5,3) {$\m A \to \m A \times \m A$};
    \node at (-5,2) {$b[a/x]$};
    \node at (3,2) {$\m b \circ \m a$};
  \end{tikzpicture}
  \caption{Syntax and semantics}
  \label{fig:internal-logic}
\end{figure}

If we fix some collection of type constructors, corresponding to a notion of structured category, then we can define a ``category of type theories'' based on these type constructors (with varying choices of constants).
Then $\Ctx$ becomes a functor from this category to an appropriate category of structured categories, and $\Lang$ into a right adjoint (or perhaps, depending on how we define the category of type theories, an inverse equivalence) to $\Ctx$; see \cref{fig:internal-logic}.
The counit of this adjunction is a functor $\Ctx(\Lang(\C)) \to \C$ that interprets the internal language of $\C$ in $\C$ itself; this gives a ``complete'' syntax for constructions in $\C$, analogous to the canonical presentation of a group $G$ involving one generator for each element and one relation for each equality.

I have just sketched an appealing general picture of the correspondence between type theories and categories.
However, proving the correctness of this picture can be exceedingly technical.
Seely's original proposal~\cite{seely:lccc-tt} contained a subtle technical flaw, later fixed by Hofmann and others~\cite{curien:subst,hofmann:ttinlccc,hofmann:ssdts,cd:lccc-tt}. 
But even now, complete proofs of the freeness of $\CtxT$ are quite involved and hard to find: they exist for some collections of type constructors~\cite{ls:hocl,jacobs:cltt,ptj:elephant,cd:lccc-tt,streicher:semtt}, and everyone expects all other cases to be analogous, but at present there is no general theorem. 
Indeed, even a precise definition of ``any type constructor'' is still lacking in the literature.
This is a current research problem, but I expect it to be solved one day; so I will say no more about this issue, except for a brief discussion of the (even more difficult) higher-categorical situation in \cref{sec:semantics}.

\subsection{Subobject classifiers and universes}
\label{sec:universes}

%

There is one class of type constructors, called \emph{universes}, that merits some individual discussion.
The simplest of these is a \emph{subobject classifier}, which categorically is a monomorphism $\true:1\to \Omega$ of which every other monomorphism is uniquely a pullback. 
In the category of sets, $\Omega = \{\true,\false\}$ is the set of truth values, and for a subset $A\subseteq B$ we have $A = \chi_A^{-1}(\true)$, where $\chi_A:B\to \Omega$ is the characteristic function of $A$.

If we identify $\true$ and $\false$ with a 1-element set and a 0-element set respectively, then up to isomorphism, the characteristic function of $A\subseteq B$ sends each $b\in B$ to its preimage under the inclusion $A\hookrightarrow B$.
This leads us to represent $\Omega$ by a type whose elements \emph{are themselves types}, with a rule like
\begin{equation}
  \inferrule{\Gamma \types P:\Omega}{\Gamma \types P\ty}\label{eq:omega-el}
\end{equation}
In particular, we have $x:\Omega \types x\ty$; and any other instance of~\eqref{eq:omega-el} can be obtained from this ``universal case'' by substitution.
Semantically, the interpretation of $\Gamma \types P:\Omega$ is a morphism $\m P : \m \Gamma \to \m \Omega$, while the interpretation of $\Gamma \types P\ty$ is an object of the slice over $\m\Gamma$, i.e.\ a projection morphism $\m{\Gamma, x:P} \to \m \Gamma$.
These two morphisms fit into a pullback square: 
\[ \xymatrix@-.7pc{
\m{\Gamma,x:P} \ar[r] \ar[d] \pullback &
\m{x:\Omega, y:x} \ar[d] \\
\m{\Gamma} \ar[r]_{\m P} & \m\Omega.
}\]
Thus, the morphism on the right (the interpretation of $x:\Omega \types x\ty$) is the universal monomorphism $\true:1\to \m\Omega$ for a subobject classifier (and in particular, $\m{x:\Omega, y:x}$ is a terminal object).

The fact that any $\Gamma \types P:\Omega$ classifies a monomorphism means equivalently that the types in $\Omega$ (the fibers of the corresponding objects of $\Ctx/\m\Gamma$) should ``have at most one element''.
We thus express it by the following rule:
\begin{equation}
  \label{eq:omega-hprop}
  \inferrule{\Gamma\types P:\Omega \\ \Gamma \types a:P \\ \Gamma \types b:P}{\Gamma\types \mathsf{tr}_P(a,b) : a=b}
\end{equation}
This says that the diagonal $\m{\Gamma,x:A} \to \m{\Gamma,x:A,y:A}$ has a section; hence it is an isomorphism and so $\m{\Gamma,x:A} \to \m{\Gamma}$ is mono.
(The notation $\mathsf{tr}_P$ stands for ``truncation''; see \cref{sec:prop-as-types}.)
The universality of $\true:1\to \m\Omega$ means that any type with ``at most one element'' is equivalent to one in $\Omega$:\footnote{The astute reader may notice that something is missing; we will return to this in \cref{sec:homotopy-type-theory}.}
\begin{gather}
  \label{eq:omega-classif}
  \inferrule{\Gamma\types P\ty \\ \Gamma,x:P,y:P \types p:x=y}
  {\Gamma \types \mathsf{Rsz}(P,xy.p):\Omega } 
  \\
  \inferrule{\Gamma\types P\ty \\ \Gamma,x:P,y:P \types p:x=y}{\Gamma\types \mathsf{rsz}_{P,xy.p} :
    (\mathsf{Rsz}(P,xy.p) \to P)\times (P\to \mathsf{Rsz}(P,xy.p))}
\end{gather}
The notations $\mathsf{Rsz}$ and $\mathsf{rsz}$ stand for ``resize'', indicating that $P$ may be ``too big to fit inside'' $\Omega$, but there is an equivalent type that does.

More generally, we can consider a \emph{universe type} $\U$, whose elements are types without any monomorphy restriction.\footnote{It is common in type theory to denote $\U$ by ``$\mathsf{Type}$'', and similarly to denote $\Omega$ by ``$\mathsf{Prop}$''.
  The latter will make more sense in \cref{sec:prop-as-types}.}
That is, we have the analogue of~\eqref{eq:omega-el}:
\[ \inferrule{\Gamma \types P:\U}{\Gamma \types P\ty} \]
but no analogue of~\eqref{eq:omega-hprop}.
The direct analogue of~\eqref{eq:omega-classif} would yield in particular $\mathsf{Rsz}(\U):\U$,
making the theory inconsistent due to Cantorian-type paradoxes.
Instead we assert that $\U$ is closed under the \emph{other} type constructors, e.g.
\begin{mathpar}
  \inferrule{\Gamma\types P:\U \\ \Gamma,x:P \types Q:\U}{\Gamma\types \tprod_{x:P} Q:\U}\and
  \inferrule{\Gamma\types P:\U \\ \Gamma,x:P \types Q:\U}{\Gamma\types \tsm_{x:P} Q:\U}\and
  \inferrule{\Gamma \types P:\U \\ \Gamma \types a:P \\ \Gamma\types b:P}{\Gamma \types (a=b):\U}
\end{mathpar}
Thus $\U$ is similar to a set-theoretic ``Grothendieck universe'' or inaccessible cardinal.
Of course, we can also have many universes $\U_i$ of different sizes.

Categorically, subobject classifiers are characteristic of elementary toposes, while universe objects arise in algebraic set theory~\cite{jm:ast,ast-website}.
But we will see in \cref{sec:homotopy-type-theory} that universes really come into their own when we pass to \io-categories and incorporate Voevodsky's \emph{univalence axiom}.

\subsection{Toposes of spaces}
\label{sec:toposes-spaces}

The preceding general theory tells us that given any category of spaces, if we choose a collection of type constructor ``packages'' corresponding to universal properties that exist in that category, the resulting type theory can be used to reason ``internally'' about that category.
Turning this around, for each ``package'' we want to include in our type theory, there is a corresponding restriction on the categories of spaces in which we can model it. 

Starting with the least restrictive case, \emph{simple} type theory (i.e.\ no dependent types) requires only a category with finite products, which includes  
practically any category of spaces. 
Nothing further is required to interpret binary product types $A\times B$ and the unit type $\unit$.
To interpret the coproduct type $A+B$ and the empty type $\emptyset$, we need a category with finite products and finite coproducts, with the former distributing over the latter (this is because the elimination rule for coproducts can be applied anywhere in a context, so that for instance $\m{x:A+B, z:C}$ also has the universal property of $\m{x:A,z:C} + \m{y:B,z:C}$).

To interpret 
dependent type theory, we require at least finite limits (since substitution into dependent types is interpreted by pullback).
This rules out a few examples such as smooth manifolds, but these can generally be embedded into larger categories having limits, such as ``generalized smooth spaces'' (see for instance Iglesias-Zemmour's chapter).
Nothing further is required to interpret the dependent sum type $\sum_{x:A} B$ and the identity type $(x=y)$.

The function type $A\to B$ in simple type theory can be interpreted in any cartesian closed category, but in \emph{dependent} type theory it requires a \emph{locally} cartesian closed category, since type constructors can be applied in any context (i.e.\ in any slice category). 
Local cartesian closure also allows us to interpret the dependent function type $\prod_{x:A} B$.

Ordinary topological spaces are not cartesian closed, but various slight modifications of them are.
The best-known of these, such as compactly generated spaces, are not \emph{locally} cartesian closed,
but there are others that are, such as subsequential spaces~\cite{ptj:topological-topos} (sets with a convergence relation between sequences and points) or pseudotopological spaces~\cite{choquet:convergences,wyler:topoi-in-topology} (similar, but using filters instead of sequences).
In such categories, 
a sequence of functions $f_n:X\to Y$ converges to $f_\oo:X\to Y$ in $Y^X$ if for any convergent sequence $x_n \leadsto x_\oo$ in $X$, the sequence $f_n(x_n)$ converges to $f_\oo(x_\oo)$ in $Y$; this is sometimes called \emph{continuous convergence}.
In fact, subsequential and pseudotopological spaces both form \emph{quasitoposes}~\cite{wyler:quasitopoi}, as do various kinds of generalized smooth spaces (see~\cite{bh:cc-smooth} and also Iglesias-Zemmour's chapter on diffeologies).


It is much less clear what sort of ``space'' could function as a subobject classifier or a universe.
One guess for a subobject classifier that doesn't work is the Sierpinski space $\Sigma$ (the set $\{\true,\false\}$ where $\{\true\}$ is open but $\{\false\}$ is not).
Continuous maps $X\to \Sigma$ classify \emph{open} subspaces of $X$; but not every mono is open. 
If instead we give the set $\{\true,\false\}$ the indiscrete topology, then it classifies arbitrary subspaces; but the monos of topological spaces (and their relatives) include all \emph{injective continuous functions}, which need not be subspace inclusions.

Thus, a subobject classifier $\Omega$ has to have sufficient structure that maps into it can encode chosen topologies on subsets.
For instance, if $Y\to X$ is mono and a sequence $(x_n)$ lying in $Y$ converges (in $X$) to a point $x_\oo$ also lying in $Y$, then it might or might not also converge to $x$ \emph{in the topology of $Y$}.
Thus, in defining a map $\chi_Y: X\to \Omega$ classifying $Y$, even after we know that $x_n$ and $x_\oo$ are sent to $\true$ (hence lie in $Y$), we need an additional degree of freedom in defining $\chi_Y$ 
to specify whether or not the convergence $x_n \to x_\oo$ is still ``present'' in $Y$.

Obviously this is impossible with classical topological spaces, 
but there are categories of spaces in which such an object exists:
the trick is to make the ``spatial structure'' into \emph{data} rather than a \emph{property}.
For instance, the quasitopos of subsequential spaces sits inside the topos of \emph{consequential spaces}.
A consequential space is a set equipped with, for every sequence $(x_n)$ and point $x_\oo$, a \emph{set} of ``reasons why'' or ``ways in which'' $(x_n)$ converges to $x_\oo$.
(Of course, this set might be empty, i.e.\ $(x_n)$ might not converge to $x_\oo$ at all.)
The axioms of a subsequential space are then promoted to operations on these ``witnesses of convergence'', which then have to satisfy their own axioms.
See~\cite{ptj:topological-topos} for details.%
\footnote{The term ``consequential space'' should not be blamed on~\cite{ptj:topological-topos}, who considered it but then discarded it.
  I have chosen to use it anyway, since I know no other term for such spaces.}

The category of consequential spaces is a topos, so it is locally cartesian closed and has a subobject classifier.
The latter has two points $\{\true,\false\}$, but many different ``witnesses'' that the constant sequence at $\true$ converges to $\true$, allowing the characteristic function of a mono $Y\to X$ to retain information about the topology of $Y$.
One might think we only need \emph{two} such witnesses, to record whether a convergent sequence $x_n\leadsto x_\oo$ in $X$ also converges in $Y$; but in fact we need to record which \emph{subsequences} of $\{x_n\}$ also converge to $x_\oo$ in $Y$.
We can exactly determine the witnesses of convergence in $\Omega$ by its universal property: they must be the sub-consequential-spaces of the ``universal convergent sequence'' $\N_\oo$ (the one-point compactification of $\N$).
See~\cite[Corollary 4.2]{ptj:topological-topos}.

The category of consequential spaces also has universe objects $\U$.
Roughly speaking, this means that the collection of all consequential spaces (bounded in size by some cardinality) can be made into a consequential space.
By the desired universal property of $\U$, a witness that a sequence of spaces $(X_n)$ converges to $X_\oo$ consists of a consistent way to make $\big(\coprod_{n} X_n\big)\sqcup X_\oo$ into a consequential space over $\N_\oo$, which roughly means giving a consistent collection of witnesses for convergence of sequences $x_n \in X_n$ to points $x_\oo\in X_\oo$.
In~\cite{es:universe-indiscrete} it is shown that for any $(X_n)$ and $X_\oo$ there is at least one such witness, so the topology of $\U$ is ``indiscrete'' in some sense (though unlike for a classical indiscrete space, interesting information is still carried by \emph{how many} such witnesses there are).

Generalizing from consequential spaces, for any small collection of spaces $\T$, we can build a topos whose objects are ``spaces'' whose ``topology'' is determined by ``probing'' them with maps out of $\T$.
More precisely, we take the category of sheaves for some Grothendieck topology on $\T$.
Consequential spaces are the case when $\T = \{\N_\oo\}$, so that a space is determined by its ``convergent sequences''.

Another reasonable choice is $\T = \{\R^n\}_{n\in \N}$, yielding spaces determined by a notion of ``continuous paths, homotopies, and higher homotopies''.\footnote{To be precise, we take the full subcategory of spaces $\T = \{\R^n\}_{n\in \N}$ with the Grothendieck topology of open covers, whereas consequential spaces are obtained from the one-object full subcategory $\{\N_\oo\}$ with the ``canonical'' Grothendieck topology.
  This difference in topology is the actual source of a significant amount of the differences between the two toposes.}
Urs Schreiber has suggested to call these \textbf{continuous sets}.
Just as a sequence in a consequential space can ``converge to a point'' in more than one way, a path in a continuous set can ``be continuous'' in more than one way.
By the adjunction, continuous paths $\R \to Y^X$ are equivalently homotopies $\R\times X\to Y$.
Similarly, in $\U$, the witnesses to continuity of a ``path of spaces'' $\{X_t\}_{t\in\R}$ are the ways to make $\coprod_{t\in\R} X_t$ into a space over $\R$, which roughly means a consistent collection of witnesses for the continuity of ``paths'' consisting of points $x_t\in X_t$ for all $t\in \R$.
Consequential spaces and continuous sets are similar in many ways, but different in others, and each has advantages and disadvantages.

If we restrict the morphisms of $\T$ to preserve some additional structure on its objects, we obtain a topos whose objects have a version of that structure.
For instance, if we choose $\T = \{\R^n\}_{n\in \N}$ but with only the \emph{smooth} functions between them, then we obtain a topos whose objects are a kind of generalized smooth space, in which any given path or homotopy has a set of ``witnesses to smoothness'' rather than witnesses to continuity; we call these \textbf{smooth sets}.
Just as the quasitopos of subsequential spaces sits inside the topos of subsequential spaces, the quasitopos of diffeological spaces (and hence also the usual category of smooth manifolds) sits inside the topos of smooth sets.

There are likewise toposes of ``algebraic'' spaces, as well as toposes containing ``infinitesimal'' or ``super'' spaces; see~\cite{bh:cc-smooth,moerdijk-reyes:sia} and Schreiber's chapter.

\subsection{Why Spaces?}
\label{sec:why-spaces}

We can now address the question mentioned in \cref{sec:introduction} of \emph{why} so many objects in mathematics come naturally with spatial structure.
The short answer is that since type theory can be a language for all of mathematics, any construction of a mathematical object can be phrased in type theory and then interpreted into any of these toposes of spaces, thereby yielding not just a set but a space.
%
%
%

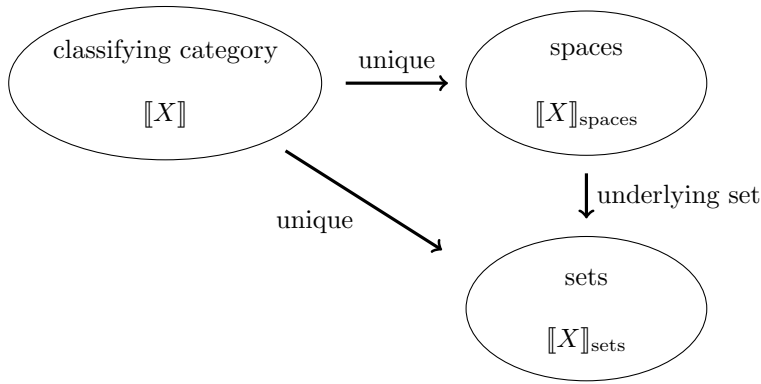
\begin{figure}
  \centering
  \begin{tikzpicture}[yscale=.6,xscale=.8]
    \draw (0,3) ellipse (2.6 and 1.7);
    \node at (0,3.7) {classifying category};
    \node at (0,2.3) {$\m X$};
    \draw[->,very thick] (3,3) -- node[auto] {unique} (4.7,3);
    \draw (7,3) ellipse (2 and 1.6);
    \node at (7,3.7) {spaces};
    \node at (7,2.3) {$\m{X}_{\text{spaces}}$};
    \draw[->,very thick] (2,1.5) -- node[auto,swap] {unique} (4.6,-.7);
    \draw[->,very thick] (7,1) -- node[auto] {underlying set} (7,0);
    \draw (7,-2) ellipse (2 and 1.6);
    \node at (7,-1.3) {sets};
    \node at (7,-2.7) {$\m{X}_{\text{sets}}$};
  \end{tikzpicture}
\caption{A hypothetical logical underlying-set functor}
\label{fig:logical-underlying}
\end{figure}

However, this answer is missing something important: to say that a set $X$ comes naturally with spatial structure is to say that we have a canonically defined space \emph{whose underlying set is} $X$, and simply being able to interpret the construction of $X$ in a topos of spaces doesn't ensure this.
If the ``underlying set'' functor from our topos of spaces to the category of sets preserved all the structure used to interpret type theory (limits, colimits, dependent exponentials, universes, and so on --- this is called a \emph{logical functor}), then this would follow from the initiality of the classifying category.
Namely, the unique functor from the classifying category to the category of sets would factor uniquely through our topos of spaces, implying that the image of $X$ in sets would be the underlying set of the image of $X$ in spaces, as shown in \cref{fig:logical-underlying}.

In general, the underlying-set functor does not preserve \emph{all} the structure.
But it does preserve quite a lot of it.
In fact, for all the categories of spaces considered above, the underlying-set functor fits into a string of adjunctions:
\[
\xymatrix{ \text{topos of spaces} \ar[d]^\Gamma \\ \text{topos of sets}
  \ar@<10mm>[u]^{\Delta} \ar@<5mm>@{}[u]|{\dashv} \ar@<-6mm>@{}[u]|{\dashv} \ar@<-10mm>[u]_{\nabla}  }
\]
where $\Gamma$ is the underlying-set functor, $\Delta$ equips a set with a ``discrete topology'', and $\nabla$ equips it with an ``indiscrete topology''.
(In consequential spaces, the discrete topology says that only eventually-constant sequences converge, while the indiscrete one says that every sequence converges (uniquely) to every point.)
Moreover, $\Delta$ preserves finite limits, and $\Delta$ and $\nabla$ are fully faithful; this makes the category of spaces into a \emph{local topos} over sets (see~\cite{jm:local-toposes} or~\cite[\S C3.6]{ptj:elephant}).

In particular, this string of adjunctions implies that $\Gamma$ preserves all limits and colimits.
In the above examples, $\Gamma$ also preserves the subobject classifier and the universes --- this is clear from our explicit descriptions of these above, and categorically it implies that the spaces are also a \emph{hyperconnected topos}.

The main thing that $\Gamma$ doesn't preserve is function-spaces.
However, it does preserve function-spaces whose domain is discrete: for any set $X$ and space $Y$ we have $\Gamma(Y^{\Delta X}) \cong (\Gamma Y)^X$, and likewise for dependent exponentials.
(This follows formally from the fact that $\Delta$ preserves finite limits.)

The upshot is that any mathematical construction can be interpreted as a space, and as long as the only functions it uses have a domain that (when interpreted as a space) is discrete, the resulting space will be a topology on the set that we originally thought we were defining.
This restriction makes sense: once we start using functions with non-discrete domain, in the world of spaces we have to consider only \emph{continuous} functions, causing a divergence from the world of sets.
Moreover, in this case often the world of spaces is the ``correct'' one: we \emph{should} restrict to continuous maps between profinite groups, and use only continuous homomorphisms in the Pontryagin dual.

The next natural question to ask is, when \emph{does} a given mathematical construction inherit a non-discrete topology?
Categorically, this means asking which constructions are (not) preserved by $\Delta$.
Since $\Delta$ is a left adjoint, it preserves colimits, and as remarked it preserves finite limits; but in general it doesn't preserve much more than this.
Thus, nontrivial topologies can arise from (1) infinite limits,
(2) function-spaces (which, in the category of sets, are just particular infinite products), or (3) the subobject classifier or a universe.

Many constructions that automatically inherit spatial structure are obtained by infinite limits.
For instance, the profinite completion is a limit of finite quotients.
In the topos of consequential spaces, this infinite limit generally has its expected topology, since ``take the convergent sequences'' is a limit-preserving functor from classical topological spaces to consequential spaces.

By contrast, for the topos of continuous sets (and similarly smooth sets), $\Delta$ \emph{does} preserve infinite limits, and indeed it has a further left adjoint (see \cref{sec:cohes-homot-type}).
Thus, the profinite completion in these toposes gets the discrete topology (intuitively, profinite topologies are incompatible with the ``manifold-like'' structure of a continuous set).
More generally, nontrivial continuous-set structures can arise \emph{only} in constructions that include the subobject classifier or a universe.
Perhaps the most important example of such a space is the real numbers $\R$;
but before discussing them, we need to talk about logic in type theory.

\section{Towards synthetic topology}
\label{sec:synthetic-topology}

In \cref{sec:syntax} I claimed that type theory presents a free ``universe of mathematics''. 
So far, we have seen that type theory contains type constructors such as products, coproducts, and exponentials that look like the standard operations on sets.
Moreover, when taken literally (rather than as ``code'' for objects and morphisms in a category), the syntax of type theory talks about \emph{elements} of types, and the rules stipulate that the elements of product types, exponential types, and so on are exactly what we would expect.
Thus, any mathematical construction that is classically performed with sets, such as building the rational and real numbers out of the integers, can be performed with types instead. 

However, there is more to mathematics than constructing things: we also like to \emph{prove} things about them.
The formal system of ZFC set theory is formulated inside of first-order logic, so that proving is the ``basic act of mathematics'':
what are called ``constructions'' are actually just existence (or existence-and-uniqueness) proofs.
But in type theory, \emph{constructions} are the ``basic act of mathematics''; so what has happened to proofs?

\subsection{Propositions as types}
\label{sec:prop-as-types}


In set theory, a property $P$ of elements of a set $A$ can be equivalently expressed as a subset of $A$, namely $\{ x\in A \mid P(x) \}$.
Categorically, this is a monomorphism into $A$, or equivalently its characteristic function $A\to\Omega$.
This provides us with the means to \emph{define logic inside of type theory}: we declare that by a \textbf{property} of elements of a type $A$ we \emph{mean} a judgment $x:A \types P:\Omega$.
Similarly, by a \textbf{proposition}\footnote{We use ``proposition'' in the logician's sense of something that one might \emph{try} to prove, rather than the other common meaning of something that \emph{has been} proven.
  Note that in~\cite{hottbook}, types with at most one element are called \emph{mere propositions} rather than just ``propositions''.}
we mean simply an element of $\Omega$: that is, a type having \emph{at most one element}.
We regard such a proposition as ``true'' when it \emph{does} have an element.
That is, the proposition corresponding to $P:\Omega$ is ``$P$ has an element''.
Similarly, by a \textbf{proof} of $P:\Omega$ we mean a construction of an element of $P$.

We have already seen an example of this way of representing properties: the equality type $x:A,y:A \types (x=y)\ty$, which
in \cref{sec:identity-types} we said represents the diagonal $\diag_A : \m A \to \m A \times \m A$ as an object of $\Ctx/(\m A \times \m A)$.
Since $\diag_A$ is a monomorphism, it is classified by a map $\m A \times \m A \to \Omega$, which is the binary relation of equality.
In the category of sets, this map sends $(x,x)$ to $\true$, since the fiber of $\diag_A$ over $(x,x)$ has one element, and sends $(x,y)$ to $\false$ if $x\neq y$, since then the fiber of $\diag_A$ over $(x,y)$ is empty.
(But see \cref{sec:homotopy-type-theory}.)

This identification of propositions with ``types having at most one element'' is close to, but not quite, the usual meaning of the phrase \emph{propositions as types}.
The latter refers to allowing \emph{any} type to be called a proposition, rather than only those with at most one element.
However, our choice (which is increasingly common) relates better to the standard practice in mathematics whereby once a proposition has been proven, the particular proof given has no further mathematical (as opposed to aesthetic or conceptual) importance.

One reason the identification of propositions with (certain) types is so con\-venient is an observation called the \emph{Curry--Howard correspondence}~\cite{curry:curry-howard,howard:curry-howard,martinlof:itt-pred,wadler:pat}: the operations of logic are \emph{already present} in type theory as constructions on types.
For instance, if $P$ and $Q$ are propositions, then so is $P\times Q$; and since it has an element just when $P$ and $Q$ do, it is natural to call it ``$P$ and $Q$''. 
Similarly, $P\to Q$ is ``if $P$ then $Q$'', since a function $f:P\to Q$ transforms the truth of $P$ into the truth of $Q$, while $\prod_{x:A} P(x)$ is ``for all $x:A$, $P(x)$'', since a dependent function $f:\prod_{x:A} P(x)$ assigns to any $x:A$ a proof of $P(x)$.
It is also reasonable to regard $P\to \emptyset$ as ``not $P$'', since a function $f:P\to \emptyset$ can only exist if $P$ is empty (i.e.\ false).

We might expect $P+Q$ to be ``$P$ or $Q$'', but $P+Q$ may not be a proposition even if $P$ and $Q$ are.
For this purpose we introduce a new type constructor called the \emph{propositional truncation}\footnote{Propositional truncation has a long history and many variations, with names such as \emph{squash type}, \emph{mono-type}, and \emph{bracket type}; see e.g.~\cite{nuprlbook,mendler:quotient-types,maietti:tt-hpretop,ab:bracket-types,hottbook}.}, whose rules are shown in \cref{fig:brck}.
(The notation $\mathsf{ptr}$ just stands for ``propositional truncation''.)
\begin{figure}
  \centering
  \begin{mathpar}
    \inferrule{\Gamma\types A\ty}{\Gamma\types \brck A: \Omega}\and
    \inferrule{\Gamma\types A\ty\\ \Gamma\types a:A}{\Gamma\types \bproj{a}:\brck{A}}\and
    \inferrule{\Gamma\types A\ty\\ \Gamma\types B:\Omega \\ \Gamma,x:A\types b:B \\ \Gamma\types a:\brck{A}}{\Gamma\types \mathsf{ptr}(x.b,a) : B}\and
    \inferrule{\Gamma\types A\ty\\ \Gamma\types B:\Omega \\ \Gamma,x:A\types b:B \\ \Gamma\types a:A}{\Gamma\types \mathsf{ptr}(x.b,\bproj a) \jdeq b[a/x]}
  \end{mathpar}
  \caption{Propositional truncation}
  \label{fig:brck}
\end{figure}
Intuitively, $\brck A$ is the proposition ``$A$ has at least one element''; while categorically, $\Gamma\types \brck A: \Omega$ presents the \emph{image} of the projection $\m{\Gamma,x:A} \to \m{\Gamma}$.
The introduction rule says that if we have an element of $A$, then $A$ has at least one element.
The elimination rule says that if we know that $A$ has at least one element, then when proving a proposition we may assume \emph{given} an element of $A$.
(Removing the hypothesis that $B$ is a proposition would imply a choice principle that is too strong even for classical mathematics; see \cref{fn:globalchoice} on page \pageref{fn:globalchoice}.)

Now we define ``$P$ or $Q$'' to be $\brck{P+Q}$, and similarly ``there exists an $x:A$ such that $P(x)$'' to be $\brck{\sum_{x:A} P(x)}$.
As observed by Lawvere~\cite{lawvere:adjointness}, this definition of the existential quantifier can be described categorically as the left adjoint to pullback between posets of subobjects $\mathrm{Sub}(\m\Gamma) \to \mathrm{Sub}(\m{\Gamma,x:A})$.
The untruncated $\sum_{x:A}$ gives the left adjoint to the pullback between slice categories $\Ctx/\m\Gamma  \to \Ctx/\m{\Gamma,x:A}$, and the truncation reflects it back into monomorphisms.
Similarly, the universal quantifier ``for all $x:A$, $P(x)$'' is the right adjoint of the same functor: since the right adjoint $\prod_{x:A}$ between slice categories already preserves monomorphisms, no truncation is necessary.

Rather than being the existential quantifier, the untruncated $\sum_{x:A} P(x)$ plays the role of the subset $\{ x\in A \mid P(x) \}$.
Its elements are pairs of an element $x:A$ and a proof that $P(x)$ holds, and since $P(x)$ is a proposition, to give an element of $P(x)$ contains no more information than that $P(x)$ ``is true''.
Thus we may consider the elements of $\sum_{x:A} P(x)$ to be ``the elements of $A$ such that $P(x)$ is true''.
The type of \emph{all} subsets of $A$, mentioned in \cref{sec:type-term-judgments}, is just ``$A\to\Omega$''.

In conclusion, instead of describing mathematics inside of logic, as is done by ZFC set theory, in type theory we ``define logic inside of mathematics''.
One advantage of this is that, as we have seen, type theory can be interpreted in suitably structured categories, yielding an intrinsic ``logic'' internal to any such category.
Another is that it allows us to draw finer distinctions, thereby actually representing informal mathematics \emph{more} faithfully, in the following way.

We mentioned above that there is a tradition in type theory that interprets ``there exists an $x:A$ such that $P(x)$'' as $\sum_{x:A} P(x)$ rather than $\brck{\sum_{x:A} P(x)}$.
In addition to its mismatch with the standard practice of mathematics, when done indiscriminately this can actually lead to inconsistencies; see~\cite[\S3.2]{hottbook} and~\cite{escardo-xu:brouwer-ch}.
However, there \emph{are} places in ordinary mathematics where a ``there exists'' statement is more naturally interpreted as $\sum_{x:A} P(x)$.
For instance, the Yoneda lemma $\mathrm{Nat}(\C(-,a),F) \cong F(a)$ doesn't mean that there merely \emph{exists} such an isomorphism but that we have specified a \emph{particular} one.
ZFC set theory, being formulated inside first-order logic, forces every ``theorem'' to be a ``mere existence'' statement; but type theory frees us from this straitjacket, allowing us to directly express ``constructions'' in addition to ``proofs''.

\subsection{Constructive logic}
\label{sec:constructive-logic}

When the logical connectives are defined according to the Curry--Howard correspondence, most of the basic laws of logic can be derived from the rules for the basic type constructors.
For instance, one of de Morgan's laws
\[ (\neg P \land \neg Q) \to \neg(P\lor Q) \]
(where as usual $\land,\lor,\neg$ mean ``and'', ``or'', and ``not'') has the following proof:
\[ \lambda x.\lambda y. \mathsf{ptr}(\case(\emptyset,\proj1(x)(y),\proj2(x)(y),y),y) : (\neg P \land \neg Q) \to \neg(P\lor Q). \]

(Of course, it would be a heavy burden to carry around such long terms whenever we want to use de Morgan's law.
But writing out a fully formalized proof in ZFC is no easier --- in fact, often it is much harder!
In both cases the formalism simply justifies our ordinary mode of mathematical writing.)

However, we cannot derive from the rules of type theory any proofs of the following classical tautologies of logic:
\begin{mathpar}
  \neg(P\land Q) \to (\neg P \lor \neg Q)\and
  \neg \neg P \to P \and
  P\lor \neg P \and
  (\neg\forall x, P(x)) \to (\exists x, \neg P(x))
\end{mathpar}
In other words, the logic we obtain from propositions-as-types is \emph{constructive} or \emph{intuitionistic logic}.
Constructive logic acquired a bad name due to some fundamentalists in the early 20th century, but in fact it is natural and unavoidable once we recognize that type theory is a syntax for categories, and that sets are not the only category in the world (see~\cite{bauer:5stages}).
The above tautologies are simply \emph{not true} in most categories of spaces. 

For example, consider the proposition $\prod_{x:A} P(x) \lor \neg P(x)$, where $A$ is a type and $P:A\to\Omega$ a property, and let us interpret this in a ``topological'' topos such as those discussed in \cref{sec:toposes-spaces}.
Now $P$ classifies a monomorphism $B\rightarrowtail A$, which as mentioned previously need not be a subspace embedding.
Similarly, $\neg P$ classifies a different monomorphism $\neg B\rightarrowtail A$, which turns out to be the ``maximal mono disjoint from $B$''.
In other words, $\neg B$ contains all the points of $A$ that are absent from $B$, and also ``all the topology'' on those points that is absent from $B$ (e.g.\ all the convergent sequences, or all the continuous paths).
However, $(\lambda x. P(x)\lor \neg P(x)) : A \to \Omega$ classifies their \emph{union} $B\cup \neg B$ as monos into $A$, which is not generally a subspace even if $B$ is: it contains all the points of $A$, but its topology is that of the disjoint union $B \sqcup \neg B$.
Thus the mono $B\sqcup \neg B \to A$ has no \emph{continuous} section, and so we cannot assert $\prod_{x:A} P(x) \lor \neg P(x)$.

In other words, constructive logic is simply {more general} than classical logic.
As always, using fewer assumptions --- here, assumptions about \emph{logic} --- leads to a more general conclusion --- here, one that applies to more categories.

There is also another sense in which constructive logic is more general: classical logic can be \emph{embedded} in constructive logic.
Specifically, the subset $\Omega_{\neg\neg}$ of $\Omega$ consisting of those $P$ such that $\neg\neg P \to P$ (formally, the $\Sigma$-type $\sum_{P:\Omega} (\neg\neg P\to P)$) admits logical operations satisfying the laws of classical logic.
In fact, $\Omega_{\neg\neg}$ is closed under all the ordinary logical operations except for ``or'' and ``there exists'', and we can define $P\mathbin{\lor'} Q$ to be $\neg\neg(P\lor Q)$ and $\exists' x:A$ to be $\neg\neg \exists x:A$.
(Note the similarity to how in \cref{sec:prop-as-types} we applied $\brck{-}$ to $+$ and $\sum$ to get $\lor$ and $\exists$.)
In categories of spaces, the subtypes whose classifying map factors through $\Omega_{\neg\neg}$ generally coincide with the subspace embeddings.

Using $\Omega_{\neg\neg}$ gives us a different ``logic'' that is always classical, but is not as well-behaved in other ways.
For instance, it fails to satisfy \emph{function comprehension} (a.k.a.\ the \emph{principle of unique choice}): ``if for all $x:A$ there is a unique $y:B$ such that $P(x,y)$, then there is $f:A\to B$ such that $P(x,f(x))$ for all $x:A$.''
However, there is a subuniverse of types where $\Omega_{\neg\neg}$-logic does behave well.
Define a type $A$ to be a \emph{$\neg\neg$-sheaf} if the ``constant functions'' map $(\lambda x.\lambda p.x):A \to (P\to A)$ is an isomorphism for any $P:\Omega$ such that $\neg\neg P$.
The world of $\neg\neg$-sheaves in constructive mathematics behaves just like the world of classical mathematics, with both classical logic and function comprehension.
In categories of spaces, the $\neg\neg$-sheaves are generally the \emph{indiscrete} spaces.
(If you were expecting to hear ``discrete'' instead of ``indiscrete'', wait for \cref{sec:modalities-cohesion}.)

On the other hand, it is always possible to add classical logic \emph{globally} to type theory as a rule:
\[ \inferrule{\Gamma\types P:\Omega}{\Gamma\types \mathsf{lem}(P) : P \vee \neg P} \]
(This rule, called the \emph{Law of Excluded Middle (LEM)}, suffices to prove all the other classical tautologies.)
This would mean restricting our syntax to apply only to ``Boolean'' categories, such as sets, and excluding most topological ones, just as adding the relations $g h = h g$ to a free group turns it into a free \emph{abelian} group, with a more restricted universal property.
Similarly, we can add a type-theoretic version of the \emph{Axiom of Choice (AC)}, which is not provable%
\footnote{The reader may have heard rumors that the axiom of choice is actually provable in type theory \emph{without} any added axioms.
It is true that one can prove a statement that looks like the axiom of choice \emph{if} arbitrary types are allowed to play the role of ``propositions'', i.e.\ if all propositional truncations are removed from the definitions of the logical operations in \cref{sec:prop-as-types}.
But from our perspective, this provable statement is not at all the axiom of choice, since its hypothesis already essentially carries along the data of a choice function.}
even after adding LEM, and whose precise formulation we leave to the reader.%
\footnote{There is also a subuniverse like the $\neg\neg$-sheaves that satisfies AC as well as LEM: one can build G\"odel's ``constructible universe'' $L$ (no relation to ``constructive logic'') inside the $\neg\neg$-sheaves.
  However, this relies on first-order logic and ZFC-style membership-based set theory; no category-theoretic or type-theoretic construction of such a subuniverse is known.}

If we wanted only to use type theory merely as a foundation for classical mathematics, 
there would be no problem with this.\footnote{There \emph{are} many toposes, other than the category of sets, that satisfy both LEM and AC.
  They are roughly the same as forcing models of ZFC.
  However, none of them are ``spatial'' in the sense we care about here.}
But in this chapter our focus is on type theory as a syntax for categories of \emph{spaces}, which frequently means that we must learn to live with constructive logic (though we will see in \cref{sec:homotopy-type-theory} that ``homotopy spaces'' can be compatible with classical logic).
Fortunately, this is usually not very difficult; often it suffices to rephrase things carefully, avoiding unnecessary negations.
For instance, constructively it is not very useful to say that a type is nonempty ($\neg (A \cong \emptyset)$ or equivalently $\neg\neg A$); instead we use the positive statement that it is ``inhabited'' ($\exists x:A$, i.e.\ $\brck{A}$).
It also often happens that a group of classically equivalent definitions are no longer the same constructively, and we have to judiciously choose the ``correct'' one.
For instance, classically a set is finite just when it is not bijective to any proper subset of itself, but constructively  
this is a weaker and less useful condition; the correct definition of ``finite'' is ``bijective to $\{k:\N \mid k<n \}$ for some $n:\N$''.

Once we get over this minor hurdle, we can develop mathematics on top of type theory in basically the same way as usual.
Formally, type-theoretic proofs and constructions involve heavy manipulation of syntax (just as in any other formal foundational system like ZFC); but as mentioned above, when actually \emph{doing} mathematics there is usually no reason to bother about this.

As a simple example, once we have the natural numbers type $\N$, we can define the integers $\Z$ as $\N+\N$ (with appropriate structure), the rational numbers $\Q$ as a subtype of $\Z\times \N$ (the ``fractions $a/b$ in lowest terms''), and the real numbers $\R$ as a subtype of $(\Q\to\Omega)\times (\Q\to\Omega)$ (the two-sided Dedekind cuts).
Recall that by a ``subtype of $A$'' we mean a type of the form $\sum_{x:A} P(x)$ where $P:A\to\Omega$ is a property; for instance, more formally we have
\[ \Q \;=\; \textstyle\sum_{r:\Z\times \N} \Big(\proj2(r)>0 \land \prod_{n:\N} \big((n\mid \proj1(r)) \land (n\mid\proj2(r)) \to n=1\big)\Big). \]
Here $\proj1$ and $\proj2$ are the projections out of a cartesian product, and $>$ and $\mid$ are relations that we have to define previously, e.g.
\begin{align*}
  (n\mid m) &=  \big(\exists p:\N, (p \cdot n = m)\big) \\
  &= \brck{\textstyle\sum_{p:\N} (p \cdot n = m)}
\end{align*}
where multiplication has been previously defined, etc.
We can then proceed to define all the usual functions and properties of numbers of all sorts, and build the rest of mathematics on top of them.

When this syntax is interpreted into the category of sets, it of course yields the usual sets of numbers.
The point, however, is that if we instead interpret it into a category of spaces, the types of numbers automatically inherit a spatial structure, and usually that spatial structure is the \emph{intended} one!
For instance, in the toposes of consequential spaces or continuous sets from \cref{sec:toposes-spaces}, the real numbers type $\R$ is interpreted by the real numbers \emph{with their usual Euclidean topology} (see~\cite[Proposition 4.4]{ptj:topological-topos} and~\cite[Theorem VI.9.2]{mm:shv-gl}).

Note that the definition of $\R$ using Dedekind cuts (which we may denote $\R_d$ for emphasis) fulfills the requirement for a nontrivial continuous-set structure, since it uses $\Omega$.
However, not all the classically-equivalent definitions of $\R$ remain equivalent constructively.
For instance, the \emph{Cauchy reals} $\R_c$, defined by taking equivalence classes of Cauchy sequences, come with an inclusion $\R_c \to \R_d$ that is not generally surjective.
In consequential spaces, we have $\R_c \cong \R_d$; but in continuous sets, $\R_c$ gets the discrete topology.
Thus, it is usually better to regard the Dedekind real numbers as ``the'' real numbers. 

Since the Dedekind reals $\R$ have their usual topology in our toposes, other types built from them, such as the circle $\topcirc = \setof{ (x,y):\R\times\R \mid x^2+y^2=1 }$, the complex numbers $\mathbb{C} \cong \R\times \R$, or matrix groups $\mathrm{GL}_n(\R) \subseteq \R^{n^2}$,
also have their usual topologies.
Furthermore, all functions definable in type theory are interpreted by continuous maps; so the fact that we can define addition of real numbers in type theory tells us automatically that $\R$ is a topological group, and so on.
Analogous facts are true for the other constructions leading to nontrivial topologies mentioned in \cref{sec:why-spaces}, such as profinite completion.
Thus, we have finally made good on our promise from \cref{sec:introduction} to provide a formal system for describing ``groups with background spatial structure'' that is sufficiently flexible to include all different kinds of spatial structure at once, with the added benefit of a uniform way of constructing the standard examples.

The fact that ``the real numbers'' defined in type theory are interpreted in some categories by the usual \emph{space} of real numbers has an interesting consequence: \textit{constructively, we cannot define\footnote{To be precise, we cannot define such a function ``in the empty context'', i.e.\ without any ambient assumptions.} any discontinuous function $\R\to\R$}.
In particular, the usual examples of discontinuous ``piecewise'' functions $\R\to\R$, such as the Heaviside step function
\[ f(x) =
\begin{cases}
  0 &\quad \text{if } x<0\\
  1 &\quad \text{if } x\ge 0
\end{cases}\]
cannot be defined constructively --- or, more precisely, their domain cannot be shown constructively to be all of $\R$ (that being tantamount to the assertion that every real number is either $<0$ or $\ge 0$, which is essentially an instance of LEM).
That is, restricting ourselves to constructive logic automatically ``notices'', and forces us to respect, a canonical and implicit topological structure on types such as $\R$.
In the next section we briefly discuss another such implicit structure.

\subsection{A digression on computation}
\label{sec:compute}

While our primary concern here is with the suitability of type theory as a ``logic of space'', historically it developed rather differently.
The first type theory, which bore little resemblance to its modern descendants, was introduced by Russell to avoid his eponymous paradox.
After other logicians such as G\"odel refined Russell's type theory in various ways, Church~\cite{Church:1940tu,Church:1941tc} combined it with his ``$\lambda$-calculus'' to obtain what today we can see as a typed functional programming language.
The dependent type theory we are using here 
is mainly due to Martin-L\"{o}f~\cite{martinlof:itt-pred,martinlof:itt}, whose intent was to give ``a full scale system for formalizing intuitionistic mathematics'' in the sense of Bishop~\cite{bishop:fca,bb:constr-analysis}.
Bishop, in turn, wanted to develop a form of mathematics in which all statements would have computational meaning, so that for instance whenever we assert something to \emph{exist} we must have a method for finding it.
This led him, following the earlier pioneering work of Brouwer, to reject the law of excluded middle, since in general there can be no \emph{method} for deciding which of $P$ or $\neg P$ holds.

Thus, type theory was originally conceived as a formal basis for a mathematics that would be ``constructive'' in this computational sense.
It is remarkable that it turned out to also be a flexible system for reasoning in arbitrary categories!
The existence of internal languages of categories was apparently first recognized for elementary toposes~\cite{ls:itt-free-topos,mitchell:topoi-sets,bj:logique-des-topos}; the ``Mitchell--Benabou language'' 
of a topos is a sort of type theory in which the only dependent types are those in $\Omega$.
The generalization to Martin-L\"{o}f's type theory was first written down by Seely~\cite{seely:lccc-tt}, and corrected and refined by others as discussed in \cref{sec:rules}.

From our category-theoretic point of view, the computational aspect of type theory is partly explained by a different class of models.
In addition to categories of ``spaces'' where every map is continuous, there are categories of ``computable objects''~\cite{hyland:eff,oosten:realiz} in which every map is computable. 
Thus, everything in type theory
must be ``potentially computable'' in addition to ``potentially continuous''.
In particular, just as any constructively definable function $\R\to\R$ must be continuous, any constructively definable function $\N\to\N$ must be computable.

However, there is more to the computational side of type theory than this: its syntax \emph{is} actually a programming language that can be executed.
That is, not only does every term $f:\N\to\N$ represent a computable function, but its definition is an algorithm for computing that function.
The ``execution'' of such a program is essentially the type-theoretic version of the ``reduction'' algorithm for free groups that simplifies $(g((hg)^{-1}((hh)g^{-1})))h^{-1}$ to $hg^{-1}h^{-1}$.

Specifically, these reductions implement the ``computation rules'' mentioned in \cref{sec:rules} (hence the name).
For instance, the first computation rule in \cref{fig:function-types} can be interpreted as a ``reduction'' or ``normalization'' step allowing us to ``simplify'' $(\lambda x.M)(N)$ to $M[N/x]$.
(This partially explains why we used a different equality symbol $\jdeq$; see also \cref{sec:equality}.)
In good cases, these reduction steps are guaranteed to terminate at a unique ``value'' or ``normal form'', analogous to the reduced words for elements of a free group.
Reduction in a free group is a fairly simple process, but since type theory is complicated enough to encode all of mathematics, its notion of ``reduction'' can serve as a general-purpose programming language.%
\footnote{This programming language is not technically ``Turing-complete'', since all its ``programs'' must terminate; otherwise we could prove a contradiction with a divergent computation.
But it can still encode all computable functions, e.g.\ with a partiality monad~\cite{capretta:genrec-coind,danielsson:opsem-partiality}.}

This makes type theory a convenient language for reasoning about computer programs, and as such it has many adherents in computer science.
Moreover, we can implement a ``compiler'' for type theory on a physical computer, which then also serves as a \emph{proof checker} for mathematical arguments.
Thus, mathematics done in type theory not only can be interpreted in many categories, but can have its correctness formally verified in this way.
Such computerized ``proof assistants'' built on type theory play an increasingly important role in computer science, and are slowly growing in importance in mathematics.

One thing to note is that the computational interpretation of type theory is rather ``fragile'': the guarantee of unique termination (called ``strong normalization'' or ``canonicity'') depends on the details and interactions of all the rules.
New axioms can cause computation to get ``stuck'' and never reach a value, while new computation rules can cause the reduction algorithm to loop or diverge.
In particular, LEM and AC cause computation to stick (and semantically, they also rule out computational models along with topological ones).
This is not surprising since they assert that certain things exist or are true without giving any way to construct them.
However, none of this prevents us from using computer proof assistants to check proofs involving such noncomputational axioms; it only means that such proofs can't then be ``executed''.

Since this book is about space rather than computation, I will not say much more about the computational side of type theory.
However, it is worth pointing out that computation and topology are actually closely related.
As first recognized by Scott~\cite{scott:cts-lattices,scott:dt-lat,scott:ttalt}, computational objects often come naturally with, or are represented by, topologies representing the fact that a finite computation can only consume a finite amount of data.
In fact, Brouwer's original intuitionism was also arguably more ``topological'' than computational.
For further discussion of these ideas, see for instance~\cite{escardo:syntop-datatypes,vickers:topology-via-logic}.

\subsection{Synthetic topology}
\label{sec:synthetic-topology-1}

So far we have seen that the types in type theory admit interpretations as various different kinds of space.
Thus, one might say that they have \emph{latent} or \emph{potential} spatial structure: they {might} be spaces, but they also might not have any nontrivial spatiality, depending on where we interpret them.
Moreover, they also have other latent structures, such as computability.

Until now we have considered mainly the aspects of types that are {independent} of their potential spatial structure, where the topology simply comes along for the ride.
However, some spatial aspects of types are visible {inside} of type theory, without needing to interpret them first in some category.
This leads to subjects called \emph{synthetic topology} and \emph{synthetic differential geometry}.

One important observation is that in many cases we can {detect} topology using structures that already exist in type theory.
For instance, following~\cite{es:universe-indiscrete}, define $\N_\oo$ to be the type of non-increasing binary sequences:
\[ \N_\oo \;=\; \textstyle\sum_{a:\N\to \bool} \prod_{n:\N} (a_{n+1} \le a_n) \]
where $\bool$ is the ``Boolean'' type with two elements $0$ and $1$.
Then we have an injection $i:\N\to\N_\oo$ where $i(m)_n = 1$ if $m<n$ and $0$ otherwise, and we also have an element ``$\oo:\N_\oo$'' defined by $\oo_n = 1$ for all $n$.
In the topos of consequential spaces, $\N_\oo$ is interpreted by the ``actual'' one-point compactification of $\N$; thus it is sensible to \emph{define} a \textbf{convergent sequence} in a type $A$ to be a map $\N_\oo \to A$.
In this way, without assuming any axioms, we see that every type {automatically} has a structure like a consequential space, and every function is {automatically} ``continuous'' in the sense of preserving convergent sequences.

Of course, the actual interpretation of $\N_\oo$ depends on the category: in the topos of continuous sets, it yields $\N\sqcup\{\oo\}$ with the discrete topology, so that every sequence ``converges'' uniquely to every point.
However, the structure of continuous sets is detectable internally in a different way: define a \textbf{continuous path} in a type $A$ to be a map $\R \to A$ out of the (Dedekind) real numbers.
This gives the expected answer for both consequential spaces and continuous sets, since in both cases $\R$ has its usual topology.
We will come back to this in \cref{sec:cohes-homot-type}.

These internally defined ``topologies'' are only ``potentially nontrivial'': e.g.\ if we assume LEM, then every sequence converges uniquely to every point and every path is uniquely continuous.
If we want to ensure that they definitely \emph{are} nontrivial, we can assert ``nonclassical axioms'' that contradict LEM, excluding the category of sets but retaining topological models.
For instance, we could assert that the only ``convergent sequences'' in $\R$ are those that converge in the $\epsilon$-$N$ sense, or that every ``continuous path'' in $\R$ is continuous in the $\epsilon$-$\delta$ sense.\footnote{The latter statement is sometimes known as \emph{Brouwer's theorem}, since Brouwer proved it in his ``intuitionistic'' mathematics using principles derived from ``choice sequences''.  It is \emph{not} a theorem of pure constructive mathematics, which unlike Brouwer's ``intuitionism'' is fully compatible with classical principles like excluded middle, though it does not assume them.}

Convergent sequences and continuous paths are ``covariant'' notions of topology, i.e.\ they are defined using maps \emph{into} a type.
We can also describe ``contravariant'' notions of topology synthetically, involving maps \emph{out} of a type.
For instance, with continuous sets in mind, 
we can define \emph{open subsets} as the preimages of open intervals under functions $A\to\R$.
Alternatively, we can construct or postulate a subtype $\Sigma\subseteq\Omega$ behaving like the Sierpinski space (usually called a \emph{dominance}), and define an \emph{open subset} to be one whose classifying map $A\to\Omega$ factors through $\Sigma$.
In a topos of sheaves on a category $\T$ of spaces, there is an obvious choice of such a $\Sigma$, namely the sheaf represented by the actual Sierpinski space (whether or not it is in $\T$).
On the other hand, if we want to construct a particular $\Sigma$ \emph{inside} type theory, one possibility is the \emph{Rosolini dominance}~\cite{rosolini:dominance}:
\[ \Sigma_{\mathrm{Ros}} \;=\; \Big(\textstyle\sum_{P:\Omega} \exists f:\N\to\bool, \Big(P \cong \exists n:\N, (f(n)=1)\Big)\Big) \]
That is, $\Sigma_{\mathrm{Ros}}$ is the type of propositions of the form $\exists n:\N, (f(n)=1)$ for some $f:\N\to\bool$.
In consequential spaces, $\Sigma_{\mathrm{Ros}}$ is the Sierpinski space, so the resulting ``open subsets'' are as we would expect.
(But in continuous sets, $\Sigma_{\mathrm{Ros}} \cong \bool$, so the only ``open subsets'' in this sense are unions of connected components.)
However we choose $\Sigma$, once chosen we can develop much of classical point-set topology with these ``synthetic open sets'', including compactness, Hausdorffness, and so on.
See~\cite{escardo:syntop-datatypes,escardo:top-hoil,taylor:lamcra,bl:met-syntop} for more examples of this sort of ``synthetic topology''.
There are also other ways to define open sets synthetically; for instance,~\cite{penon:thesis} defines $U\subseteq A$ to be open if $\forall x:A,\forall y:A, (x\in U \to (y\in U \lor x\neq y))$.

It is more subtle to obtain a synthetic theory of \emph{smoothness}, since smoothness does not arise automatically from constructive logic the way that continuity and computation do.
Indeed, it is easy to define non-differentiable functions $\R\to\R$ in constructive mathematics, such as the absolute value.
Semantically, the type $\R_d$ of Dedekind reals interpreted in the topos of \emph{smooth} sets actually yields the sheaf of \emph{continuous} (not necessarily smooth) real-valued functions on the domain spaces.%
\footnote{In particular, while it is true in a sense that ``everything is smooth'' in this topos, what this actually means is that each object comes with a ``notion of smoothness'', which could in some cases happen to coincide with mere continuity.}

Clearly a more interesting smooth set than this is the sheaf of \emph{smooth} real-valued functions, which is equivalently the usual smooth manifold $\R$ regarded as a diffeological space.
In the internal language of smooth sets, this appears as a type $\R_s$ of ``smooth reals'' living strictly in between the ``discrete (Cauchy) reals'' $\R_c$ and the ``continuous (Dedekind) reals'' $\R_d$.
It seems unlikely that there is any type definable in type theory whose interpretation in smooth sets is $\R_s$, but we can at least write down some axioms that $\R_s$ satisfies, such as being a subring of $\R_d$, or more generally closed under the action of all ``standard smooth functions'' (see~\cite{fourman:smooth-reals} for one way to make this precise).

A more transformative approach is to make the notion of ``smoothness'' synthetic as well, rather than relying on the classical limit definition of derivative.
Following Grothendieck's insight into the importance of nilpotent elements in algebraic geometry, we can enhance the category $\T = \{\R^n\}_{n\in \N}$ by replacing each $\R^n$ by its algebra of smooth functions $C^\oo(\R^n)$ and turning the arrows around to obtain a category of $\R$-algebras, then adding new algebras that are ``deformations'' of some $C^\oo(\R^n)$ containing nilpotents.
Whatever the details, the resulting topos will contain an internal ring $\bR$ that enhances $\R_s$ to include nilpotent ``infinitesimals'', with $\R_s$ the quotient by these:
\[\xymatrix@-1pc{ & \bR \ar@{->>}[d] \\
\R_c \ar[r] & \R_s \ar[r] & \R_d. } \]
Nilpotents allow a synthetic definition of differentiation:
if $D\subseteq \bR$ is defined by $D = \{ d:\bR \mid d^2=0 \}$, then for any $f:D\to\bR$ there is a unique $f'(0):\bR$ such that $f(d) = f(0) + f'(0)\cdot d$ for all $d:D$.
(This is sometimes called the \emph{Kock--Lawvere axiom}.)
In particular, all functions $f:\bR\to\bR$ are differentiable in this synthetic sense.
The theory resulting from this and similar axioms that hold in the above sheaf toposes is called \textit{synthetic differential geometry}; for further reading see~\cite{bell:sdg,kock:sdg,moerdijk-reyes:sia,lavendhomme:sdg} and Kock's chapter herein.

This section has been a very brief sketch of some ways to access the spatial structure of types internally.
Rather than pursue any of these avenues in detail, I want to describe (in \cref{sec:cohes-homot-type}) a newer approach to synthetic topology that leverages more of the categorical and type-theoretic structure, involving \emph{higher modal operators}.
However, first we must move sideways to consider another very different latent structure in type theory: that of \emph{homotopy spaces} or \emph{\oo-groupoids}.

\section{Homotopy type theory}
\label{sec:homotopy-type-theory}

\subsection{The mystery of identity types}
\label{sec:myst-ident-types}

For many years, the most mysterious part of Martin-L\"{o}f's type theory was the identity types ``$x=y$''.
As mentioned in \cref{sec:identity-types,sec:prop-as-types}, the semantic idea is that the dependent type $x:A,y:A\types (x=y)\ty$ represents the diagonal $\diag_A:\m A \to \m A\times \m A$, regarded as an object of $\Ctx/(\m A\times \m A)$. 
Of course, $\diag_A$ is automatically present in the classifying category, as defined in \cref{sec:type-term-judgments}; but without the identity type it isn't represented by any dependent type.

The rules of the identity type are difficult to understand at first, but essentially they use a universal property to express the fact that it is the diagonal.
In fact, any object $f:B\to A$ of a slice category $\C/A$ has a universal property: it is the image of the terminal object $1_B:B\to B$ of $\C/B$ under the left adjoint $f_!:\C/B\to\C/A$ to pullback along $f$.
In other words, for any object $g:C\to A$ of $\C/A$, morphisms $f\to g$ in $\C/A$ are in natural bijection with sections of the pullback $f^*(g) : f^*C \to B$:
\[ \vcenter{\xymatrix{ f^*C \ar[d] \ar[r] \pullback & C \ar[d]^g \\ B \ar[r]^f \ar@/^3mm/@{.>}[u] & A }}
\quad\iff\quad
\vcenter{\xymatrix@C=1pc{ B \ar[dr]_f \ar@{.>}[rr] && C \ar[dl]^g \\ & A }}
\]
Thus, $\diag_A:A\to A\times A$ is characterized in $\C/(A\times A)$ by saying that for $g:C\to A\times A$, morphisms $\diag_A \to g$ are naturally bijective to sections of $\diag_A^*(g)$:
\[ \vcenter{\xymatrix{ \diag_A^*C \ar[d] \ar[r] \pullback & C \ar[d]^g \\ A \ar[r]^-{\diag_A} \ar@/^3mm/@{.>}[u] & A\times A }}
\quad\iff\quad
\vcenter{\xymatrix@C=1pc{ A \ar[dr]_{\diag_A} \ar@{.>}[rr] && C \ar[dl]^g \\ & A\times A }}
\]
If we represent $g$ by a dependent type $x:A,y:A \types C\ty$, then the pullback on the left corresponds to substitution of the same variable for both $x$ and $y$, e.g.\ $C[w/x,w/y]$ in context $w:A$.
The section on the left is then a term $w:A \types c:C[w/x,w/y]$; whereas the induced map on the right 
corresponds to a term $x:A,y:A,p:x=y \types d:C$.
If we represent the latter in Yoneda form, as we did for \eqref{eq:coprod-nondepelim}, we obtain the following rule:
\begin{equation}
  \inferrule{\Gamma,x:A,y:A \types C\ty \\ \Gamma, w:A \types c : C[w/x,w/y] \\\\
    \Gamma \types a:A \\ \Gamma \types b:A \\ \Gamma \types p:a=b}
  {\Gamma \types \J(C,c,p) : C[a/x,b/y]}\label{eq:id-rec}
\end{equation}
The corresponding introduction rule is just the unit of this adjunction, saying that the diagonal of $A$ has a specified section when pulled back along itself.
Logically, it expresses the reflexivity of equality:
\[ \inferrule{\Gamma \types A\ty \\ \Gamma \types a:A}{\Gamma \types \refl_a:a=a}\]
In the classifying category, this rule is a section as on the left below:
\[ \vcenter{\xymatrix{ \diag_A^*\Id_A \ar[d] \ar[r] \pullback & \Id_A \ar[d] \\
    \m A \ar[r]^-{\diag_A} \ar@/^3mm/@{.>}[u] & \m A\times \m A }}
\quad\iff\quad
\vcenter{\xymatrix@C=1pc{ \m A \ar[dr]_{\diag_A} \ar@{.>}[rr]^{\refl} && \Id_A \ar[dl] \\ & \m A\times \m A }}
\]
where $\Id$ denotes $\m{x:A,y:A,p:x=y}$.
Thus, $\refl$ gives a morphism as on the right above; the universal property should make this an isomorphism $\m A \cong \Id_A$.

This adjoint characterization of equality is due to Lawvere~\cite{lawvere:comprehension}, but it is closely related to Leibniz's ``indiscernibility of identicals''.
Specifically, given a property $x:A \types P(x):\Omega$, we can form $x:A, y:A \types (P(x)\to P(y)) :\Omega$.
Taking this as $C$, we have the second hypothesis of~\eqref{eq:id-rec} given by $w:A \types (\lambda p.p) : P(w) \to P(w)$.
Thus, from~\eqref{eq:id-rec} we get
\[x:A, y:A, p:x=y \types \J(\cdots) :  P(x) \to P(y).\]
This says that if $x$ and $y$ are equal (``identical'') then any property that holds of $x$ also holds of $y$ (``indiscernible'').
The function $P(x) \to P(y)$ induced by $p$ is often denoted $p_*$ and called \emph{substitution} or \emph{transport}.


In dependent type theory, we need to enhance~\eqref{eq:id-rec} to allow $C$ to depend on a ``witness of equality'' $p:x=y$ as well, and also add a computation rule relating it to the introduction rule $\refl_a$.
This yields the rules shown in \cref{fig:id}, which are due to Martin-L\"{o}f~\cite{martinlof:itt,martinlof:itt-pred}.
Just as for coproducts, we expect this enhancement to ensure the full universal property of the desired adjunction.
When phrased informally in terms of points, the stronger elimination rule says that if we want to perform a construction or proof involving a general element $p:x=y$ (for general $x$ and $y$), it suffices to consider the case when $y$ is $x$ and $p$ is $\refl_x$.
This is formally analogous to the elimination rule for (say) coproducts, which says that to perform a construction or proof involving a general element $z:A+B$, it suffices to consider the two cases when $z$ is $\inl(x)$ and $\inl(y)$.


\begin{figure}
  \centering
  \begin{mathpar}
  \inferrule{\Gamma \types A\ty \\ \Gamma \types a:A \\ \Gamma\types b:A}{\Gamma \types (a=b)\ty}\and
  \inferrule{\Gamma \types A\ty \\ \Gamma \types a:A}{\Gamma \types \refl_a:a=a}\and
  \inferrule{\Gamma,x:A,y:A,e:x=y \types C\ty \\ \Gamma, w:A \types c : C[w/x,w/y,\refl_w/e] \\
    \Gamma \types a:A \\ \Gamma \types b:A \\ \Gamma \types p:a=b}
  {\Gamma \types \J(C,c,p) : C[a/x,b/y,p/e]}\and
  \inferrule{\Gamma,x:A,y:A,e:x=y \types C\ty \\\\ \Gamma, x:A \types c : C[x/y,\refl_x/e] \\
    \Gamma \types a:A}
  {\Gamma \types \J(C,c,\refl_a)\jdeq c[a/x]}
\end{mathpar}
\caption{The rules for identity types}
\label{fig:id}
\end{figure}

It may seem odd to introduce a new notation and system of rules for an object (the identity type) that turns out to be isomorphic to something we already had (the diagonal map).
The point is that dependent type theory has two different ways of representing a morphism $A\to B$, depending on whether we view it as simply a morphism between two objects or as an object of the slice category over its codomain $B$.
The same is true in set theory: we can have a function $f:A\to B$ between sets, or we can have a $B$-indexed family of sets $\{A_b\}_{b\in B}$, and up to isomorphism the two are equivalent by $A= \coprod_b A_b$ and $A_b = \{ a\in A \mid f(a)=b \}$.
As mentioned in \cref{sec:type-term-judgments}, enabling us to think in terms of ``indexed families'' in an arbitrary category is actually one of the \emph{strengths} of dependent type theory, and
the identity type and $\Sigma$-type are exactly what supply the isomorphism between functions and families in type theory. 

However, although the rules in \cref{fig:id} are well-motivated categorically from a universal property, they appear to have a serious problem.
Specifically, if the identity type satisfies these rules only, and we define the classifying category as in \cref{sec:type-term-judgments}, then it does not have pullbacks!

Consider, for instance, how we might try to pull back a dependent type $\m{y:B, z:C(y)} \to \m{B}$ along a morphism $\m f : \m A \to \m B$.
In \cref{sec:type-term-judgments} we claimed that such a pullback should be obtained by substituting $f(x)$ for $y$ in $C(y)$.
To check the universal property, we would consider a diagram as below:
\[ \xymatrix{
  \m T \ar@/_5mm/[ddr]_{g} \ar@/^5mm/[drr]^{h} \ar@{.>}[dr]^{?}\\
  & \m{x:A, z:C(f(x))} \ar[r] \ar[d] \pullback &
  \m{y:B, z:C(y)} \ar[d] \\
  & \m{A} \ar[r]_{\m f } & \m B. } \]
Then $g$ is a term $t:T \types g(t):A$, while $h$ is determined by two terms $t:T \types h_1(x):B$ and $t:T \types h_2(t) : C(h_1(x))$.
Now it seems as though the commutativity of the square would force $h_1(t)$ to be $f(g(t))$, so that $h_2$ would be a term $t:T \types h_2(t): C[f(g(t))/y]$, or equivalently $t:T \types h_2(t): C[f(x)/y][g(t)/x]$, inducing the dotted morphism.

This appealing argument stumbles on the fact that we quotiented the morphisms in $\Ctx$ by an equivalence relation induced by the identity type.
Thus, to say that the above square commutes doesn't mean that $h_1(t)$ is literally $f(g(t))$, only that we have a term $t:T \types p : h_1(t) = f(g(t))$.
This is not by itself the end of the world, because $p$ induces a transport function $p_* : C(h_1(t)) \to C(f(g(t)))$, so we can define the dotted morphism as $t:T \types p_*(h_2(t)): C(f(g(t)))$.
The real problem is that this morphism depends on the choice of the term $p$, but the term $p$ is not \emph{specified} by the mere \emph{fact} that the outer square commutes; thus the dotted factorization is not unique.

I have stated this problem in a form that may make its solution seem obvious to a modern reader with a certain background.
However, for many years after Martin-L\"{o}f no one with this background looked at the problem; nor was it stated in this way.
Instead the question was whether we can prove \emph{inside type theory} that the identity type $x:B, y:B \types (x=y)\ty$ is a ``proposition'' in the sense of \cref{sec:universes,sec:prop-as-types}, i.e.\ any two of its elements are equal.
Categorically, this would mean proving that the projection $\m{x:B, y:B, p:x=y} \to \m{x:B,y:B}$ is a monomorphism (as we expect, if it is to be the diagonal); while syntactically it would mean constructing, given $p:x=y$ and $q:x=y$, a term of $p=q$.
If this were the case, different choices of $p$ would result in terms $p_*(h_2(t))$ that are equal, so that uniqueness for the pullback would be restored.

It turns out, however, that we \emph{cannot} prove that the identity type is always a proposition.
Thus, people (starting with Martin-L\"{o}f) considered adding this statement as an extra standalone axiom:
\begin{equation}
  \inferrule{
    \Gamma\types p:a=b \\\Gamma\types q:a=b}{\Gamma\types e:p=q}\label{eq:uip}
\end{equation}
called \emph{Uniqueness of Identity Proofs (UIP)}.
Type theory with UIP is sometimes called \emph{extensional},\footnote{More precisely, ``propositionally extensional''; see \cref{sec:equality}.} while type theory without UIP is called \emph{intensional}.
Thus, as mentioned in \cref{sec:rules}, it is \emph{extensional} Martin-L\"{o}f type theory for which our previous construction of $\CtxT$ presents free locally cartesian closed categories (or other kinds of structured categories).
However, baldly assuming UIP is unsatisfying, since it doesn't fit into the system of rule packages motivated by universal properties, as described in \cref{sec:rules}.
Moreover, this approach provides no insight into why UIP might be true, or why it isn't provable.

But to a modern reader with a background in homotopy theory, the above problem looks familiar: it is the same reason why the homotopy category of spaces doesn't have pullbacks.
In that case, we consider instead \emph{homotopy} pullbacks, where the factorization morphism is \emph{required} to depend on a choice of homotopy filling the square.
This suggests we should regard $\Ctx$ as a ``homotopy theory'' or $\io$-category, thereby explaining why we cannot prove that $x=y$ is a proposition: diagonals in an $\io$-category are \emph{not} in general monic.
For instance, in the 2-category of groupoids, the monomorphisms are the fully faithful functors, but a diagonal $G\to G\times G$ is not generally full: its functorial action on hom-sets
\[\hom_G(x,y) \longrightarrow \hom_{G\times G}((x,x),(y,y)) = \hom_G(x,y) \times \hom_G(x,y)\]
is not an isomorphism if $\hom_G(x,y)$ has more than one element.
The first model of type theory using this idea was constructed by~\cite{hs:gpd-typethy} using groupoids; later authors~\cite{aw:htpy-idtype,klv:ssetmodel} generalized it using homotopy theory.

This situation should be compared with the remarks about constructive logic in \cref{sec:constructive-logic}.
In both cases we have a rule (LEM or UIP) that seems reasonable given one model or class of models (the category of sets, or all 1-categories).
But this rule turns out not to be provable, because type theory admits more general models, in some of which the rule is false.
This provides us with the proper attitude towards the rule: assuming it simply means restricting the class of categories we are interested in, whereas declining to assume it allows us to use type theory as a syntax for a wider class of models.

Unfortunately, there is much unresolved subtlety in an $\io$-categorical interpretation of type theory.
One hopes for an analogue of the 1-categorical situation, with a ``classifying $\io$-category'' that is free in some \io-category of structured \io-categories, setting up an \io-adjunction, but it is quite difficult to make this precise.
In \cref{sec:semantics} I will sketch the current state of the art; until then I will just assume that the problem will be solved somehow, as I believe it will be.

\subsection{Types as \oo-groupoids}
\label{sec:oo-group-struct}

In \cref{sec:synthetic-topology-1} we saw that if we decline to assume LEM, we can detect the potential spatial structure of types internally. 
Similarly, if we decline to assume UIP, then types have potential ``homotopy space'' or ``\oo-groupoid'' structure, and the natural way to try to detect this is by using the identity types.
But what does $x=y$ mean when it is not just a proposition?

In higher category theory we have a notion of \emph{$n$-groupoid}, which is an \oo-groupoid containing no interesting information above dimension $n$.
This can be defined inductively: a $0$-groupoid is an \oo-groupoid that is equivalent to a discrete set, while an $(n+1)$-groupoid is one all of whose hom-\oo-groupoids $\hom_A(u,v)$ are $n$-groupoids.
Moreover, we can extend the induction downwards two more steps: an \oo-groupoid is a $0$-groupoid just when each $\hom_A(u,v)$ is empty or contractible, so it makes sense to define a \emph{$(-1)$-groupoid} to be an \oo-groupoid that is either empty or contractible.
Similarly, an \oo-groupoid is a $(-1)$-groupoid just when its homs are \emph{all} contractible, so we can define a \emph{$(-2)$-groupoid} to be a contractible one.
(See, for instance,~\cite[\S2]{bs:ncats-cohom}.)



In particular, when we regard a set as an \oo-groupoid, the proposition that two elements $u,v$ are equal turns into the $(-1)$-groupoid $\hom_A(u,v)$.
Thus, the homs of an \oo-groupoid generalize the notion of equality for elements of a set, so it is natural to expect the \emph{type} $u=v$ to behave like $\hom_A(u,v)$.
This is correct: we can derive all the composition structure on these hom-objects that should be present in an \oo-groupoid 
from the rules in \cref{fig:id}~\cite{pll:wkom-type,bg:type-wkom}.
For instance, we can construct the composition law
\[ (x=y) \times (y=z) \to (x=z) \]
by applying the eliminator to $p:x=y$ to assume that $y$ is $x$ and $p$ is $\refl_x$, in which case the other given $q:y=z$ has the same type as the goal $x=z$. 
(This is the same as the proof of transitivity of equality in extensional type theory.)

Other aspects of homotopy theory can also be defined using the identity type.
For instance, the \emph{loop space} $\Omega(A,a)$ of a type $A$ at a point $a$ is just the identity type $(a=a)$.
Voevodsky also showed that we can mimic the above inductive definition of $n$-groupoids, also called \emph{homotopy $n$-types}\footnote{Voevodsky's terminology~\cite{vv:unimath} is ``type of h-level $n+2$''.}: a type $A$ is an $(n+1)$-type if for all $x:A$ and $y:A$ the type $(x=y)$ is a $n$-type.
We can start at $n=-1$ with the propositions as defined in \cref{sec:prop-as-types}, i.e.\ types $A$ such that for all $x:A$ and $y:A$ we have $x=y$.
We can also start at $n=-2$ with the \emph{contractible types}, which are just the propositions that have an element.
Note that the homotopy $0$-types, also called \emph{sets}, are those that satisfy UIP; so UIP
could equivalently be phrased as ``all types are sets''.

Two types are \emph{homotopy equivalent} if we have $f:A\to B$ and $g:B\to A$ such that $g\circ f = 1_A$ and $f\circ g = 1_B$.
However, the type of such data
\begin{equation}
  \textstyle\sum_{f:A\to B} \sum_{g:B\to A} (g\circ f = 1_A)\times (f\circ g = 1_B)\label{eq:qinv}
\end{equation}
is not a correct definition of \emph{the type of homotopy equivalences}. 
(It is correct if $A$ and $B$ are sets, in which case we generally say \emph{isomorphism} or \emph{bijection} rather than \emph{equivalence}.)
The problem is that given $f:A\to B$, the rest of~\eqref{eq:qinv}:
\begin{equation}
  \textstyle\sum_{g:B\to A} (g\circ f = 1_A)\times (f\circ g = 1_B)\label{eq:isquinv}
\end{equation}
may not be a proposition, whereas we want ``being an equivalence'' to be a mere \emph{property} of a morphism. 
For instance, if $f$ is the identity map of the homotopical circle $S^1$ (see \cref{sec:synth-homot-theory}), then~\eqref{eq:isquinv} is equivalent to $\mathbb{Z}$.
Thus if we took~\eqref{eq:qinv} as our definition of equivalence, there would be ``infinitely many self-equivalences of $S^1$'', which is not correct: up to homotopy there should be only two.

Many equivalent ways to correct~\eqref{eq:qinv} are now known; here are a few: 
\begin{gather}
  \textstyle\sum_{f:A\to B} \sum_{g:B\to A} \sum_{h:B\to A} (g\circ f = 1_A)\times (f\circ h = 1_B) \label{eq:hiso}\\
  \textstyle\sum_{f:A\to B} \sum_{g:B\to A} \sum_{\eta:g\circ f = 1_A}\sum_{\epsilon: f\circ g = 1_B} f(\eta) = \epsilon_f.\label{eq:hae}\\
  \textstyle\sum_{f:A\to B} \brck{\sum_{g:B\to A} (g\circ f = 1_A)\times (f\circ g = 1_B)}\\
  \textstyle\sum_{f:A\to B} \prod_{y:B} \Big(\big(\sum_{x:A} f(x)=y\big) \text{ is contractible}\Big)
\end{gather}
Each of these admits maps back and forth from~\eqref{eq:qinv},
while the data after the $\sum_{f:A\to B}$ form a proposition.
We can think of these as building contractible cell complexes.
For instance, in~\eqref{eq:hiso} we glue on two 1-cells $g,h$ each with a contracting 2-cell, giving a contractible space; whereas in~\eqref{eq:qinv} we glue on \emph{one} 1-cell with \emph{two} contracting 2-cells, giving a non-contractible 2-sphere.
And in~\eqref{eq:hae} we add to~\eqref{eq:qinv} a 3-cell filler, getting a contractible 3-ball.

When making definitions of this sort, we generally think of types as \oo-groupoids (or homotopy spaces), just as in \cref{sec:constructive-logic} we thought of types as sets.
However, since type theory presents an initial structured \io-category, these definitions can also be interpreted in any structured \io-category, yielding ``classifying spaces'' for $n$-types and equivalences.
For example, given $\Gamma\types A\ty$ and $\Gamma\types B\ty$, if $\Gamma\types\equiv(A,B)\ty$ denotes the type of equivalences (with any of the corrected definitions above), then the object $\m{\equiv(A,B)} \to \m\Gamma$ of the slice category has the universal property that for any map $f:X\to\m\Gamma$, lifts of $f$ to $\m{\equiv(A,B)}$ are equivalent to homotopy equivalences $f^*\m A \simeq f^* \m B$ over $X$.
(See~\cite[\S3]{klv:ssetmodel} or~\cite[\S4]{shulman:elreedy}.)
In other words, $\m{\equiv(A,B)}$ is a ``classifying space for equivalences between $\m A$ and $\m B$''.

\subsection{Extensionality and Univalence}
\label{sec:voev-univ-axiom}

In \cref{sec:synthetic-topology-1} we ``actualized'' the potential spatial structure of types by adding axioms such as a dominance or a type of smooth reals.
Similarly, we can add axioms ensuring that there really are types with higher groupoid structure --- i.e.\ that not all types are sets, or that UIP fails.
The only serious contender for such an axiom at present is Voevodsky's \emph{univalence axiom}. 

To explain univalence, let us return to the type constructors discussed in \cref{sec:rules,sec:universes}.
We claimed that the rules for coproduct types, function types, subobject classifiers, and so on express their desired categorical universal properties.
For most types this is literally true, but there are a couple of cases%
\footnote{What these cases have in common is that they are ``mapping in'' universal properties.
  ``Mapping out'' universal properties, like that for coproducts, can be expressed more powerfully in type theory using a dependent output such as in~\eqref{eq:coprod-depelim}, enabling us to derive their full universal property from the basic rules.
  (Of course, the cartesian product also has a ``mapping in'' universal property, but it doesn't have this problem; formally this is because the classifying category is better described as a sort of ``cartesian multicategory'' in which the cartesian product also has a ``mapping out'' universal property.)}
in which something is missing from our discussion so far.

Firstly, the universal property of an exponential object requires that any map $X\times \m A \to \m B$ factors through a \emph{unique} map $X\to {\m B}^{\m A}$.
For $\m{A\to B}$ to be ${\m B}^{\m A}$, therefore, requires that if $\Gamma\types f:A\to B$ and $\Gamma\types g:A\to B$ and $\Gamma,x:A \types h: f(x)=g(x)$, then also $\Gamma \types e:f=g$ (because elements of the identity type induce equalities of morphisms in $\Ctx$, or homotopies in $\cCtx$).
This is not derivable from the rules in \cref{fig:function-types}; it is an extra axiom called \emph{function extensionality}.
Informally, it says that two functions are equal if they take equal values.
(Dependent function types 
require a similar axiom.)

Secondly, the universal property of a subobject classifier requires that a mono $M \to \m\Gamma$ is classified by a \emph{unique} map $\m\Gamma\to\Omega$; or equivalently, two maps $\m\Gamma\to\Omega$ classifying the same subobject of $\m\Gamma$ are equal.
Here ``the same'' means \emph{isomorphism} in $\Ctx/\m\Gamma$; a classifying map only determines a mono up to isomorphism anyway.
Type-theoretically, this means that if $\Gamma\types P:\Omega$ and $\Gamma\types Q:\Omega$ and $\Gamma\types h:\equiv(P,Q)$, then $\Gamma\types e:P=Q$.
This is not derivable from the rules in \cref{sec:universes}; it is an extra axiom called \emph{propositional extensionality}.
(Note that two propositions are equivalent as soon as each implies the other.)

When we try to generalize propositional extensionality for $\Omega$ to a statement about type universes $\U$, things become more subtle.
We still expect classifying maps to classify only up to isomorphism --- or better, up to homotopy equivalence, which leads us towards homotopical classifying spaces.
In traditional homotopy theory (e.g.~\cite{may:csf}), \emph{homotopy classes} of maps into a classifying space correspond to \emph{homotopy equivalence classes} of fibrations over it.
But in an $(\oo,1)$-category, it is more natural to ask directly that the \oo-groupoid $\hom(X,\U)$ is equivalent to a full sub-\oo-groupoid of the slice category over $X$; this gives the notion of an \emph{object classifier}~\cite[\S6.1.6]{lurie:higher-topoi}.
Type-theoretically, the corresponding condition is that for $\Gamma\types A:\U$ and $\Gamma\types B:\U$, the type $\Gamma \types (A=B) \ty$ (i.e.\ the \oo-groupoid of homotopies between classifying maps) is \emph{equivalent} to the type $\equiv(A,B)$ of homotopy equivalences as in \cref{sec:oo-group-struct}.
More precisely, identity-elimination yields a function $\mathsf{idtoeqv}_{A,B} : (A=B) \to \equiv(A,B)$, and we should require that this map is itself an equivalence.

This axiom is due to Voevodsky, who dubbed it \emph{univalence}.%
\footnote{By analogy with function extensionality and propositional extensionality, univalence could be called \emph{typal extensionality}.
  In particular, like function extensionality and propositional extensionality, univalence is an ``extensionality'' property, meaning that ``types are determined by their behavior''.
  For this reason, it is unfortunate that the phrase ``extensional type theory'' has come to refer to type theory with UIP, which is incompatible with univalence.
  Historically, the special case of univalence when $A$ and $B$ are sets, in which case one can use~\eqref{eq:qinv} without ``correction'', was proposed by Hofmann and Streicher~\cite{hs:gpd-typethy} under the name ``universe extensionality'', but it didn't attract much attention.}
Univalence clearly implies propositional extensionality, 
while Voevodsky showed~\cite{vv:unimath} that it also implies function extensionality; see e.g.~\cite[\S4.9]{hottbook}.
Univalence also does indeed ensure that not all types are sets (i.e.\ homotopy 0-types).
For instance, if $B:\U$ has a nontrivial automorphism, such as $B=\unit+\unit$, then $\equiv(B,B)$ is not a proposition.
Hence neither is the equality type $B=B$ in $\U$, so $\U$ is not a set.
More generally, with a hierarchy of universes $\U_n$ with $\U_n : \U_{n+1}$, each $\U_n$ is not an $n$-type~\cite{ks:u-not-ntype}.

In particular, for models in the category of sets, or more generally in any 1-category, univalence must be \emph{false}.
For instance, any ``Grothendieck universe'' in ZFC set theory can be used as a type-theoretic universe $\U$ in $\mathbf{Set}$; but it is not univalent, since it would be a set (a 0-type), whereas the above argument shows no univalent universe containing a 2-element set can be a set.

Formally, univalence is an axiom like UIP and LEM that cuts down our collection of models, only now in a way that \emph{excludes} all 1-categories.
Just as the topological axioms from \cref{sec:synthetic-topology-1} are incompatible with LEM, univalence is incompatible with UIP.
These two oppositions are essentially independent: the topos of sets satisfies both LEM and UIP, the toposes of consequential spaces and continuous sets satisfy UIP together with topological axioms instead of LEM, and the $\io$-category of $\oo$-groupoids satisfies univalence instead of UIP but still satisfies LEM.
(
In particular, LEM does not rule out \emph{all} ``spatial'' interpretations of type theory, at least if we regard $\oo$-groupoids as a kind of ``space''.)
Finally, in \cref{sec:cohes-homot-type} we will mention some $\io$-categories that combine topological axioms with univalence, thus satisfying neither LEM nor UIP.%
\footnote{\label{fn:globalchoice}Most ``classicality'' properties such as the axiom of choice behave similarly to LEM in this way, but some very strong choice principles do conflict with univalence, such as the existence of a global ``Hilbert choice'' operator; see~\cite[\S3.2 and Exercise 3.11]{hottbook}.}

To a homotopy theorist or higher category theorist, assuming univalence instead of UIP is obviously the right move; but it can be a difficult step for those used to thinking of types as sets.
However, univalence can also be motivated from purely type-theoretic considerations, as giving a ``correct'' answer to the question ``what are the identity types of a universe?'', just as function extensionality answers ``what are the identity types of a function type?''.
And from a philosophical point of view, univalence says that all properties of types are invariant under equivalence, since we can make any equivalence into an equality and apply transport;
thus it expresses a strong ``structural'' nature of type theory~\cite{awodey:struc-invar-ua,corfield:structure-of-a}, in contrast to ZFC-style set theory.

For the homotopy theorist, univalence is one place where we start to see the advantage of the type-theoretic syntax.
Inside type theory, the ``elements'' of a universe type $\U$ \emph{are themselves types}, in contrast to the classical construction of classifying spaces whose ``points'' often lack a meaning directly connected to the things being classified.
This enables us to define other classifying spaces and operations between them in a very intuitive way.
For instance, if $G$ is a group (meaning a \emph{set}, a 0-type, with a group structure), we can define its classifying space to be ``the type of free transitive $G$-sets'':
\[ \mathbf{B}G = \textstyle\sum_{A:\U} \sum_{a:G\times A\to A} (\text{$A$ is a set and $a$ is a free transitive action}) \]
That is, an element of $\mathbf{B}G$ is a tuple $(A,a,\dots)$ consisting of a type, an action of $G$ on that type, and witnesses of the truth of the necessary axioms.
It turns out that $\mathbf{B}G$ is a connected 1-type with $\Omega(\mathbf{B}G) = G$.
If $G$ is abelian, we can define an operation $\mathbf{B}G \times \mathbf{B}G \to \mathbf{B}G$ by taking the ``tensor product'' of $G$-sets, and so on.

This definition of $\mathbf{B}G$ also immediately defines the objects it classifies: a ``torsor'' over a type $X$ is just a function $X \to \mathbf{B}G$.
The first component of such a function is a map $X\to \U$, corresponding to a dependent type $x:X \types A\ty$, and hence a map $\m A \to \m X$.
The rest of the classifying map equips this with the usual structure of a torsor over $X$.

In fact, \emph{any} definition of a structure in type theory automatically defines the classifying space for such structures, and therefore also automatically the corresponding notion of ``bundle of structures''.
For instance, a group can be considered a tuple $(G,e,m,\dots)$ of a set, an identity, a multiplication, and proofs of the axioms, giving a definition of ``the type of groups'':
\[ \mathsf{Group} = \textstyle \sum_{G:\U}\sum_{e:G} \sum_{m:G\times G\to G} (\text{$G$ is a set and $(m,e)$ is a group structure}). \]
We then automatically obtain a notion of a ``family of groups'', namely a function $X\to \mathsf{Group}$.
This turns out to correspond precisely to a \emph{local system} of groups in the sense of classical homotopy theory.
Similarly, we can define a \emph{spectrum} to be a sequence of pointed types $(X_n,x_n)$ each of which is the loop space of the next; thus ``the type of spectra'' is
\[ \mathsf{Spectrum} = \textstyle\sum_{X:\N\to \U} \sum_{x:\prod_{n:\N} X_n} \prod_{n:\N} \big((X_n,x_n) = \Omega(X_{n+1},x_{n+1})\big). \]
This yields automatically a notion of ``parametrized spectrum'', namely a function $X\to \mathsf{Spectra}$.
The homotopy groups of a spectrum are functions $\pi_n:\mathsf{Spectrum} \to \mathsf{AbGroup}$, while the Eilenberg--Mac Lane construction is a function $H : \mathsf{AbGroup} \to \mathsf{Spectrum}$; these then act by simple composition to relate parametrized spectra and local systems.
Thus, type theory automatically handles generalizations to ``parametrized spaces'', which in classical homotopy theory and category theory have to be done by hand.

\subsection{Higher inductive types}
\label{sec:hits}

In \cref{sec:rules} 
we mentioned type constructors corresponding to coproducts, products, exponentials, initial and terminal objects, diagonals, and natural numbers objects.
Combining dependent sum types with the identity type yields all finite limits; for instance, the pullback of $f:A\to C$ and $g:B\to C$ is
\[ \tsm_{x:A} \tsm_{y:B} (f(x)=g(y)). \]
Moreover, with the natural numbers type we can express certain infinite limits, e.g.\ the limit of a sequence $\cdots \xrightarrow{p_2} A_2 \xrightarrow{p_1} A_1 \xrightarrow{p_0} A_0$ is
\[ \tsm_{f:\prod_{n:\N} A_n} \textstyle\prod_{n:\N} p_n(f(n+1)) = f(n). \]

\begin{figure}
  \centering
  \begin{mathpar}
  \inferrule{\Gamma \types f:A\to B \\ \Gamma \types g:A\to B}{\Gamma\types \coeq(f,g)\ty}\and
  \inferrule{\Gamma\types N:B}{\Gamma\types \classof{N}:\coeq(f,g)}\and
  \inferrule{\Gamma \types M:A}{\Gamma \types \ceq(M) : \classof{f(M)} = \classof{g(M)}}\and
  \inferrule{\Gamma,z:\coeq(f,g)\types C\ty \\ \Gamma,y:B \types c_B : C[\classof{y}/z] \\ \Gamma,x:A \types c_A : \trans{\ceq(x)}{c_B[f(x)/y]} = c_B[g(x)/y] \\ \Gamma \types P:\coeq(f,g)}{\Gamma\types \cind(C,c_B,c_A,P) : C[P/z]}\and
    \inferrule{\vdots}{\Gamma\types \cind(C,c_B,c_A,\classof{N}) = c_B[N/y]}
  \and
  \inferrule{\vdots}{\Gamma\types \ap{\cind(C,c_B,c_A)}(\ceq(M)) = c_A[M/x]}
\end{mathpar}
\caption{The rules for coequalizer types}
\label{fig:coeq-types}
\end{figure}

However, to represent \emph{colimits} other than coproducts 
we need new type constructors.
For instance, the rules for the coequalizer type $\coeq (f,g)$ are shown in \cref{fig:coeq-types}.
The first is formation: any $f,g:A\to B$ have a coequalizer.
The next two are introduction: there is a map $B\to \coeq(f,g)$, and the two composites $A \rightrightarrows B \to \coeq(f,g)$ are equal.
The third is the elimination rule, which is analogous to the case analysis rule~\eqref{eq:coprod-depelim}.
To understand this, consider first the simpler version analogous to~\eqref{eq:coprod-nondepelim}, where $C$ does not depend on $\coeq(f,g)$:
\[  \inferrule{\Gamma\types C\ty \\ \Gamma,y:B \types c_B : C \\ \Gamma,x:A \types c_A : c_B[f(x)/y] = c_B[g(x)/y] \\ \Gamma \types P:\coeq(f,g)}{\Gamma\types \cind(C,c_B,c_A,P) : C}. \]
This expresses the existence part of the universal property of a coequalizer: given a map $B\to C$ such that the composites $A\rightrightarrows B \to C$ are equal, there is an induced map $\coeq(f,g) \to C$.

As with coproducts, the more general version in \cref{fig:coeq-types} also implies the uniqueness part of the universal property.
It contains one new aspect: if $C$ depends on $\coeq(f,g)$, then $c_B[f(x)/y]$ and $c_B[g(x)/y]$ have different types $C[f(x)/y]$ and $C[g(x)/y]$, so we cannot write ``$c_B[f(x)/y] = c_B[g(x)/y]$''.
But we have $\ceq(x) : f(x)=g(x)$, so the types $C[f(x)/y]$ and $C[g(x)/y]$ ought to be ``the same''; but formally we need to ``transport'' $c_B[f(x)/y]$ along $\ceq(x)$ (using identity-elimination) to get an element of $C[g(x)/y]$ that we can compare to $c_B[g(x)/y]$.
This is what the notation $\trans{\ceq(x)}{c_B[f(x)/y]}$ means.

Finally, the last two rules (in which I have omitted the premises for brevity) are the computation rules.
The first says that when a map $\coeq(f,g)\to C$ is induced by the universal property, the composite $B \to \coeq(f,g) \to C$ is indeed the original map $B\to C$.
The second says similarly that the ``induced equality'' between the composites $A\rightrightarrows B \to \coeq(f,g) \to C$ is the originally given one.
(Don't worry about the notation; it's not important for us.)

If this ``equality of equalities'' sounds weird, recall that in \emph{homotopy} type theory, the type $x=y$ represents the hom-$\oo$-groupoid, and hence can have many different elements.
Thus, it makes sense to ask whether two such ``equalities'' are equal.
In fact, when we regard type theory as presenting an $\io$-category rather than a 1-category, the type $\coeq(f,g)$ represents an $\oo$-categorical coequalizer, a.k.a.\ homotopy coequalizer.
From this we can build all finite (homotopy) colimits, as in~\cite[Corollary 4.4.2.4]{lurie:higher-topoi}.
We also obtain certain infinite colimits: e.g.\ the coproduct of a countably infinite family $A:\N\to\U$ is just $\sum_{n:\N} A_n$, and the co\-limit of a sequence $A_0 \xrightarrow{f_0} A_1 \xrightarrow{f_1} A_2 \xrightarrow{f_2} \cdots$ is the coequalizer of two maps $(\sum_{n:\N} A_n) \rightrightarrows (\sum_{n:\N} A_n)$.

The rules for $\coeq(f,g)$ do not \emph{require} any $\oo$-categorical behavior, and are perfectly consistent with UIP.
In particular, adding them to \emph{extensional} MLTT yields a type theory for locally cartesian closed categories with finite colimits.
Nevertheless, types such as $\coeq(f,g)$ were not widely studied prior to the advent of homotopy type theory;
they are known as \emph{higher inductive types}. 

In general, an \textbf{inductive type} $W$ is specified by a list of \emph{constructors}, which are (possibly dependent) functions into $W$.
For instance, the coproduct $A+B$ is the inductive type specified by two constructors $\inl : A \to A+B$ and $\inr : B \to A+B$.
(The empty type $\emptyset$ is inductively specified by \emph{no} constructors.)
The constructors are the introduction rules of the resulting type, while the elimination rule says that to define a map out of the inductive type it is sufficient to specify its behavior on the constructors.

A \textbf{higher inductive type (HIT)} is similar, but the constructors can also be functions into equality types of the HIT.
For instance, $\coeq(f,g)$ is specified by two constructors $\classof{\blank} : B\to \coeq(f,g)$ and $\ceq : \tprod_{x:A} (\classof{f(x)} = \classof{g(x)})$.

The word ``inductive'' comes from the fact that in general, the type being defined is allowed to appear in the \emph{domains} of its constructors in certain limited ways.
For instance, the natural numbers are the inductive type specified by two constructors $0:\N$ (a 0-ary function) and $\succ:\N\to\N$.
Informally, this means that the elements of $\N$ are generated by applying the constructors successively any number of times; thus we have $0$, $\succ(0)$, $\succ(\succ(0))$, and so on.

When combined with higher constructors, this additional feature is quite powerful; for instance, the propositional truncation $\brck{A}$ from \cref{fig:brck} is the HIT specified by two constructors $\bproj{\blank} : A\to\brck{A}$ and $\mathsf{tprp}:\tprod_{x,y:\brck{A}} (x=y)$.
We can similarly construct an \textbf{$n$-truncation} that is the universal map $A\to \trunc n A$ into a homotopy $n$-type (i.e.\ its $n^{\mathrm{th}}$ Postnikov section).
In particular, the $0$-truncation $\trunc 0 A$ is the ``set of connected components''. 

``Recursive'' HITs of this sort can also be used to construct more exotic objects, such as homotopical localizations.
Given a map $f:S\to T$, we say that a type $A$ is \textbf{$f$-local} if the map $(\blank\circ f):(T\to A)\to (S\to A)$ is an equivalence.
The \textbf{$f$-localization} is the universal map from a type $X$ into an $f$-local type $L_f X$.
In classical homotopy theory, constructing localizations in general requires a fairly elaborate transfinite composition.
But in homotopy type theory, we can simply define $L_f X$ to be the HIT generated by the following constructors:
\begin{compactitem}
\item A map $\eta : X\to L_f X$.
\item For each $g:S\to L_f X$ and $t:T$, an element $\mathsf{ext}(g,t):L_f X$.
\item For each $g:S\to L_f X$ and $s:S$, an equality $\mathsf{ext}(g,f(s)) = g(s)$.
\item For each $g:S\to L_f X$ and $t:T$, an element $\mathsf{ext}'(g,t):L_f X$.
\item For each $h:T\to L_f X$ and $t:T$, an equality $\mathsf{ext}'(h\circ f,t)=h(t)$.
\end{compactitem}
The last four constructors combine to lift $(\blank\circ f)$ to an element of~\eqref{eq:hiso}.  (This is why we have both $\mathsf{ext}$ and $\mathsf{ext}'$; if we collapsed them into one we would only get~\eqref{eq:qinv}.)
This is one example of how homotopy type theory gives a ``direct'' way of working with objects and constructions that in classical homotopy theory must be laboriously built up out of sets.
For more examples and theory of higher inductive types and their applications, see~\cite[Chapter 6]{hottbook}.

\subsection{Synthetic homotopy theory}
\label{sec:synth-homot-theory}

With HITs we can define many familiar spaces from classical homotopy theory.
For instance, in the $\oo$-category of $\oo$-groupoids, the circle $\hocirc$ is the homotopy coequalizer of $\unit \rightrightarrows \unit$; thus we expect the corresponding coequalizer type to behave like an ``internal $\hocirc$'' in homotopy type theory.
It is equivalently the HIT generated by two constructors $\base:\hocirc$ and $\lloop :\base=\base$; its elimination rule (universal property) says roughly that to give a map $\hocirc \to C$ is equivalent to giving a point $c:C$ and a loop $l:c=c$.

Since HITs are consistent with UIP, this ``circle'' may not behave as expected: in a 1-category, the coequalizer of $\unit \rightrightarrows \unit$ is just $\unit$.
But if we also assume univalence, type theory becomes a powerful tool for working directly with $\oo$-groupoids such as $\hocirc$.
By the universal property of $\hocirc$, to give a dependent type $C:\hocirc \to \U$ we must give a type $B:\U$ and an equality $B=B$; but by univalence the latter is the same as an \emph{autoequivalence} of $B$.
For instance, if $B$ is $\Z$, we can use the autoequivalence ``$+1$''; the resulting dependent type is then a version of the \emph{universal cover} of $\hocirc$.
With a little extra work~\cite{ls:pi1s1}, we can adapt the classical calculation of $\pi_1(\hocirc)$ to show, in type theory, that $\Omega \hocirc \simeq \Z$.

This is the first theorem of a growing field known as \emph{synthetic homotopy theory}, more of which can be found in~\cite[Chapter 8]{hottbook} and recent work such as~\cite{lf:emspaces,ffll:blakers-massey,brunerie:thesis}.
Just as in the \emph{synthetic topology} of \cref{sec:synthetic-topology-1} the types come automatically with topological structure, which we can then study ``synthetically'' rather than breaking it down into a set equipped with a topology, in synthetic homotopy theory the types come automatically with homotopical or $\oo$-groupoid structure, which we can then study synthetically rather than breaking it down into any explicit definition of an $\oo$-groupoid.
Thus it is a ``model-independent'' language for homotopy theory, avoiding the need to choose (say) topological spaces or simplicial sets as a definition of ``$\oo$-groupoid.''

It is too early to say how useful this will be to classical homotopy theory.
In its very short existence so far, synthetic homotopy theory has not led to proofs of any new theorems, but it has shown an impressive ability to produce new proofs of old theorems:
as of this writing, synthetic homotopy theorists have calculated $\pi_n(S^n)=\Z$, $\pi_k(S^n)=0$ for $k<n$, $\pi_3(S^2)=\Z$, and $\pi_4(S^3)=\Z/2\Z$, and proven numerous foundational results such as the Freudenthal suspension theorem, the Blakers--Massey connectivity theorem, and the Serre spectral sequence.

More importantly, the theorems of synthetic homotopy theory are more general than those of classical homotopy theory, because (modulo subtleties to be mentioned in \cref{sec:semantics}) they apply in any well-behaved $\io$-category, 
including any \emph{$\oo$-topos}~\cite{rezk:homotopy-toposes,lurie:higher-topoi}.
(Some $\oo$-toposes of interest to classical homotopy theorists include equivariant and parametrized homotopy theory.)
A particularly interesting example is the Blakers--Massey theorem, for which no purely homotopical proof applicable to $\oo$-toposes was known prior to the synthetic one~\cite{ffll:blakers-massey}; the latter has now been translated back into categorical language~\cite{rezk:hott-blakersmassey}.

Finally, synthetic homotopy theory gives a new way to think about the ``homotopy hypothesis'' of Grothendieck~\cite{baez:homotopy-hypothesis} that $\oo$-groupoids describe the homotopy theory of spaces.
Rather than looking for an equivalence between some notions of $\oo$-groupoid and space, we have a \emph{synthetic} theory of $\oo$-groupoids that is \emph{modeled by} classical homotopy spaces --- but also other things.
(In fact, Brunerie has observed that the $\oo$-groupoid structure of types in homotopy type theory looks almost exactly as it was envisioned by Grothendieck~\cite{maltsiniotis:groth-oocat}, rather than like any of the definitions of \oo-groupoid used more commonly today.)
In \cref{sec:cohes-homot-type} I will sketch a particular context in which this extra generality is useful.

\begin{subappendices}

\subsection{The classifying \io-category}
\label{sec:semantics}

In this appendix to \cref{sec:homotopy-type-theory} I will describe the ``classifying \io-category'' informally, then give a precise definition of it, and end with some remarks about the current state of knowledge as regards its freeness.
This appendix and \cref{sec:equality} are provided to satisfy the curious reader, but can be skipped without consequence.

Let $\fT$ be an intensional type theory.
We define the objects and morphisms of its classifying $\io$-category $\cCtxT$ just as we did for the classifying 1-category $\CtxT$ in \cref{sec:type-term-judgments}.
However, we do \emph{not} quotient the morphisms by terms in the equality type.
Instead we will use those to define the 2-morphisms, as well as 3-morphisms, 4-morphisms, and so on.

The idea is to generalize the representation of equalities using diagonals to a characterization of 2-morphisms.
Given morphisms $f,g:A\to B$ in an $\io$-category, their ``equalizer'' is a morphism $e:E\to A$ equipped with a 2-morphism $f e \cong g e$ that is ``universal'' among such 2-morphisms.
In particular, to give a 2-morphism $f\cong g$ is equivalent to giving a section of $e$ (that is, a morphism $s:A\to E$ and a 2-morphism $e s \cong 1_A$).
As in the 1-categorical case, this equalizer can be constructed, up to equivalence, as the pullback of the diagonal $B\to B\times B$ along $(f,g):A\to B\times B$.
Thus, assuming that the identity type of $B$ still presents the diagonal (up to the appropriate sort of $\io$-categorical equivalence), and substitution still presents pullback, 2-morphisms $f\cong g$ should be equivalent to terms of the form
\begin{equation}
  x:A \types p(x):f(x)=_B g(x).\label{eq:ctx-2mor}
\end{equation}
Hence we simply \emph{define} a 2-morphism in $\cCtxT$ to be a term of this sort.
Similarly, we define a 3-morphism $p\cong q$ to be a term in an iterated identity type $x:A \types h(x):p(x)=_{(f(x)=_B g(x))} q(x)$, and so on.

To make this precise, we need to choose a method of presenting \io-categories.
In principle there are many options, but at present the method of choice for defining $\cCtxT$ is to use \emph{fibration categories}~\cite{brown:ahtgsc}.
A fibration category (or ``category of fibrant objects'') is a 1-category with two classes of morphisms called \emph{weak equivalences} and \emph{fibrations}, satisfying certain axioms, e.g.\ pullbacks of fibrations exist and preserve weak equivalences.
The most important axiom is that every diagonal $A \to A\times A$ factors as a weak equivalence followed by a fibration, with the intermediate object called a \emph{path object} $P A$ for $A$.

Of course, this is an abstraction of a common situation from homotopy theory: fibrations of topological spaces, Kan simplicial sets, or chain complexes (and more generally the fibrations between fibrant objects in any Quillen model category) all have these properties.
Generalizing these examples, in any fibration category we define a \emph{homotopy} between $f,g:A\rightrightarrows B$ to be a lift of $(f,g):A \to B\times B$ to a map $A \to P B$.
We can similarly define higher homotopies and thereby construct a more explicit notion of \io-category (such as a quasicategory), although the combinatorics are somewhat involved; see~\cite{szumilo:hothy-ccplhot}.

Now, if in the definition of $\CtxT$ from \cref{sec:type-term-judgments} we omit the quotient of morphisms,\footnote{Technically, we replace it with a different quotient; see \cref{sec:equality}.} we obtain a fibration category $\sCtxT$.
Its fibrations are the composites of projections $\m{\Gamma,x:A} \to \m\Gamma$, its weak equivalences are the homotopy equivalences defined in \cref{sec:oo-group-struct}, and its path objects are the identity types $P\m A = \m{x:A, y:A, p:x=y}$.
(The fact that identity types satisfy the axioms of path objects was one of the central insights of Awodey and Warren~\cite{aw:htpy-idtype,warren:thesis}.)
With this definition, homotopies in the fibration-category sense correspond bijectively to terms of the form~\eqref{eq:ctx-2mor}: 
The former are lifts as on the left below, whereas the latter are sections as on the right.
\[\xymatrix{ & \Id_B \ar[d] \\ \m A \ar[r]_-{(f,g)} \ar@{.>}[ur] & \m B \times \m B }\qquad
\xymatrix{\m{x:A, p:f(x)=g(x)} \ar[d] \ar[r] \pullback & \Id_B \ar[d] \\ \m A \ar[r]_-{(f,g)} \ar@{.>}@/^4mm/[u] & \m B \times \m B }
\]

Thus, we may define $\cCtxT$ to be the \io-category presented by $\sCtxT$.
Here we see the second advantage of syntax mentioned in \cref{sec:syntax}: giving a presentation of a free object (here, an \io-category) that is actually stricter (here, a fibration category) than one would expect from only its universal property.

However, although this $\cCtxT$ has some of the expected structure~\cite{kapulkin:lccqcat-tt,kl:hot-tt},
no one has yet proven its \io-categorical freeness. 
Instead, to interpret type theory in \io-categories, we use the fact that $\sCtxT$ is free in a category of structured fibration categories.
The latter have various names like
``contextual categories''~\cite{cartmell:gatcc}, ``comprehension categories''~\cite{jacobs:compr-cat}, ``categories with families''~\cite{dybjer:internal-tt}, ``categories with attributes''~\cite{cartmell:gatcc}, ``display map categories''~\cite[\S8.3]{taylor:pracfdn}, ``type-theoretic fibration categories''~\cite{shulman:invdia}, ``tribes''~\cite{joyal:tribes}, ``C-systems''~\cite{voevodsky:subquot-csys}, and so on.
Although this approach has proven more tractable, it is still quite difficult, for two reasons.
One is that, as mentioned for the 1-categorical case in \cref{sec:rules}, complete proofs of the freeness of $\sCtxT$ have been given only for a few particular type theories~\cite{streicher:semtt}.
Everyone expects these proofs to generalize to all other type theories, but actually writing down such a generalization, and in a useful amount of generality, is a current research problem.

Another difficulty is that this approach incurs a new proof obligation.
In principle, a type theory $\fT$ should be interpreted in an \io-category $\C$ by means of the unique functor $\cCtxT \to \C$ determined by the universal property of $\cCtxT$.
If we stick with the 1-categorical universal property of $\sCtxT$, then to interpret $\fT$ in $\C$ we need to also present $\C$ by a fibration category of the appropriate sort.
This is a sort of ``coherence theorem'' for structured $\io$-categories --- which, again, is known in some particular cases, but a fully general version of which is a current research problem; the state of the art includes~\cite{klv:ssetmodel,ak:htmtt,gb:topsimpid,shulman:invdia,shulman:elreedy,shulman:eiuniv,stekelenburg:modestkan,gk:univlcc,kapulkin:lccqcat-tt,lw:localuniv}.
(Part of this coherence theorem is showing that pullbacks of fibrations can be made strictly functorial and preserve all the type operations strictly, which is nontrivial even for 1-categories~\cite{curien:subst,hofmann:ttinlccc,hofmann:ssdts,cd:lccc-tt}.)

I have chosen not to dwell on these issues because I have faith that they will eventually be resolved.
Instead I want to focus on the picture that such a resolution will make possible (and which is \emph{substantially} achievable even with current technology).
Thus one might call this chapter a ``programme'' for homotopy type theory and its higher-categorical semantics.
In \cref{sec:equality} I will briefly discuss another technical detail; in \cref{sec:cohes-homot-type} we will return to the programme.

\subsection{Judgmental Equality}
\label{sec:equality}

In \cref{sec:syntax} we described both the ``tautological'' and the ``reduced-words'' presentation of a free group using ``rules'' in the style of type theory. 
For the reduced-words description, this is the end of the definition; but for the tautological description, we need to describe the equivalence relation to quotient by.
This can also be defined inductively by the rules shown in \cref{fig:fg-taut-jdeq}, which essentially say that it is the smallest equivalence relation imposing the group axioms and compatible with the operations.
We also remarked that there is an algorithm for ``reducing'' any word from the tautological presentation, so that two terms are related by $\jdeq$ precisely when they reduce to the same result.
Finally, in \cref{sec:compute} we mentioned that type theory includes an analogous ``reduction algorithm'' making it into a general-purpose programming language.

\begin{figure}
  \centering
\begin{mathpar}
  \inferrule{X\elt}{X\jdeq X}\and
  \inferrule{X\jdeq Y}{Y\jdeq X}\and
  \inferrule{X\jdeq Y \\ Y\jdeq Z}{X\jdeq Z}\and
  \inferrule{X\jdeq X' \\ Y\jdeq Y'}{(X Y) \jdeq (X' Y')}\and
  \inferrule{X \jdeq Y}{X^{-1} \jdeq Y^{-1}}\and
  \inferrule{X\elt\\ Y\elt \\ Z\elt}{(X (Y Z))\jdeq ((X Y) Z)}\and
  \inferrule{X\elt}{(X e) \jdeq X}\and
  \inferrule{X\elt}{(e X) \jdeq X}\and
  \inferrule{X\elt}{(X X^{-1}) \jdeq e}\and
  \inferrule{X\elt}{(X^{-1} X) \jdeq e}\and
\end{mathpar}
\caption{Equality rules for free groups}
\label{fig:fg-taut-jdeq}
\end{figure}

Taken together, these remarks suggest that there should be two forms of type theory, one involving an equivalence relation $\jdeq$ and one not, with a ``reduction algorithm'' mapping the first to the second. 
This is more or less correct, but it turns out to be quite fiddly to describe the second type theory without reference to the first.
It is sometimes possible~\cite{hl:lf}, but more common is to describe only a type theory involving $\jdeq$, with the reduction algorithm an endofunction of its terms, and then \emph{define} the ``canonical forms'' to be those that are ``fully reduced''.
This is also more flexible, since we can add new $\jdeq$ axioms without knowing whether there is a corresponding reduction algorithm that terminates at a canonical form (or even knowing that there isn't!).
The relation $\jdeq$ is known as \emph{judgmental equality} or \emph{definitional equality} or \emph{substitutional equality}.\footnote{Technically, these three terms have slightly different meanings, but in the most common type theories they all turn out to refer to the same thing.}

Just as for free groups, when defining the corresponding free object we have to quotient by the relation $\jdeq$.
For the classifying 1-category $\CtxT$, this quotient is included in the quotient by terms of the identity type.
But for the fibration category $\sCtxT$, where we omitted the latter quotient, we do still have to impose a quotient by judgmental equality --- or, if our type theory has a terminating reduction algorithm (the technical term is ``strongly normalizing''), use only the canonical forms to represent objects and morphisms.

The puzzling thing, of course, is how this equality $\jdeq$ is related to the equality \emph{type} $x:A,y:A \types (x=y)\ty$.
Formally, the difference between these ``two equalities'' is analogous to the difference between the variables $x:A$ occurring in a context and the ``meta-variables'' such as $\Gamma$ that we use in describing the operations of the theory.
Any inductive definition uses ``meta-variables'' and can have an inductively defined equivalence relation; type theory is special because \emph{internal} to the theory there are also notions of ``variable'' and ``equality''.
The identity \emph{type} is defined by a universal property, just like most other types; whereas judgmental equality, like the equivalence relation on words in a free group, is inductively defined as the smallest equivalence relation imposing the desired axioms (the computation rules from \cref{sec:rules}, which we denoted with $\jdeq$ for this very reason) and respected by all the other judgments.
The latter condition means we have additional rules such as:
\begin{mathpar}
  \inferrule{\Gamma \types a:A}{\Gamma\types a\jdeq a}\and
  \inferrule{\Gamma \types a\jdeq b}{\Gamma\types b\jdeq a}\and
  \inferrule{\Gamma \types a\jdeq b \\ \Gamma\types b\jdeq c}{\Gamma\types a\jdeq c}\and
\end{mathpar}
\begin{equation}
  \label{eq:coe}
  \inferrule{\Gamma \types a:A \\ \Gamma\types A\jdeq B}{\Gamma\types a:B}
\end{equation}

This formal description, however, does not really explain \emph{why} we need two equalities, or what they mean intuitively.
To start with, it cannot be emphasized strongly enough that \emph{it is the identity type that represents mathematical equality}.
Equality in mathematics is a proposition, and in particular something that can be \emph{hypothesized} and \emph{proven} or \emph{disproven}.
Judgmental equality cannot be hypothesized (added to a context), nor can it be proven (inhabited by a term) or disproven (we cannot even state internally a ``negation'' of judgmental equality).
In its simplest form, judgmental equality is simply the algorithmic process of expanding definitions (hence the name ``definitional equality''): for instance, the function $\lambda x.x^2$ is \emph{by definition} the function that squares its argument, so $(\lambda x.x^2)(y+1)$ is \emph{by definition} equal to $(y+1)^2$.
But even the simplest equalities with mathematical content, such as the theorem that $x+y=y+x$ for $x,y:\N$, are not a mere matter of expanding definitions but require proof.

What, then, \emph{can} we do with judgmental equality?
The main property it has that the identity type doesn't is~\eqref{eq:coe}: given $a:A$ and $A\jdeq B$, the \emph{same term} $a$ is also an element of $B$ (hence the name ``substitutional equality'').
In particular, if $a\jdeq b$, then $(a=a)\jdeq (a=b)$, so that $\refl_a : a=b$;
thus judgmental equality implies mathematical equality.
By contrast, given $a:A$ and a mathematical equality $e:A=B$, it is possible to obtain a term of $B$, but that term is not syntactically equal to $a$; instead it is $e_*(a)$, involving the transport operation.

This need for explicit transports is somewhat annoying, so it is tempting to eliminate it by collapsing the two equalities with a \emph{reflection rule}
\[ \inferrule{\Gamma \types p:a=b}{\Gamma\types a\jdeq b}. \]
Unfortunately, this makes it impossible to detect $\jdeq$ using a reduction algorithm, since questions of mathematical equality cannot be decided algorithmically.
The reflection rule also turns out to imply UIP, which is 
a dealbreaker if we want to talk about \io-categories.
The \io-categorical point of view also makes clear why we need to notate $e$ in $e_*(a)$: since the type $A=B$ is (by univalence) the type of equivalences from $A$ to $B$, it could have many \emph{different} elements, so that $e_*(a)$ really does depend on the choice of $e$.

One might then be tempted to go to the other extreme and try to eliminate judgmental equality entirely.
We could in principle express all the computation rules from \cref{sec:rules} using elements of identity types rather than judgmental equalities.
However, the resulting proliferation of transport operations would be so extreme as to render the theory essentially unusable.
We need a happy medium, with a judgmental equality as strong as feasible but no stronger.

The intuitive meaning of judgmental equality is not entirely clear, although in some ways it is analogous to Frege's ``equality of sense'' (with mathematical equality analogous to ``equality of reference'').
Categorically, judgmental equality is analogous to the ``point-set-level'' or ``strict'' equality occurring in strict or semistrict models for higher categories, such as Quillen model categories or Gray-categories.
This finds a formal expression in the fibration-category approach to semantics, where we need a ``semistrictification'' theorem presenting any \io-category by a fibration category satisfying all the judgmental equalities of our type theory strictly.
Finding the right balance of strictness and weakness here is an active frontier of research.

\end{subappendices}

\section{Cohesive homotopy type theory}
\label{sec:cohes-homot-type}

\subsection{Spaces versus \oo-groupoids}
\label{sec:spaces-oo-groupoids}

Twice now we have encountered something called a ``circle'': in \cref{sec:synthetic-topology-1} we mentioned that $\topcirc = \setof{ (x,y):\R\times\R \mid x^2+y^2=1 }$ has the correct topology, and in \cref{sec:synth-homot-theory} we mentioned that $\hocirc = \mathrm{coeq}(\unit \rightrightarrows \unit)$ has the correct fundamental group.
However, these two types $\topcirc$ and $\hocirc$ are very different!
The first $\topcirc$ is a \emph{set} in the sense of \cref{sec:oo-group-struct}; 
whereas $\hocirc$ is definitely not, since its loop space is $\Z$.
On the other hand, $\hocirc$ is \emph{connected}, in the sense that its 0-truncation $\trunc 0 \hocirc$ is contractible; whereas since $\topcirc$ is a set, it is its own 0-truncation.

What is happening is that classical homotopy theory has led us to confuse two different things in our minds.
On one hand, a \emph{topological space} is a set with a notion of ``cohesion'' enabling us to define continuous functions and paths.
The nearby points of a continuous path are ``close'' in some sense, but they are still distinct.
On the other hand, an \emph{\oo-groupoid} has a collection of ``points'' or ``objects'', plus for each pair of objects a collection $\hom(x,y)$ of equivalences or ``ways in which $x$ and $y$ are the same'', plus for each $f,g\in \hom(x,y)$ a collection $\hom_{\hom(x,y)}(f,g)$ of ways in which $f$ and $g$ are the same, and so on.
When $\hom(x,y)$ is nonempty, $x$ and $y$ really \emph{are} the same to \oo-groupoid theory, just as in plain category theory we do not distinguish between isomorphic objects.

The relation between topological spaces and \oo-groupoids is that from any space $X$ we can \emph{construct} an \oo-groupoid $\shape X$, called its \emph{fundamental \oo-groupoid} or \emph{shape}.%
\footnote{The symbol $\shape$ is not an integral sign ($\int$) but an ``esh'', the IPA sign for a voiceless postalveolar fricative (English \textit{sh}); in \LaTeX\ it is available as \verb|\esh| with the package \texttt{phonetic}.
  An alternative notation is $\Pi_\oo$, but the letter $\Pi$ is overworked in type theory already.
  The term ``shape'' comes from ``shape theory'', which also studies generalizations of $\shape$ for ill-behaved topological spaces.}
The objects of $\shape X$ are the points of $X$, the objects of $\hom(x,y)$ are the continuous paths from $x$ to $y$, the objects of $\hom_{\hom(x,y)}(f,g)$ are the continuous endpoint-preserving homotopies from $f$ to $g$, and so on.
The confusion arises because we can study $\shape X$ without actually constructing it (or even having a definition of ``\oo-groupoid''), by working with $X$ itself and ``doing everything up to homotopy''; and historically, people did this for a long time before they even thought of defining \oo-groupoids.
Thus, algebraic topologists came to use the word ``space'' for objects that were actually being treated as \oo-groupoids.%
\footnote{Arguably, therefore, \oo-groupoids do not even belong in a book about notions of space.
However, tradition is weighty, and moreover \oo-groupoids do share some important attributes of notions of space, notably their ability to be present as ``background structure'' in the sense described in the introduction.
It is to emphasize this aspect, but also their distinctness from other notions of space, that I sometimes call them \emph{homotopy spaces}.}

Homotopy type theory forcibly brings the distinction between topological spaces and \oo-groupoids front and center, since it allows us to talk about \oo-groupoids directly in a foundational system that is also strong enough to study topological spaces.
In particular, we have the previously noted contrast between the types $\topcirc$ and $\hocirc$.
The relation between the two \emph{ought} to be that $\hocirc = \shape \topcirc$; but how are we to express this in type theory?

\subsection{Combining topology with homotopy}
\label{sec:comb-topol-with}

The description of $\shape X$ in \cref{sec:spaces-oo-groupoids} treats both topological spaces and \oo-groupoids as structures built out of sets.
However, we have seen that in type theory we can treat \emph{both} of them synthetically, suggesting that $\shape$ ought also to have a synthetic description.
This requires combining the perspectives of \cref{sec:synthetic-topology,sec:homotopy-type-theory}, obtaining a type theory in which topology and homotopy are synthetic \emph{at the same time}.
That is, we allow some types to have ``intrinsic topology'', like $\topcirc$, and also some types to have ``intrinsic homotopy'', like $\hocirc$.
It follows unavoidably that there must also be types with \emph{both} nontrivial topology and nontrivial homotopy.

At this point the advantages of a synthetic treatment become especially apparent.
Classically, to combine structures in this way we have to define a new structure called a ``topological \oo-groupoid'' or a ``topological \oo-stack'': an \oo-groupoid equipped with a ``topology'' on its objects, another on its morphisms, and so on.
If such a gadget has no nontrivial morphisms it reduces to a topological space, while if all the topologies are discrete it reduces to an ordinary \oo-groupoid.
Formally, we might define these to be \oo-stacks on one of the sites $\{\N_\oo\}$ and $\{\R^n\}_{n\in \N}$ from \cref{sec:toposes-spaces}, comprising \oo-toposes of \emph{consequential \oo-groupoids} and \emph{continuous \oo-groupoids} (or \emph{smooth \oo-groupoids}).
We would then need to develop a whole theory of such objects.

In type theory, however, we have seen that types ``potentially'' have \emph{both} topological and homotopical structure, which we can draw out by asserting axioms such as Brouwer's theorem or Voevodsky's univalence axiom.
Thus, to obtain a synthetic theory of ``topological \oo-groupoids'' is simplicity itself: we simply assert \emph{both} groups of axioms at the same time.
Of course, to model the theory in classical mathematics we still need to construct topological \oo-groupoids, but we don't need to bother about that when working \emph{in} the theory.

Schreiber's chapter argues that topological \oo-groupoids (or some enhancement thereof) are the correct context in which to formulate modern theories of physics.
(For more general discussion of stacks, see the chapter of Mestrano and Simpson.)
The type theory modeled by \io-categories of this sort is an active field of current research called \emph{cohesive homotopy type theory}~\cite{lawvere:cohesion,schreiber:dcct,ss:qgftchtt,shulman:bfp-realcohesion}.
I will conclude by sketching some of its most appealing features.

\subsection{Modalities and cohesion}
\label{sec:modalities-cohesion}

The synthetic description of $\shape$ involves a different way to access the latent topological structure of types, based on Lawvere's ideas of \emph{cohesion}~\cite{lawvere:cohesion}.
Recall from \cref{sec:why-spaces} that most ``topological'' toposes come with a string of adjunctions
\[
\xymatrix{ \text{topos of spaces} \ar[d]^\Gamma \\ \text{topos of sets}
  \ar@<10mm>[u]^{\Delta} \ar@<5mm>@{}[u]|{\dashv} \ar@<-6mm>@{}[u]|{\dashv} \ar@<-10mm>[u]_{\nabla}  }
\]
where $\Gamma$ is the underlying-set functor, $\Delta$ constructs discrete spaces, and $\nabla$ constructs indiscrete spaces, and $\Delta$ and $\nabla$ are fully faithful. 
If we restrict our attention to the topos of spaces, then what is left of this adjoint triple is a monad $\sharp = \nabla\Gamma$ that reflects into the subcategory of indiscrete types, a comonad $\flat = \Delta\Gamma$ that coreflects into the category of discrete types, and an adjunction $\flat\dashv\sharp$ such that the induced transformations $\sharp\flat\to\sharp$ and $\flat\to\flat\sharp$ are isomorphisms.

We can incorporate $\sharp$ and $\flat$ in type theory as \emph{higher modalities}.
Traditional ``modal logic'' studies propositional modalities, most famously ``it is necessary that $P$'' (usually written $\Box P$) and ``it is possible that $P$'' (usually written $\lozenge P$), but also others such as ``so-and-so knows that $P$'', ``it will always be the case that $P$'', and so on.
Since these often have monad- or comonad-like properties (e.g.\ $\Box P \to P$ and $\Box P \to \Box\Box P$), and propositions are particular types (see \cref{sec:prop-as-types}), we may consider monads and comonads acting on all types as ``higher-categorical modalities''.
I refer to type theory with $\flat$ and $\sharp$ as \emph{spatial type theory}, since it is designed for ``topological'' models such as consequential, continuous, and smooth sets or \oo-groupoids.
We will not state its rules precisely here since they involve some technicalities,
but the practical upshot is that $\flat$ and $\sharp$ behave as described above; see~\cite{shulman:bfp-realcohesion} for a more extensive discussion.

As an example, these modalities allow us to state ``discontinuous'' versions of classicality axioms, such as LEM, that \emph{do} hold in these intended models.
The usual version of LEM is $\prod_{P:\Omega} \brck{P + \neg P}$, which is false in the topological models, as discussed in \cref{sec:constructive-logic}, because a space is not generally the disjoint union of a subspace and its complement.
But $\prod_{P:\Omega} \sharp\brck{P + \neg P}$ and $\prod_{P:\flat\Omega} \brck{P + \neg P}$ \emph{are} true in these models: both equivalently express the true statement that any space is the smallest \emph{subspace} of itself containing both any given subspace and its complement.
They imply in particular that the (equivalent) subuniverses of discrete and indiscrete types satisfy ordinary LEM, and thus are a place for classical reasoning inside synthetic-topological type theory.
(Recall from \cref{sec:constructive-logic} that the indiscrete spaces are also usually the $\neg\neg$-sheaves.
This often follows automatically in spatial type theory; see~\cite{shulman:bfp-realcohesion}.)

Now, in many cases the functor $\Delta$ also has a \emph{left} adjoint, i.e.\ the discrete spaces are \emph{reflective} as well as coreflective.
A map from a space $A$ into a discrete set $\Delta B$ breaks $A$ up as a coproduct of one disjoint piece for each element of $B$.
Thus if $A$ is a coproduct of ``connected components'', any map $A\to \Delta B$ is uniquely determined by where each connected component goes, i.e.\ by a map $\pi_0(A)\to B$.
Thus $\pi_0$ is left adjoint to $\Delta$, or more precisely any left adjoint to $\Delta$ deserves the name $\pi_0$.
Note that this ``$\pi_0(A)$'' is not the same as the $0$-truncation $\trunc 0 A$ discussed in \cref{sec:homotopy-type-theory}; the latter treats types as \oo-groupoids while this one treats them as topological spaces.
In a moment we will see that $\pi_0(A) = \trunc 0{\shape A}$.

Such a left adjoint $\pi_0$ exists for continuous sets and smooth sets, though not for consequential spaces (because the latter contain spaces, like $\N_\oo$ that are not locally connected, hence not a coproduct of connected components).
A topos with an adjoint string $\pi_0\dashv \Delta\dashv \Gamma\dashv \nabla$ where $\Delta$ and $\nabla$ are fully faithful and $\pi_0$ preserves finite products (and perhaps more; see~\cite{lawvere:cohesion,johnstone:punctual-lc,shulman:bfp-realcohesion}) is called \emph{cohesive}.

Finally, this all works basically the same in the \oo-case: ``cohesive \oo-toposes'', such as continuous and smooth \oo-groupoids, are related to the \oo-topos of \oo-groupoids by a string of \oo-adjunctions, which can be represented by modalities in type theory.\footnote{Of course, the formal connection between cohesive \oo-toposes and cohesive type theory is at least as difficult as the ordinary case discussed in \cref{sec:semantics}; indeed the cohesive case has not yet been studied formally at all.
However, the cohesive type theory at least is fully rigorous as a formal system in its own right, with reference to \oo-toposes only for motivation.}
For intuition, a ``discrete'' cohesive \oo-groupoid is one whose topologies are discrete at all levels, i.e.\ neither its points, nor its equalities between points, etc., have any interesting topology.
It could still have interesting \oo-groupoid structure; for instance, $\hocirc$ is discrete (but $\topcirc$ is not!).

The magical thing is that for $\oo$-toposes, a left adjoint of $\Delta$ is no longer just $\pi_0$; instead, it deserves to be called the shape functor $\shape$ discussed above!
To prove this is technical (see~\cite[Proposition 4.3.32]{schreiber:dcct} or \cite[\S3]{carchedi:hotyorb}), but we can get a feel for it with examples.

First of all, by comparing universal properties, we see that (denoting a left \oo-adjoint of $\Delta$ by $\shape$) the set $\trunc 0{\shape A}$ is a reflection of $A$ into discrete \emph{sets} (i.e.\ homotopy 0-types).
Thus, the 1-categorical argument above implies that $\trunc 0{\shape A}$ deserves the name $\pi_0(A)$, which is what we expect for the shape of $A$.

Secondly, 
we have seen that $\hocirc = \mathrm{coeq}(\unit\rightrightarrows\unit)$, and since the discrete types are closed under colimits (being coreflective), $\hocirc$ is also discrete.
On the other hand we have $\topcirc = \mathrm{coeq}(\R \rightrightarrows \R)$, where one map $\R\to\R$ is the identity and the other is ``$+1$''.
Since left adjoints preserve colimits, we will have $\shape\topcirc = \hocirc$ as long as $\shape \R=\unit$.
This is true for continuous $\oo$-groupoids (an analogous fact about the smooth reals is true for smooth $\oo$-groupoids).

In fact, the discrete objects in continuous $\oo$-groupoids are essentially \emph{defined} by the property that $\shape\R=\unit$.
More specifically, a type $A$ is discrete if and only if every map $\R\to A$ is constant, or more precisely if the map $\mathsf{const}:A \to (\R\to A)$ is an equivalence.
This axiom is called \emph{real-cohesion}~\cite{shulman:bfp-realcohesion}; it immediately implies that $\shape \R=\unit$.
(The real-cohesion axiom also allows us to \emph{construct} $\shape$ as a higher inductive type, by ``localizing'' in the sense of \cref{sec:hits} at the map $\R\to \unit$.)

We can make similar arguments in other examples.
For instance, the topological 2-sphere $\topsph$ is the pushout of two open discs (each isomorphic to $\R^2$) under an open strip (isomorphic to $\topcirc \times \R$).
Thus, as long as $\shape$ preserves products, $\shape\topsph$ is the pushout of two copies of $\unit$ under $\hocirc$, i.e.\ the homotopy-theoretic \emph{suspension} of $\hocirc$, which is one definition of the homotopical 2-sphere $\hosph$.
Many familiar spaces can be presented as ``open cell complexes'' of this sort, thereby identifying their shapes with the expected discrete \oo-groupoids.

We do have to avoid the more classical ``closed cell complexes'' that glue intervals and closed discs along boundaries.
Gluing the endpoints of the unit interval $[0,1]$ in the topos of continuous sets does not produce $\topcirc$, but rather a circle with a ``speed bump'' at which any continuous path must stop for a finite amount of time before proceeding.
This problem is avoided by consequential spaces, but as remarked previously that topos fails to have $\shape$.
In fact, as discussed briefly in~\cite{ptj:topological-topos}, it seems impossible to have both closed cell complexes and $\shape$.

This description of $\shape$ enables synthetic arguments that involve both topological spaces and homotopy spaces and their relationship.
For instance, in~\cite{shulman:bfp-realcohesion} I used $\shape\topcirc = \hocirc$ to prove the Brouwer fixed point theorem synthetically.
This is a theorem about the \emph{topological} closed disc $\topdisc$ (whose boundary is $\topcirc$), but its classical proof uses a \emph{homotopical} argument, constructing a retraction $\topdisc \to \topcirc$ which is impossible since $\topcirc$ is not \emph{homotopically} contractible.
Synthetically, the proof can be done in almost exactly the same way, inserting $\shape$ at the last step, and using the fact that $\Omega\hocirc = \Z$ mentioned in \cref{sec:synth-homot-theory} (which uses the univalence axiom) so that $\shape\topcirc$ (being $\hocirc$) is not a retract of $\shape \topdisc$ (being $\unit$).

At a more advanced level, Schreiber's chapter in this book shows that smooth \oo-groupoids --- and, by extension, cohesive type theory --- are a natural setting for differential cohomology and gauge field theory, which involve the interaction between smooth spaces and homotopy spaces.
The synthetic approach to $\shape$ is thus not just a conceptual way to explain the difference between topological and homotopy spaces, but a practical tool for combining them in applications.

\section{Conclusion}
\label{sec:conclusion}

What does the future hold for type theory and synthetic mathematics?
Current research divides into two threads.
One is ``internal'': developing mathematics in type theory.
This includes both ordinary mathematics in constructive logic without LEM or AC, so as to be valid in all toposes (\cref{sec:constructive-logic}); and also more novel synthetic mathematics using of nonclassical structure (\cref{sec:synthetic-topology-1,sec:homotopy-type-theory,sec:cohes-homot-type}).

The constructivization of ordinary mathematics has a long history, but plenty of fundamental questions remain unanswered, due in part to a tradition among some ``constructivists'' of neglecting propositional truncation and assuming countable choice.
Synthetic mathematics is newer: synthetic differential geometry is several decades old but not well-known outside topos theory, while synthetic homotopy theory is only a handful of years old~\cite{ls:pi1s1,hottbook}, and synthetic topology is in between.
Thus, there are many open questions regarding which results of ``analytic'' mathematics can be reproduced synthetically.



The other thread of current research is ``meta-theoretic''.
As mentioned in \crefrange{sec:semantics}{sec:equality}, there are many unsolved problems in the $\io$-categorical semantics of type theory.
There are also purely syntactic open problems, such as reconciling the topological/homotopical point of view with the computational one from \cref{sec:compute}.
For instance, can we make HITs and univalence ``compute'' (i.e.\ not break the fragile computational interpretation mentioned in \cref{sec:compute})?
(At present the most successful approaches to this use ``cubical'' methods, e.g.~\cite{coquand:ctt0915,cohen:cubicaltt,ahw:chtt-i,ah:chtt-ii}.)

Some problems involve both syntax and semantics.
For instance, homotopy type theory is an excellent synthetic language for higher groupoids, but what about higher \emph{categories}?
Any classical definition of \io-category (such as quasicategories) can be repeated inside the \emph{sets} of type theory, but that would not be what we want: a good definition of \io-category in homotopy type theory should use the synthetic notion of \oo-groupoid provided by the types.
The most promising approach is something like Rezk's ``complete Segal spaces''~\cite{rezk:css}; this can be done for 1-categories~\cite{aks:rezk}, but for the \oo-case it would require a notion of ``coherent simplicial type'', which so far has proven elusive.


This is a special case of another open question that I call the ``problem of infinite objects'', which also applies to other homotopy-theoretic notions like $A_\oo$-spaces and structured ring spectra.
Classically, such infinite coherence structures involve strict point-set-level equalities. 
For instance, $A_\oo$-spaces are \emph{strict} algebras for a topological operad; the weakness is in the operad.
But homotopy type theory, in its most common form, severely restricts the use of strict equality: it can be accessed only using dependent types (e.g.\ terms belonging to a dependent type are strict sections of a fibration) and judgmental equality.
This is good because it makes everything automatically homotopy-invariant, but it means we lack a flexible way to assemble arbitrary higher coherence structures.
(In particular, while synthetic homotopy theory can do a lot, further technical advances are needed before it could reproduce all of classical homotopy theory.)
This problem might be solvable completely internally, but it might also require modifying the syntax, leading to a whole host of new meta-theoretic problems.


Let me end with some remarks about the \emph{philosophical} implications of synthetic mathematics.
I have presented type theory in a way intended to seem useful and unobjectionable to a classical mathematician: as a syntax for reasoning about structured categories in a familiar language.
Crucial to the usefulness of this syntax is the fact that it, like ZFC set theory, is general enough to encode all of mathematics, and therefore anything we can prove (constructively) in ordinary mathematics is automatically also ``true internally'' in any category.

This leads naturally to a slightly different question: can we \emph{actually} use type theory as \emph{the} foundation for mathematics?
That is, must we consider the objects of mathematics to ``really'' be built out of sets, with ``types'' just a convenient fiction for talking about such structures?
Or can we consider \emph{types} to be the basic objects of mathematics, with everything else built out of \emph{them}?

The answer is undoubtedly yes: the ``sets'' in type theory can encode mathematics just like the sets of ZFC can.
Of course, there are subtleties.
On one hand, if our type theory is constructive, 
we need to do our mathematics constructively. 
On another hand, type theory often suggests different ways to do things, using the synthetic spatial or homotopical structure of types instead of analytic topological spaces or $\oo$-groupoids.
\footnote{In particular, for type theory to be an autonomous foundation for mathematics, it ought to suffice for its own metatheory, including the freeness of its own classifying $\io$-category;
but we don't yet even know how to \emph{define} $\io$-categories in homotopy type theory.}
Both of these involve their own open problems; but they are only potential \emph{enhancements} or \emph{refinements} of ordinary mathematics, so regardless of how they turn out, it is certainly \emph{possible} to treat type theory as a foundation for all of mathematics.

The real question, therefore, is not ``can we?''\ but ``should we?''
This is where things get more philosophical.
Over the past century, mathematicians and philosophers have become accustomed to the fundamental objects of mathematics being discrete sets, with no spatial or homotopical structure.
However, \textit{a priori} there is no reason this has to be the case.
Indeed, some of the early-20th-century constructivists, notably Brouwer, can (with a bit of hindsight) be read as arguing for the intrinsically spatial nature of mathematical objects.

But can spaces really be \emph{fundamental} like sets are?
A discrete set certainly seems simpler, and hence more fundamental, than a set equipped with spatial structure.
But this argument merely begs the question, since if spaces are fundamental objects then they are \emph{not} just sets ``equipped with spatial structure''.
In spatial type theory there is no obvious non-tautological ``structure'' with which we can equip the discrete set of reals $\flat\R$ that determines the space of reals $\R$.
Is $\flat\R$ ``simpler'' than $\R$?
When we consider all the pathological nowhere-continuous functions supported by $\flat\R$ but not $\R$, it seems at least consistent to believe that $\R$ is the simpler.
Moreover, discrete sets are just a particular kind of space; so even if they are simpler, that doesn't necessarily argue that non-discrete spaces can't be fundamental.
The empty set $\emptyset$ is probably simpler than $\aleph_\omega$, but in ZFC they are equally fundamental objects (i.e.\ sets).

Similar arguments apply to homotopy spaces, i.e.\ $\oo$-groupoids.
One of the central insights of category theory and homotopy theory is that no class of mathematical objects should be considered without the corresponding notion of isomorphism or equivalence:
we study groups up to isomorphism, spaces up to homeomorphism, categories up to equivalence, and so on.
Thus, all mathematical collections naturally form groupoids, or more generally $\oo$-groupoids, when equipped with the relevant ``notion of sameness''.
(See~\cite{shulman:synhott} for further philosophical discussion of this point.)
The \emph{set}\footnote{or ``proper class''} of all groups is much less tractable, and much less interesting, than the \emph{category} of all groups; so even though the former is ``simpler'' in the sense of containing no nontrivial automorphisms, it is reasonable to regard the latter as being at least as fundamental.

One possible objection to treating spaces as fundamental is to ask how we should decide \emph{which} rules our ``spaces as fundamental'' should satisfy.
Indeed, we have already seen that there are different kinds of synthetic topology adapted for different purposes, modeled respectively by consequential, continuous, or smooth \oo-groupoids.
Moreover, other kinds of synthetic mathematics, such as synthetic domain theory, synthetic differential geometry, and other fields waiting to be developed, will have their own toposes and their own type theories.

However, if we shift perspective a bit, we can see that this is a feature rather than a bug.
Why must we insist on singling out some particular theory as ``the'' foundation of mathematics?
The idea of a ``foundation for mathematics'' stems from the great discovery of 20th century logic that we can encode mathematics into various formal systems and study those systems mathematically.
But in the 21st century, we are sufficiently familiar with this process that we no longer need to tie ourselves to only \emph{one} such system.%
\footnote{In particular, it is meaningless to ask whether statements like the Continuum Hypothesis are ``true''; they are simply true in some systems and false in others.
This perspective is very natural to a category theorist, but has recently made inroads in set theory as well~\cite{hamkins:multiverse}.}
Even ZFC has a role from this point of view: it is a synthetic theory of well-founded membership structures!

Bell~\cite{bell:topos-lst} makes an excellent analogy to Einstein's theory of relativity.
In Newtonian physics, there is a special absolute ``rest frame'', relative to which all motion can be measured.
There are moving observers, of course, but they are second-class citizens: the standard laws of physics do not always apply to them.
They feel ``fictitious forces'', like the centrifugal force and Coriolis force on a spinning merry-go-round or planet, that are not \emph{really} forces but just manifestations of ``truly'' inertial motion in a non-inertial reference frame.

By contrast, Einsteinian physics can be formulated equally well in \emph{any} reference frame and obeys the same laws in each, with consistent rules for transforming between reference frames. 
Some frames, called ``(locally) inertial'', lead to a simpler formulation of the laws; but often this is outweighed by the relevance of some other frame to a particular problem (such as the non-inertial reference frame of the Earth's surface).
The centrifugal and Coriolis forces are exactly as real as any other force; in fact they are simply instances of gravitational force!
To an observer on the Earth's surface, an inertial observer in a spaceship flying by is the one who is spinning (along with the rest of the universe), thereby feeling ``fictitious'' forces that cancel out these gravitational ones.

Similarly, in ZFC orthodoxy there is an absolute notion of ``set'' out of which everything is constructed.
Spaces exist, but they are second-class citizens, ultimately reducible to sets, and the basic axioms of set theory don't apply to them.
But from a pluralistic viewpoint, mathematics can be developed relative to any topos, obeying the same general rules of type theory.
We have consistent rules for translating between toposes along functors, 
and there are some toposes in which mathematics looks a bit simpler (those satisfying LEM or UIP).
However, there is no justification for regarding any particular topos or type theory as the ``one absolute universe of mathematics''.
An observer in a topos of classical mathematics can construct a topos of spaces in which all functions are continuous, thereby explaining its different behavior.
But an observer inside a topos of spaces can also construct a topos of classical mathematics as the ``discrete'' or ``indiscrete'' objects, whose different behavior is explained by the triviality of their cohesion --- and \emph{both points of view are equally valid}.
Just as modern physicists switch reference frames as needed, modern mathematicians should be free to switch foundational systems as appropriate.

This is particularly relevant for physicists and other scientists interested in \emph{using} mathematics rather than debating its Platonic existence.
If a particular synthetic theory is useful in some application domain (see e.g.\ Schreiber's chapter), we are free to take it seriously rather than demanding it be encoded in ZFC.
Set theory and 20th century logic were a crucial stepping-stone to bring us to a point where we can survey the multitude of universes of mathematics; but once there, we see that there is nothing special about the route we took.

\bibliographystyle{plain}
\bibliography{all}

\end{document}